\documentclass[12pt]{article}
\usepackage{geometry}                % See geometry.pdf to learn the layout options. There are lots.
\usepackage{graphicx}
\usepackage{amssymb}
\usepackage{epstopdf}
\begin{document}

\title{{\bf The quartic Fermat equation in Hilbert class fields of imaginary quadratic fields}}         % Enter your title between curly braces
\author{Rodney Lynch and Patrick Morton}        % Enter your name between curly braces
\date{}          % Enter your date or \today between curly braces
\maketitle

\begin{abstract} \noindent It is shown that the quartic Fermat equation $x^4 +y^4=1$ has nontrivial integral solutions in the Hilbert class field $\Sigma$ of any quadratic field $K=\mathbb{Q}(\sqrt{-d})$ whose discriminant satisfies $-d \equiv 1$ (mod 8).  A corollary is that the quartic Fermat equation has no nontrivial solution in $K=\mathbb{Q}(\sqrt{-p})$, for $p$ $( > 7)$ a prime congruent to $7$ (mod 8), but does have a nontrivial solution in the odd degree extension $\Sigma$ of $K$.  These solutions arise from explicit formulas for the points of order 4 on elliptic curves in Tate normal form.  The solutions are studied in detail and the results are applied to prove several properties of the Weber singular moduli introduced by Yui and Zagier.
\end{abstract}

\section{Introduction.}       % Enter section title between curly braces

In the paper [1] Aigner proved that the only quadratic field in which the quartic Fermat equation $x^4+y^4=z^4$ has a nontrivial solution is the field $K=\mathbb{Q}(\sqrt{-7})$, and that all solutions in this field reduce to the one solution

$$\left( \frac{1+\sqrt{-7}}{2} \right)^4 + \left( \frac{1- \sqrt{-7}}{2} \right)^4 = 1.\eqno{(1.1)}$$ \smallskip

\noindent In this paper we will show that this is actually one of an infinite number of solutions of the Fermat quartic in the Hilbert class fields of a certain family of imaginary quadratic fields; namely, those fields in which the rational prime $2$ splits into two prime ideals.  For any number field $L$, let $R_L$ denote the ring of integers in $L$.  \bigskip

\noindent {\bf Theorem 1.1.} Let $K=\mathbb{Q}(\sqrt{-d})$, where $d>0$ and $-d \equiv 1$ (mod 8).  Further, let $(2)=\wp_2 \wp_2'$ be the factorization of $2$ into conjugate prime ideals $\wp_2$ and $\wp_2'$ in $K$.  In the Hilbert class field $\Sigma$ of $K$ there are generators $\pi$ and $\xi$ of the principal ideals $\wp_2R_\Sigma=(\pi)$ and $\wp_2'R_\Sigma=(\xi)$ for which
$$ \pi^4 + \xi^4 = 1.\eqno{(1.2)}$$
The generators $\pi$ and $\xi$ may be chosen to be conjugates over $\mathbb{Q}$, and for $d > 7$ we have $\pm \pi \pm \xi \neq 1$.  Moreover, if $\displaystyle \tau=\left( \frac{\Sigma/K}{\wp_2} \right)$ is the automorphism associated to the ideal $\wp_2$ by the Artin map for $\Sigma/K$, then $\displaystyle \xi=\frac{\pi^{\tau^{2}}+1}{\pi^{\tau^{2}}-1}$. \bigskip

Note that the solution $(x, y, z ) = (\pi, \xi, 1)$ of $x^4 + y^4 =z^4$ in this theorem has the property that no non-zero power of $x$ or $y$ lies in $\mathbb{Q}$.  Since $\pm \pi \pm \xi \neq 1$ for $d > 7$, these solutions do not lie on the intersection of the Fermat curve with the line $x +y =1$, except for the solution (1.1).  In fact, for $d > 7$ the point $(\pi, \xi)$ does not lie on {\it any} rational line.  Thus, these solutions are nontrivial in several senses.  (See [22, p. 400] and [21].)  We note that the numbers $\pi$ and $\xi$ generate $\Sigma$ over $\mathbb{Q}$ (as do $\pi^4$ and $\xi^4$), so their degrees are equal to $2h(-d)$, where $h(-d)$ is the class number of the ring of integers $R_K$ of $K$.  Thus, the degrees of the solutions in Theorem 1.1 satisfy $2h(-d) \ge 4$ for $d>7$.  \medskip

For example, in the Hilbert class field $\Sigma=\mathbb{Q}(\sqrt{-3},\sqrt{5})$ of $K=\mathbb{Q}(\sqrt{-15})$ we have the equation
$$\left( 1- \frac{\sqrt{-3}+\sqrt{5}}{2} \right)^4 + \left( 1 +  \frac{\sqrt{-3}-\sqrt{5}}{2} \right)^4 = 1;$$ \smallskip
\noindent while in the Hilbert class field $\Sigma=\mathbb{Q}(\sqrt{-3},\sqrt{13},\sqrt{-50+14\sqrt{13}})$ of $K=\mathbb{Q}(\sqrt{-39})$ we have $\pi^4 +\xi^4=1$ with

$$\pi=\frac{-3-3\sqrt{-3}+\sqrt{13}+\sqrt{-39}}{4}+\frac{-3+4\sqrt{-3}+\sqrt{-39}}{12}\sqrt{-50+14\sqrt{13}},$$

$$\xi=\frac{-3+3\sqrt{-3}+\sqrt{13}-\sqrt{-39}}{4}+\frac{3+4\sqrt{-3}+\sqrt{-39}}{12}\sqrt{-50+14\sqrt{13}}.$$ \smallskip

\noindent There is also the solution $(\pi, \xi)$ in the Hilbert class field $\Sigma=\mathbb{Q}(\pi)$ of $\mathbb{Q}(\sqrt{-23})$ in which $\pi$ is a root of the irreducible polynomial
$$b_{23}(x)=x^6+x^5+9x^4-13x^3+18x^2-16x+8$$ \smallskip
\noindent and $\xi$ is the conjugate of $\pi$ given as
$$\xi=-\frac{6}{7}+\frac{15}{7}\pi-\frac{5}{4}\pi^2 +\frac{85}{56}\pi^3+\frac{3}{14}\pi^4+\frac{9}{56}\pi^5.$$ \smallskip
\noindent This last solution can also be written in the form

$$\pi=\left(-1-\frac{7}{23}\sqrt{-23} \right) \gamma^2+\left(1-\frac{\sqrt{-23}}{23}\right) \gamma+\frac{1}{2}+\frac{17}{46}\sqrt{-23},$$
$$\xi=\left(2-\frac{2}{23}\sqrt{-23} \right) \gamma^2+\left(-1+\frac{3}{23}\sqrt{-23} \right) \gamma-\frac{3}{2}-\frac{5}{46}\sqrt{-23},$$ \smallskip

\noindent where $\gamma$ is a root of the polynomial $x^3-x-1$.  This solution is especially interesting, in that (1.2) has no nontrivial solution in $K=\mathbb{Q}(\sqrt{-23})$ but does have a nontrivial solution in the cubic extension $\Sigma$ of $K$.   Theorem 1.1 provides infinitely many examples of this type, since when $d=p \equiv 7$ (mod 8) is a prime, the class number $h(-d)=[\Sigma : K]$ is odd.  \bigskip

\noindent {\bf Theorem 1.2.} If $p \equiv 7$ (mod 8) is a prime $> 7$, then $x^4+y^4=z^4$ has no nontrivial solution in the field $K=\mathbb{Q}(\sqrt{-p})$, but does have a nontrivial solution in the odd degree extension $\Sigma$ of $K$, where $\Sigma$ is the Hilbert class field of $K$. \bigskip

This is reminiscent of the examples given in [3], except that there quartic diophantine equations are considered which have {\it no} solution in $\mathbb{Q}$ but do have solutions in cubic extensions of $\mathbb{Q}$.  \bigskip

In the proof we consider the points of order 4 on the elliptic curve
$$E_1: \hspace{.1 in} Y^2+XY+ bY = X^3 + bX^2,$$ \smallskip
\noindent the Tate normal form for a curve with a point of order 4, where $b=1/\alpha^4$.  We find expressions for all the points in $E_1[4]$ by introducing the diophantine condition
$$ Fer_4: \quad 16 \alpha^4+16 \beta^4=\alpha^4 \beta^4.\eqno{(1.3)}$$ \smallskip
\noindent We show that the non-zero coordinates of points in $E_1[4]$ are in the multiplicative group generated by $-1, 2$ and the linear fractional quantities $\displaystyle \beta_n=\frac{\beta + 2i^n}{2\beta}$, ($1\le n \le 4$, $ i=\sqrt{-1}$) (see Proposition 4.1 and [26]).  The expressions we give are analogous to the nice form for the points of order 4 on the Jacobi normal form, but are more convenient for our purposes.  The formulas in Proposition 4.1 are also the analogue of the formulas for the points of order $3$ on the Deuring normal form given in [28], which are expressed in terms of solutions of the Fermat cubic $27X^3+27Y^3=X^3Y^3$.  See [26] for similar expressions for the points in $E_1[8]$.  \medskip

With a judicious choice of $\alpha$ in terms of the Dedekind $\eta$-function, namely

$$\alpha= \zeta_8^{-1} \frac{\eta(w/4)^2}{\eta(w)^2},\hspace{.2 in} \zeta_8 = e^{2\pi i/8},$$ \smallskip

\noindent for a certain element $w=(v+\sqrt{-d})/2$ of $K$, the curve $E_1$ has complex multiplication by the maximal order $R_K$ in the field $K=\mathbb{Q}(\sqrt{-d})$.  It is known [32, p. 159] that $\alpha^4 \in \Sigma$, the Hilbert class field of $K$; it follows that the $X$-coordinates of points of order $4$ on $E_1$ lie in the ray class field $\Sigma_4=\Sigma(i)$ over $K$ of conductor $f=4$.  Using our formulas for the points of order $4$ on $E_1$ and similar formulas for the corresponding points on a second elliptic curve $E_3$ (see Sections 2 and 5) we show that $\beta \in \Sigma$.  This leads to a solution of (1.2) with
$$\pi = \frac{\beta}{\zeta_8^j \alpha}, \quad \xi=\frac{\beta}{2}, \quad \textrm{for} \ \textrm{some} \ j \in \{1, 3, 5, 7\}.$$  \smallskip
\noindent In Section 8 we normalize $\zeta_8^j \alpha, \beta$ by choosing $j$ and the sign of $\beta$ so that $\xi=\beta/2$ is conjugate to $\pi$ and to $(\xi+1)/(\xi-1)$ over $\mathbb{Q}$.  Then we prove algebraically that $\alpha_1=\zeta_8^j \alpha$ and the unit
$$\gamma = \frac{\beta(\beta+2)}{4(\beta-2)}$$
are squares in $\Sigma$, hence that $\zeta_8^{(j-1)/2} \eta(w/4)/\eta(w) \in \Sigma$.  This also implies that the elliptic curve $Y^2=X(X^2-4)$ has the integral solution $P=(\beta, 2(\beta-2) \sqrt{\gamma})$ in $\Sigma$.  Whenever the discriminant $-d=-p$ is prime, a simple argument shows that the sum of the conjugates of the point $P$ for $\Sigma/K$ is a point of infinite order on $Y^2=X(X^2-4)$, defined over $K$ (see Theorem 8.11).  (The curves $E_1$ and $E_3$ in the above argument are isogenous, and for the given choices of $d$ and $\alpha$, the isogeny $(E_1 \rightarrow E_3)$ represents a {\it Heegner point} on the modular curve $X_0(4)$.) \medskip

The same arguments show that for every odd, natural number $f$ there is a solution $(\pi_f,\xi_f)$ of (1.2) in the ring class field $\Omega_f$ of $K$ with conductor $f$, and that $\Omega_f=\mathbb{Q}(\pi_f)=\mathbb{Q}(\xi_f)$.  Thus, there is a whole lattice of solutions to the Fermat quartic in ring class fields of $K$ with odd conductor.  Furthermore, if $f | f'$, the numbers $\pi_{f'}/\pi_f$ and $\xi_{f'}/\xi_f$ are units in $\Omega_{f'}$.  (See Section 6 and Remark 4 in Section 12.)  These solutions reflect the existence of the associated Heegner points of $X_0(4)$ defined over the ring class fields $\Omega_f$. \medskip

In Section 10 we give explicit expressions for the solution $(\pi_f, \xi_f)$ in terms of the Schl\"afli modular functions used by Weber (see [5, p. 256] and [31], [32]).  We use a result of [36] to give an explicit value of $j \in \{1, 3, 5, 7\}$ in Theorem 10.6 for which $\displaystyle \zeta_8^{(j-1)/2} \frac{\eta(w/4)}{\eta(w)}$ lies in the ring class field $\Omega_f$, for $\displaystyle w=\frac{v+\sqrt{-d}}{2}$ with $v= 1$ or $3$ and $v^2 \equiv -d$ (mod 16).   This is similar to, but not contained in, results of [17] and [31] or [32, Chs. 5-6], because the above $\eta$-quotient is evaluated at the ideal $\wp_2^2$, in the notation of [17].  The fact that this ideal is not relatively prime to $2$ necessitates the alternate approach taken here.  \medskip

The numbers $\pi$ and $\xi$ in Theorem 1.1 turn out to be fourth roots of numbers $\lambda$ for which the Legendre normal form with parameter $\lambda$ has complex multiplication by $R_K$.  Using (1.3), this observation leads to the following result.  \bigskip

\noindent {\bf Theorem 1.3.} Let $\lambda$ be a complex number for which the Legendre normal form
$$E_\lambda: \quad Y^2 = X(X-1)(X-\lambda)$$ \smallskip
\noindent  has complex multiplication by the order $\textsf{R}_{-f^2d}$ of $K=\mathbb{Q}(\sqrt{-d})$ with odd conductor $f$ and discriminant $-f^2d \equiv 1$ (mod 8).  Then $\lambda$ lies in the ring class field $\Omega_f$ of $K$ with conductor $f$.  If either $\lambda$ or $1/\lambda$ is an algebraic integer, then $\lambda$ is a fourth power in the field $\Omega_f$.  If neither $\lambda$ nor $1/\lambda$ is an algebraic integer, then $-\lambda$ is a fourth power in $\Omega_f$.  The numbers $\lambda$ and $\lambda-1$ are only divisible by prime ideal divisors of 2 in $\Omega_f$. \medskip

This theorem is related to Theorem 1.2 of [27], which says that the values of $\bar \lambda$ for which the Legendre normal form $E_{\bar \lambda}$ is supersingular in characteristic $p$ are fourth powers in the finite field $\mathbb{F}_{p^2}$.  (See also Landweber [24].)  Theorem 1.3 is an analogue of this result in characteristic 0.  Since some square-root of the form $\sqrt{-d}$ with $d \equiv 7$ (mod 8) embeds into the endomorphism ring of any supersingular curve in characteristic $p>2$, the theorem of [27] follows from Theorem 1.3 by reduction (see Theorems 7.1 and 9.1). \medskip

A secondary purpose of this paper is to relate our results to the paper [36] of Yui and Zagier, in which they give generators of small height for ring class fields $\Omega_f$ of the quadratic fields $K=\mathbb{Q}(\sqrt{-d})$, with $-d=-f^2d_1 \equiv 1$ (mod $8$) and $(d,3)=1$.  Yui and Zagier define the Weber singular modulus $f_{w} (Q)$, for any primitive quadratic form $Q(x,y)=ax^2+bxy+cy^2$ of discriminant $-d$, as one of the Schl\"afli modular functions $\mathfrak{f}, \mathfrak{f}_1, \mathfrak{f}_2$ (depending on the parity of $a$ and $c$) evaluated at the root of $Q(x,1)=0$ which lies in the upper half-plane, times a well-chosen $48$-th root of unity.  They show that $f_w(Q)=f_w(\mathcal{A})$ only depends on the $SL_2(\mathbb{Z})$ equivalence class of the form $Q$, or in other words, only on the ideal class $\mathcal{A}$ in the order $\textsf{R}_{-d}$ to which the form $Q$ corresponds.  In Section 10 we prove their conjecture [36, p. 1648] that the numbers $f_w(\mathcal{A})$ are conjugate over $\mathbb{Q}$, as $\mathcal{A}$ ranges over the ideal classes in $\textsf{R}_{-d}$, and use this fact to prove the following. \bigskip

\noindent {\bf Theorem 1.4.} Let $\mathcal{T}$ be the ideal class of $\textsf{R}_{-d}$ containing the ideal $\wp_2 \cap \textsf{R}_{-d}=2\mathbb{Z}+ \left(\frac{v+\sqrt{-d}}{2} \right)\mathbb{Z}$, with $v^2+d \equiv 0$ (mod $16$), and assume $(d,3)=1$.  If $\mathcal{A}$ represents any ideal class in the order $\textsf{R}_{-d}$, the number $2f_w(\mathcal{A})/f_w(\mathcal{TA})^2$ is the $Y$-coordinate of an integral point $(X,Y)$ on the curve
$$Y^4=X(X^2+4)$$
with coordinates in $\Omega_f$.  The corresponding $X$-coordinate is the negative conjugate $X=-\beta_c=-\beta^{\tau_c^{-1}}$ of $\beta$, where $\tau_c=(\Omega_f/K,R_K \mathfrak{c})$ with $\mathfrak{c} \in \mathcal{A}$.  \bigskip

We relate the numbers $\alpha_1$ and $\beta$ that we study to the Weber singular moduli of Yui and Zagier by showing that the number $4\alpha_1/\beta^2=4\zeta_8^j \alpha/\beta^2$ is the $6$-th power of the unit $f_w(Q _2)$ considered by Yui and Zagier, where $Q_2(x,y)=2x^2-vxy+\left(\frac{v^2+d}{8}\right) y^2$, whether or not $3$ divides $d$.  Furthermore, we show that if $3 \mid d$ and $Q_1(x,y)=x^2-vxy+\left(\frac{v^2+d}{4}\right) y^2$, then the product $f_w(Q_1)f_w(Q_2)$ of the Weber singular moduli for the forms $Q_1$ and $Q_2$ is a generator of the ring class field $\Omega_f$ over $K$ (for $d>15$), although neither of the numbers $f_w(Q_1), f_w(Q_2)$ lies in $\Omega_f$.  Examples show that this number has small height and relatively small discriminant, as is the case for the numbers $f_w(Q_1)$ when $(d,3)=1$.  This extends the discussion in the paper [36] for the case $3 \mid d$ and gives a systematic procedure for producing generators with small height and discriminant in this case.  We give a number of examples of the minimal polynomials of these generators in Sections 11 and 12.  \medskip

In our discussion we rely mostly on basic facts from class field theory and the theory of complex multiplication, as well as a few facts about the Dedekind eta-function and the Schl\"afli functions, which can be found in Weber's treatise [35].  (See also [32], [33], and [36].)  In particular, we do not use any of the advanced theory of Heegner points in our discussion.  We take a mostly algebraic approach in Sections 3-9; this is advantageous when we turn to analytic representations of our solutions in Sections 10 and 11.  This combined approach leads to the proof of Theorem 1.4 (see Theorems 10.4 and 10.5) and to the proof that the numbers $f_w(Q_1)f_w(Q_2)$ above have degree $h(-d)$ over $\mathbb{Q}$ in Theorem 11.1.

\section{Weber's functions and the $j$-invariant.}

\noindent In this section we will use classical functions considered by Weber [35] to recognize specific values of the $j$-function as $j$-invariants of certain elliptic curves in Tate normal form.  \medskip

We start with Weber's functions

$$x=\mathfrak{f}_2(w)^8=16\frac{\eta(2w)^8}{\eta(w)^8}, \hspace{.2 in} x_0=\mathfrak{f}_1(w)^8=\frac{\eta(w/2)^8}{\eta(w)^8},$$
$$x_1=-\mathfrak{f}(w)^8= e^{2 \pi i/3}\frac{\eta((w+1)/2)^8}{\eta(w)^8}=\frac{\eta((w+3)/2)^8}{\eta(w)^8},$$ \smallskip

\noindent where $\eta(z)$ is the Dedekind eta-function.  (See [35, p. 250], [30], [31], and and Section 10 for the definitions of the Schl\"afli functions $\mathfrak{f}, \mathfrak{f}_1, \mathfrak{f}_2$.)  These functions satisfy the equation
$$X^3-X \gamma_2(w)+16=0, \hspace{.2 in} \gamma_2(w) = j(w)^{1/3}.\eqno{(2.1)}$$ \smallskip
\noindent Thus, we have
$$j(w)=\gamma_2(w)^3=\bigg( \frac{x^3+16}{x} \bigg)^3.\eqno{(2.2)}$$ \smallskip
\noindent Replacing $w$ by $w/2$ transforms the function $x$ into $16/x_0$; it follows that 

$$16^2 -x_0^2 \gamma_2\bigg(\frac{w}{2}\bigg) + x_0^3 = 0,\eqno{(2.3)}$$
and therefore
$$j\bigg(\frac{w}{2}\bigg) = \gamma_2\bigg(\frac{w}{2}\bigg)^3 = \frac{(x_0^3+16^2)^3}{x_0^6}.\eqno{(2.4)}$$ \smallskip
Now set 
$$y(w) = \frac{\eta(w/4)^4}{\eta(w)^4}=\mathfrak{f}_1(w/2)^4 \mathfrak{f}_1(w)^4,$$ \smallskip
so that $y^2=x_0(w/2)x_0(w)$.  Replace $w$ once again by $w/2$ in (2.1) for the function $x_0$ in place of $x$ and use $x_0^2\gamma_2(w/2)=16^2+x_0^3$ from (2.3).  This gives

$$0=y^6-y^2x_0^2\gamma_2\bigg(\frac{w}{2}\bigg)+16x_0^3=y^6-y^2(16^2+x_0^3)+16x_0^3$$
$$\hspace{2.1 in}=(y^2-16)(y^4+16y^2-x_0^3).$$ \smallskip

\noindent It follows that $x_0^3=y^4+16y^2$; putting this into (2.2) (with $x_0$ for $x$) and (2.4) gives

$$j(w)=\frac{(y^4+16y^2+16)^3}{y^4+16y^2},\hspace{.2 in} j\bigg(\frac{w}{2}\bigg)=\frac{(y^4+16y^2+256)^3}{(y^4+16y^2)^2}.\eqno{(2.5)}$$ \smallskip

\noindent Putting $w/2$ for $w$ in (2.3) leads similarly to the equation

$$16^2x_0^3-y^4x_0\gamma_2\bigg(\frac{w}{4}\bigg)+y^6=0.$$ \smallskip

\noindent Rearranging and cubing gives

$$16^2(y^4+16y^2)+y^6=y^4x_0\gamma_2\bigg(\frac{w}{4}\bigg)$$
$$(y^6+16^2y^4+16^3y^2)^3=y^{12}(y^4+16y^2)j\bigg(\frac{w}{4}\bigg).$$
Hence,
$$j\bigg(\frac{w}{4}\bigg)=\frac{(y^4+256y^2+4096)^3}{y^8(y^2+16)}.\eqno{(2.6)}$$ \smallskip

\noindent Thus, $j(w/4)$ arises from $j(w)$ by the substitution $y \rightarrow 16/y$. \medskip

Now we set

$$\alpha=\zeta_8^{-1} \sqrt{y(w)} = \zeta_8^{-1} \frac{\eta(w/4)^2}{\eta(w)^2},\hspace{.2 in} \zeta_8 = e^{2\pi i/8},\eqno{(2.7)}$$ \smallskip

\noindent so that $\alpha^2 = -iy(w)$.  From (2.5) and (2.6) we obtain

$$j(w) = \frac{(\alpha^8-16\alpha^4+16)^3}{\alpha^8-16\alpha^4},\hspace{.2 in} j\bigg(\frac{w}{2}\bigg)=\frac{(\alpha^8-16\alpha^4+256)^3}{\alpha^8(\alpha^4-16)^2},\eqno{(2.8)}$$ \smallskip

$$ j\bigg(\frac{w}{4}\bigg)=\frac{(\alpha^8-256\alpha^4+4096)^3}{\alpha^{16}(16-\alpha^4)}.\eqno{(2.9)}$$ \smallskip

\noindent From the calculations of [27, pp.253-254] we see that $j(w)$ is the $j$-invariant of the elliptic curve

$$E_1(\alpha): \hspace{.1 in} Y^2+XY+\frac{1}{\alpha^4} Y = X^3 + \frac{1}{\alpha^4} X^2,\eqno{(2.10)}$$ \smallskip

\noindent which is the Tate normal form for a curve with a point of order $n=4$; and $j(w/2)$ is the $j$-invariant of the elliptic curve

$$E_2(\alpha): \hspace{.1 in} Y^2 + XY + \frac{2}{\alpha^4} Y = X^3 + \frac{4}{\alpha^4} X^2 - \frac{1}{\alpha^8},\eqno{(2.11)}$$ \smallskip

\noindent which is 2-isogenous to $E_1(\alpha)$ by the map $\psi=(\psi_1,\psi_2):E_1(\alpha) \rightarrow E_2(\alpha)$ with

$$\psi_1(X,Y)=\frac{X^2}{X+b}, \quad \psi_2(X,Y)=\frac{-b^2}{X+b}+\frac{X(X+2b)Y}{(X+b)^2}, \quad b= \frac{1}{\alpha^4}.$$ \smallskip

\noindent  Still with $b = 1/\alpha^4$, the isogeny

$$X_1 = \phi_1(X,Y) = \frac{X^2-b}{X+4b}, \hspace{.2 in} $$

$$Y_1=\phi_2(X,Y) = \frac{bX^2+(b-8b^2)X+3b^2-32b^3}{(X+4b)^2}+\frac{X^2+8bX+b}{(X+4b)^2} Y$$ \smallskip

\noindent maps $E_2(\alpha)$ to the curve

$$E_3(\alpha): \hspace{.1 in} Y_1^2 + X_1Y_1 + \frac{4}{\alpha^4}Y_1 = X_1^3 + \frac{16}{\alpha^4}X_1^2 +\frac{6}{\alpha^4}X_1 + \frac{\alpha^4-4}{\alpha^8},\eqno{(2.12)}$$ \smallskip

\noindent and the $j$-invariant of this curve is the expression $j(w/4)$ in (2.9).  It is not difficult to verify that the isogeny $(E_1(\alpha) \rightarrow E_3(\alpha))$ is cyclic, and therefore represents a point on the modular curve $X_0(4)$.

\section{Determining the field of definition.}

Now let $w = (v+\sqrt{-d})/2$, where $-d \equiv 1$ (mod 8) is the discriminant of the quadratic field $K=\mathbb{Q}(\sqrt{-d})$ and $v$ is an odd integer satisfying $v^2 \equiv -d$ (mod 16).  Then $(2,w)=\wp_2$ is a prime divisor of $2$ in the ring of integers $R_K$ of $K$ with the basis $\{2, w\}$ and $\wp_2^2$ has the basis $\{4,w\}$.  It follows that $j =j(w)$,  $k=j(w/2)$ and $l=j(w/4)$ are roots of the class equation $H_{-d}(x)=0$ and elements of the Hilbert class field $\Sigma=\mathbb{Q}(\sqrt{-d},j(w))$ of $K$. (For the application in Sections 10 and 11, we note that we could just as well take $w=(v+\sqrt{-d})/(2a)$, where $16a | (v^2 + d)$, in which case $\{a, aw\}$ is a basis for an ideal $\mathfrak{a}$ and $\{4a, aw\}$ is a basis for $\wp_2^2 \mathfrak{a}$.) \medskip

From (2.8) and (2.9) we see that $t=\alpha^4$ is a root of three sixth degree equations over $\Sigma$ :

$$f(x)=(x^2-16x+16)^3-j(x^2-16x), $$
$$g(x)=(x^2-16x+256)^3-k(x^2-16x)^2.$$ 
$$h(x) = (x^2-256x+4096)^3-lx^4(16-x).$$\smallskip

\noindent We would like to show that the ideal $I=(f(x),g(x),h(x))$ in $R_\Sigma[x]$ contains a non-constant linear polynomial.  We perform the first two steps of the Euclidean algorithm on $f(x)$ and $g(x)$.  First we set
$$q_1(x) = g(x)-f(x) = (720-k)x^4+(-23040+32k)x^3+(j+380160-256k)x^2$$
$$+(-3133440-16j)x+16773120.$$ \smallskip
\noindent Next, take the remainder of $f(x)$ with respect to $q_1(x))$:
$$(k-720)^2 f(x) - \bigg((720-k)x^2+(-11520+16k)x-48k-j-161280\bigg) q_1(x) = $$
$$q_2(x) = (768k^2-k^2 j+25067520k+1488kj-161280j+19906560000+j^2)x^2$$
$$+(-12288k^2+16k^2 j-401080320k-23808kj+2580480j-318504960000-16j^2)x$$
$$+4096k^2+799211520k+2707292160000+16773120j.$$ \smallskip
Setting
$$a_2(j,k)=768k^2-k^2 j+25067520k+1488kj-161280j+19906560000+j^2,$$ \smallskip

\noindent we calculate the remainder of

$$q_3(x)=h(x)-g(x)=(-720+l)x^5+(k-16l+207360)x^4+(-32k-23040000)x^3$$
$$\hspace{1.5 in}+(256k+855244800)x^2-12881756160x+68702699520$$ \smallskip

\noindent with respect to $q_2(x)$ and find

$$a_2(j,k)^2 q_3(x) - u(x) q_2(x) = a_3(j,k,l)x+b_3(j,k,l),$$ \smallskip

\noindent where $u(x)$ is cubic and

$$a_3(j,k,l)=-2^{20}(9331200000-226800j+j^2+21945600k+1488kj+752k^2-k^2j)$$
$$\times (4095j^2+6093360kj-663390000j-4095k^2j+4095jl+102511008000k$$
$$+660960000l+3144240k^2+195120kl+k^2l+81041472000000)$$

$\hspace{.5 in}=-2^{20}A_1(j,k)A_2(j,k,l).$ \bigskip

\noindent We take the resultant of the first factor of $a_3(j,k,l)$ and the modular equation

$$\Phi_2(j,k)= k^3+1488k^2j-k^2j^2-162000k^2+8748000000k+1488j^2k+40773375kj$$
$$+j^3-162000j^2+8748000000j-157464000000000$$ \smallskip

\noindent with respect to $k$ and find

$$\textrm{Res}_k(A_1(j,k),\Phi_2(j,k))=-2^{32}(j-54000)(j+3375)^2(j^2+191025j-121287375)^2$$
$$=-2^{32}H_{-12}(j)H_{-7}(j)^2H_{-15}(j)^2.$$ \smallskip

\noindent Next we take the resultant of the second factor of $a_3(j,k,l)$ and $\Phi_2(k,l)$ with respect to $l$ and then take the resultant of that polynomial and $\Phi_2(j,k)$ with respect to $k$: \bigskip

$\textrm{Res}_k(\Phi_2(j,k),\textrm{Res}_l(\Phi_2(k,l),A_2(j,k,l)))=$
$$-2^{85}(-574852565088000000-1645137216000j-4519800j^2+2047j^3)$$
$$\times(j-16581375)^2(j^2-37018076625j+153173312762625)^2$$
$$\times(j+3375)^6(j^2+191025j-121287375)^6$$ \smallskip

$\hspace{1 in} =-2^{85}A_3(j)H_{-28}(j)^2 H_{-60}(j)^2 H_{-7}(j)^6 H_{-15}(j)^6.$ \bigskip

\noindent Now the roots of the cubic $A_3(j)$ are not algebraic integers, and $j=j(w)$ is only a root of the class equation $H_{-d}(x)$.  Hence, if $d \neq 7, 15$, the above resultants cannot be zero, so that $a_3(j,k,l) \neq 0$.  Thus the ideal $I$ contains a non-constant linear polynomial in $R_\Sigma[x]$, and it follows, since $t=\alpha^4$ is a root of each of the polynomials $f,g,h$, that $\alpha^4 \in \Sigma$.  We summarize this in \bigskip

\noindent {\bf Proposition 3.1.} Let $-d \equiv 1$ (mod 8) be the discriminant of the quadratic field $K=\mathbb{Q}(\sqrt{-d})$.  If
$$w = \frac{v+\sqrt{-d}}{2}, \hspace{.05 in} \textrm{with} \hspace{.05 in} v^2 \equiv -d  \hspace{.05 in} (\textrm{mod} \hspace{.05 in} 16),$$ \smallskip
\noindent then the number $\displaystyle -\alpha^4= \frac{\eta(w/4)^8}{\eta(w)^8}$ lies in the Hilbert class field $\Sigma$ of $K$. \bigskip

In the case $d = 7$ we can use (2.8) and (2.9) along with $H_{-7}(x)=x+3375$ to compute directly that

$$t = \alpha^4 = \frac{31\pm 3 \sqrt{-7}}{2}.$$ \smallskip

\noindent Similarly, when $d=15$ we can use $H_{-15}(x)=x^2+191025x-121287375$ to compute that $t = \alpha^4$ is a root of the polynomial

$$r_{15}(t)=t^4-17t^3+33t^2-4352t+65536.$$ \smallskip

\noindent This polynomial has degree $4 = 2h(-15)$, so it is easy to see that its roots lie in the Hilbert class field of $\mathbb{Q}(\sqrt{-15})$.  This verifies the proposition in the two remaining cases $d=7, 15$. \bigskip

Proposition 3.1 also follows from Theorem 6.6.4 in  Schertz's monograph [32, p. 159].  As a corollary we have \bigskip

\noindent {\bf Proposition 3.2.} With $d$ and $\alpha$ as in Proposition 3.1, the curves $E_1(\alpha), E_2(\alpha)$, and $E_3(\alpha)$ in (2.10)-(2.12) have complex multiplication by the ring of integers $R_K$ of $K=\mathbb{Q}(\sqrt{-d})$ and are defined over the Hilbert class field $\Sigma$ of $K$.  \bigskip

If $-d=-f^2 d_1$ with odd $f>1$ and $d_1$  square-free, then there is an analogous result for the order $\textsf{R}_{-d}$ of discriminant $-d$ in $K$, since the ideal $\wp_{2,-d}=\wp_2 \cap \textsf{R}_{-d}$ in $\textsf{R}_{-d}$ of elements divisible by $\wp_2=(2,w)$  has the $\mathbb{Z}$-basis $\{2,(v+\sqrt{-d})/2\}$ in $\textsf{R}_{-d}$ (cf. [10, pp. 19-20]).  With $\alpha$ defined as in (2.7), for $w=(v+\sqrt{-d})/2$ and $v^2 \equiv -d$ (mod 16), all the same formulas hold, and the numbers $j(w), j(w/2), j(w/4)$ are roots of $H_{-d}(x)=0$.  Each of these numbers generates the ring class field $\Omega_f$ of conductor $f$ over $K$, and exactly as before we obtain:  \bigskip

\noindent {\bf Proposition 3.3.} If $-d=-f^2d_1$ with odd $f$ and square-free $d_1 \equiv 7$ (mod 8), $\alpha$ as in (2.7), and $w=(v+\sqrt{-d})/2$ as above, the curves $E_1(\alpha), E_2(\alpha)$, and $E_3(\alpha)$ in (2.10)-(2.12) have complex multiplication by the order $\textsf{R}_{-d}$ of $K=\mathbb{Q}(\sqrt{-d_1})$ and are defined over the ring class field $\Omega_f$ of conductor $f$ over $K$.  \bigskip

This proposition shows that for the given $d$ and $\alpha$, the isogenies $(E_1(\alpha) \rightarrow E_3(\alpha))$ are examples of {\it Heegner points} on the curve $X_0(4)$.  See [2] and [16].

\section{Points of order 2 and 4 on $E_1$.}

\noindent For the application in Section 6 we derive formulas for all the points of order 4 on the curve $E_1(\alpha)$ in terms of the number $\beta$ defined by the diophantine condition that $(X, Y)=(\alpha, \beta)$ lies on the curve
$$Fer_4: \hspace{.1 in} 16X^4+16Y^4 = X^4Y^4.$$ \smallskip
We let $E_1=E_1(\alpha)$ and set $b=1/\alpha^4$, so that

$$E_1: \hspace{.1 in} Y^2+XY+ b Y = X^3 +b X^2, \quad b = \frac{1}{\alpha^4}.\eqno{(4.1)}$$ \smallskip

\noindent  The doubling formula for the $X$-coordinate of a point $P=(x,y)$ on $E_1$ is

$$X(2P)=\frac{x^4-bx^2-2b^2x-b^3}{(x+b)(4x^2+x+b)}, \quad x=X(P).$$ \smallskip

\noindent Iterating this formula, we find that the denominator of $X(4P)$ is the polynomial

$$(x+b)(4x^2+x+b) x^2 (x+2b)^2 (2x^4 + x^3 + 3bx^2 +4b^2 x+2b^3)^2.$$ \smallskip

\noindent The roots of the first two factors are the $X$-coordinates of the points of order 2, which are given by $X=-b = -1/\alpha^4$ and

$$X=\frac{-1}{8} \pm \frac{1}{8} \sqrt{1-16b} = \frac{-1}{8} \pm \frac{1}{8} \sqrt{\frac{\alpha^4-16}{\alpha^4}}= \frac{-1}{8} \pm \frac{1}{8} \frac{4}{\beta^2}=\frac{-\beta^2 \pm 4}{8\beta^2}, \eqno{(4.2)}$$  \smallskip

\noindent on using the defining condition $16\alpha^4 + 16\beta^4 = \alpha^4 \beta^4$ for $\beta$. \medskip

The $X$-coordinates of the points of order 4 on $E_1$ are solutions of the equation

$$x (x+2b) (2x^4 + x^3 + 3bx^2 +4b^2 x+2b^3)=0.$$ \smallskip

\noindent The first factor gives the points $(0,0)$ and $(0,-b)$.  The second factor gives the points

$$x = -2b, \quad y = \frac{b}{2}(1 \pm  \sqrt{1-16b}).$$ \smallskip

\noindent Using $16\alpha^4 + 16 \beta^4 = \alpha^4 \beta^4$ we can write the $y$-coordinates of these points in the form

$$y = \frac{b}{2}(1 \pm  \sqrt{1-16b}) = \frac{b}{2}\left(1 \pm  \sqrt{\frac{\alpha^4-16}{\alpha^4}}\right) = \frac{b}{2}\left(1 \pm \frac{4}{\beta^2}\right).\eqno{(4.3)}$$ \smallskip

\noindent Finding the remaining eight points of order 4 requires solving the quartic equation $2x^4+x^3+3bx^2+4b^2x+2b^3=0$, which we do by introducing the indeterminate $z$ and completing the square:

$$ \left(x^2+{\frac{1}{4}}x+z\right)^2={\frac{1}{16}}(32z+1-24b)x^2+{\frac{1}{2}}(z-4b^2)x+z^2-b^3. $$ \smallskip

\noindent We choose $z=\frac{b}{4}$ to make the right hand side of the last equation a square, and find that

$$\left(x^2+{\frac{1}{4}}x+{\frac{b}{4}}\right)^2=\left({\frac{1}{4}}\sqrt{1-16b}\, (x+b)\right)^2.$$ \smallskip

\noindent This leads to the quadratic equations
$$x^2+{\frac{1}{4}}(1\pm\sqrt{1-16b}\,)x+{\frac{b}{4}}(1\pm\sqrt{1-16b}\,)=0,$$ \smallskip
\noindent whose solutions are

$$x=-{\frac{1}{8}}\left(1-\sqrt{1-16b}\,\pm\sqrt{2}\sqrt{1-16b-(1-8b)\sqrt{1-16b}}\,\right),$$
$$-{\frac{1}{8}}\left(1+\sqrt{1-16b}\,\pm\sqrt{2}\sqrt{1-16b+(1-8b)\sqrt{1-16b}}\,\right). \eqno{(4.4)}$$ \smallskip

\noindent The corresponding values of $y$ can be found by solving $y^2+xy+by=x^3+bx^2$ for $y$:

$$y=-{\frac{1}{2}} \left(x+b\pm\sqrt{(x+b)(4x^2+x+b)} \right). \eqno{(4.5)}$$ \smallskip

We shall write these solutions in terms of the solution $(\alpha,\beta)$ of $Fer_4$.  Note that

$$ 
b=\frac{1}{\alpha^4}={\frac{\beta^4-16}{16\beta^4}}
=\left({\frac{\beta+2i}{2\beta}}\right)
\left({\frac{\beta-2}{2\beta}}\right)
\left({\frac{\beta-2i}{2\beta}}\right)
\left({\frac{\beta+2}{2\beta}}\right)=\beta_1\beta_2\beta_3\beta_4
$$
\smallskip

\noindent where $i=\sqrt{-1}$ and

$$\beta_n={\frac{\beta+2i^n}{2\beta}},\quad 1\le n\le 4.\eqno{(4.6)}$$ \smallskip
\noindent Next, note that
$$ 1-8b={\frac{\beta^4+16}{2\beta^4}},\quad 1-16b={ \frac{16}{\beta^4}}.\eqno{(4.7)}$$ \smallskip

\noindent We will take $\displaystyle \sqrt{1-16b}={\frac{4}{\beta^2}}$, since in all the formulas below this expression appears with both plus and minus signs.  
For the points of order 2 we have from (4.2):

$$X= \frac{-\beta^2+4}{8 \beta^2}=-\frac{1}{2} \left( \frac{\beta-2}{2 \beta} \right) \left( \frac{\beta+2}{2 \beta} \right)=-\frac{1}{2} \beta_2\beta_4,$$ \smallskip

\noindent and the corresponding $Y$-coordinate is
\begin{eqnarray*}
Y={\frac{1}{16}}(1-8b-\sqrt{1-16b}\,)={\frac{1}{16}}\left({\frac{\beta^4+16}{2\beta^4}}-{\frac{4}{\beta^2}}\right)={\frac{1}{2}}\left({\frac{\beta^2-4}{4\beta^2}}\right)^2\\
=\frac{1}{2} \left( \frac{\beta-2}{2\beta}\right)^2 \left(\frac{\beta+2}{2\beta} \right)^2=\frac{1}{2} \beta_2^2\beta_4^2. \hspace{.2 in}
\end{eqnarray*}

\noindent Similarly, taking the minus sign in (4.2) leads to the coordinates

$$X={\frac{1}{8}}(-1-\sqrt{1-16b}\,)=-{\frac{1}{2}}\beta_1\beta_3,\quad
Y={\frac{1}{16}}(1-8b+\sqrt{1-16b}\,)=  {\frac{1}{2}}\beta_1^2\beta_3^2.$$ \smallskip
\noindent Thus the points of order 2 on $E_1$ are given by
$$(-\beta_1\beta_2\beta_3\beta_4,0), \left(-\frac{1}{2} \beta_1\beta_3,\frac{1}{2} \beta_1^2\beta_3^2\right) ,  \left(-\frac{1}{2} \beta_2\beta_4,\frac{1}{2} \beta_2^2\beta_4^2\right) .\eqno{(4.8)}$$ \smallskip

\noindent For the $Y$-coordinates of the points of order 4 in (4.3) we have

$$\frac{b}{2}(1+\sqrt{1-16b})=\frac{b}{2} \left(1+\frac{4}{\beta^2}\right)=2b\left(\frac{\beta+2i}{2\beta} \right)\left(\frac{\beta-2i}{2\beta} \right)=2\beta_1^2\beta_2\beta_3^2\beta_4.
$$ \smallskip

\noindent Similarly,  $\displaystyle {\frac{b}{2}}\left(1-\sqrt{1-16b} \right)=2\beta_1\beta_2^2\beta_3\beta_4^2$.  This gives the following points of order 4:
$$(0,0),(0,-\beta_1\beta_2\beta_3\beta_4),(-2\beta_1\beta_2\beta_3\beta_4,2\beta_1^2\beta_2\beta_3^2\beta_4),(-2\beta_1\beta_2\beta_3\beta_4,2\beta_1\beta_2^2\beta_3\beta_4^2).\eqno{(4.9)}$$ \smallskip

The simplifications that occur for the other eight points of order 4 are remarkable.  First notice that 
\begin{eqnarray*}
1-16b-(1-8b)\sqrt{1-16b}={\frac{16}{\beta^4}}-\left({\frac{\beta^4+16}{2\beta^4}}\right)\left({\frac{4}{\beta^2}}\right)\\
=-{\frac{2}{\beta^6}}(-8\beta^2+\beta^4+16)=-2\left({\frac{\beta^2-4}{\beta^3}}\right)^2.
\end{eqnarray*}
\smallskip

\noindent Similarly, $\displaystyle 1-16b+(1-8b)\sqrt{1-16b}=2\left({\frac{\beta^2+4}{\beta^3}}\right)^2$.  We will work out the corresponding expression completely for

$$ x=-\frac{1}{8} \left(1-\sqrt{1-16b}\,+\sqrt{2}\sqrt{1-16b-(1-8b)\sqrt{1-16b}}\,\right),$$ \smallskip
 
\noindent  the other cases being similar.  Making use of the above substitutions yields

\begin{eqnarray*}
x=-{\frac{1}{8}}\left(1-{\frac{4}{\beta^2}}+\sqrt{2}\sqrt{2} i \left({\frac{\beta^2-4}{\beta^3}}\right)\right)=-{\frac{1}{8\beta^3}}\left(\beta^3+2i\beta^2-4\beta-8i\right)\\
=-{\frac{1}{8\beta^3}}(\beta+2i)(\beta-2)(\beta+2)=-\beta_1\beta_2\beta_4. \hspace{1.9 in}
\end{eqnarray*}
\smallskip
\noindent Now we use (4.5) to find the corresponding $Y$-coordinates.  Using

$$-1+\beta_3=-\beta_1,\quad 4\beta_2\beta_4-1=-{\frac{4}{\beta^2}}$$ \smallskip
\noindent we can simplify the expression under the square root in (4.5):
\begin{eqnarray*}
x+b=-\beta_1\beta_2\beta_4+\beta_1\beta_2\beta_3\beta_4=\beta_1\beta_2\beta_4(-1+\beta_3)=-\beta_1^2\beta_2\beta_4,\\
4x^2+(x+b)=4\beta_1^2\beta_2^2\beta_4^2-\beta_1^2\beta_2\beta_4=\beta_1^2\beta_2\beta_4(4\beta_2\beta_4-1)=-\frac{4\beta_1^2\beta_2\beta_4}{\beta^2}.
\end{eqnarray*}

\noindent We now have

\begin{eqnarray*}
-\frac{1}{2} \left(x+b \pm \sqrt{(x+b)(4x^2+x+b)}\right)=-\frac{1}{2} \left(-\beta_1^2\beta_2\beta_4\pm \frac{2\beta_1^2\beta_2\beta_4}{\beta} \right)\\
=\frac{\beta_1^2\beta_2\beta_4}{2} \left(1 \mp \frac{2}{\beta}\right)=\beta_1^2\beta_2\beta_4\left(\frac{\beta\mp 2}{2\beta} \right).
\end{eqnarray*}
\smallskip

\noindent Depending on the $\mp$ sign, the last expression is either $\beta_1^2\beta_2^2\beta_4$ or $\beta_1^2\beta_2\beta_4^2$.  These are the two values of $y$ that correspond to $x=-\beta_1\beta_2\beta_4$.  In a similar manner, we obtain the other six points of order 4.  We summarize these results: \bigskip

\noindent {\bf Proposition 4.1.} (See [26].)  If $(\alpha,\beta)$ is a point on the curve $Fer_4$ and $\beta_n$ is defined by (4.6), the 4-torsion points on the elliptic curve 

$$E_1(\alpha):  Y^2+XY+bY=X^3+bX^2 \hspace{.1 in} \textrm{with} \hspace{.1 in}  b={\frac{1}{\alpha^4}}=\beta_1\beta_2\beta_3\beta_4$$

\noindent  and base point $O$, are given by

\begin{eqnarray*}
O,(-\beta_1\beta_2\beta_3\beta_4\,,0),
\left(-{\frac{1}{2}}\beta_1\beta_3,{\frac{1}{2}}\beta_1^2\beta_3^2\right),
\left(-{\frac{1}{2}}\beta_2\beta_4,{\frac{1}{2}}\beta_2^2\beta_4^2\right),\\
(0,0),(0,-\beta_1\beta_2\beta_3\beta_4),(-2\beta_1\beta_2\beta_3\beta_4,2\beta_1^2\beta_2\beta_3^2\beta_4),(-2\beta_1\beta_2\beta_3\beta_4,2\beta_1\beta_2^2\beta_3\beta_4^2),\\
(-\beta_1\beta_2\beta_3,\beta_1^2\beta_2^2\beta_3),(-\beta_1\beta_2\beta_3,\beta_1\beta_2^2\beta_3^2),
(-\beta_1\beta_2\beta_4,\beta_1^2\beta_2^2\beta_4),(-\beta_1\beta_2\beta_4,\beta_1^2\beta_2\beta_4^2),\\
(-\beta_1\beta_3\beta_4,\beta_1^2\beta_3\beta_4^2),(-\beta_1\beta_3\beta_4,\beta_1\beta_3^2\beta_4^2),
(-\beta_2\beta_3\beta_4,\beta_2^2\beta_3^2\beta_4),(-\beta_2\beta_3\beta_4,\beta_2\beta_3^2\beta_4^2).
\end{eqnarray*}

\bigskip

\noindent {\bf Remarks.} 1. The correctness of the above formulas can be verified by noting that if $\beta$ is replaced by $i\beta$, then $\beta_n$ in (4.6) becomes $\beta_{n-1}$.  Thus, the last eight points are permuted among themselves in two orbits of order $4$ on iterating the map $\beta \rightarrow i\beta$.  \smallskip

\noindent 2. From the definition (4.6) of $\beta_n$ we have $\displaystyle 2\beta_n - 1=\frac{2i^{n}}{\beta}$, so $\displaystyle (2\beta_n - 1)^4 = {\frac{16}{\beta^4}}=1-16b.$ \bigskip

Note that for any solution $(\alpha,\beta)$ of $Fer_4$ we have $\beta^4 \neq 16$, so that none of the $\beta_n$ are zero.  Hence, we have \bigskip

\noindent {\bf Proposition 4.2.} The non-zero coordinates of the points in $E_1[4]$ lie in the multiplicative group generated by -1, 2, and $\beta_n$, for $1 \le n \le 4$. \bigskip

\noindent {\bf Remark.} In [26], the first author shows that the points of order $8$ on $E_1$ are also expressible in terms of quantities defined in terms of a solution $(\alpha,\beta)$ of $Fer_4$.  The resulting expressions are similar to those for the points listed in Proposition 4.1, in that the $X$- and $Y$-coordinates are given as products of certain quantities that depend algebraically on $\alpha$ and $\beta$.

\section{Points of order 4 on $E_3$.}

\noindent We shall now find linear fractional expressions for the $X$-coordinates of some of the points of order 4 on the curve
$$E_3(\beta): Y^2 + XY +4bY=X^3 + 16bX^2 + 6bX +b-4b^2, \quad b = \frac{1}{\beta^4}.\eqno{(5.1)}$$

\noindent On this curve we have the doubling formula

$$X(2P)=\frac{x^4-16bx^2-8bx-b}{(4x+1)(x^2+16bx+4b)}, \quad x = X(P).\eqno{(5.2)}$$ \smallskip
We set
$$\gamma_n=\frac{\alpha}{\alpha+\beta\zeta_8^{2n-1}},\quad 1\le n\le 4, \eqno{(5.3)}$$ \smallskip
and note that
$$
b=\frac{1}{\beta^4}=\frac{\alpha^4}{16(\alpha^4+\beta^4)}=\frac{1}{16}\left(\frac{\alpha}{\alpha + \beta\zeta_8} \right)
\left(\frac{\alpha}{\alpha + \beta\zeta_8^3}\right)
\left(\frac{\alpha}{\alpha + \beta\zeta_8^5}\right)
\left(\frac{\alpha}{\alpha + \beta\zeta_8^7}\right).
$$ \smallskip

We now iterate (5.2) and find that the $X$-coordinates of the points of order 4 are zeros of the polynomial

$$x(2x+1)(x^4+32bx^3+24bx^2+8bx+b) .$$ \smallskip

\noindent We shall focus on the roots of the quartic factor.  Treating this factor as we did in Section 4, we complete the square twice to obtain

$$(x^2+16bx+z)^2=(256b^2-24b+2z)x^2+(32bz-8b)x+z^2-b$$ \smallskip

\noindent with the indeterminate $z$.  Setting $z=4b$, the right side becomes a square:

$$(x^2+16bx+4b)^2=b(16b-1)(4x+1)^2;$$ \smallskip

\noindent and using $\displaystyle b(16b-1)= -\frac{16}{\alpha^4 \beta^4}$ we have the two quadratic equations
$$(x^2+16bx+4b) \pm \frac{4i}{\alpha^2\beta^2}(4x+1)=0,$$ \smallskip
or
$$x^2+16 \left(b\pm{\frac{i}{\alpha^2\beta^2}}\right)x+4 \left(b \pm \frac{i}{\alpha^2\beta^2} \right)=0.$$ \smallskip

\noindent Note that
\begin{eqnarray*}
b+\frac{i}{\alpha^2\beta^2}=\frac{1}{\beta^4} \left( \frac{\alpha^2+\beta^2i}{\alpha^2} \right)=b \left( \frac{\alpha+\beta\zeta_8^3}{\alpha} \right) \left( \frac{\alpha+\beta\zeta_8^7}{\alpha} \right)\\
=\frac{1}{16}\gamma_1\gamma_2\gamma_3\gamma_4\cdot \frac{1}{\gamma_2\gamma_4}=\frac{1}{16} \gamma_1\gamma_3. \hspace{1.15 in}
\end{eqnarray*}
\smallskip

\noindent Similarly, $\displaystyle b-{\frac{i}{\alpha^2\beta^2}}={\frac{1}{16}}\gamma_2\gamma_4$.  Thus, the two quadratic equations become

$$x^2+\gamma_1\gamma_3x+{\frac{1}{4}}\gamma_1\gamma_3=0,\quad x^2+\gamma_2\gamma_4x+{\frac{1}{4}}\gamma_2\gamma_4=0.$$ \smallskip

\noindent Using the relations $\gamma_1-2\gamma_1 \gamma_3 + \gamma_3=0$ and $\gamma_2-2\gamma_2 \gamma_4 + \gamma_4=0$ it can be checked directly that the roots of the first equation are $x=-\frac{1}{2} \gamma_1, -\frac{1}{2} \gamma_3$, and the roots of the second are $x=-\frac{1}{2} \gamma_2, -\frac{1}{2} \gamma_4$.  This proves \bigskip

\noindent {\bf Proposition 5.1.} The values

$$x=-\frac{1}{2} \gamma_n= -\frac{1}{2} \frac{\alpha}{\alpha+\beta\zeta_8^{2n-1}}, \quad 1 \le n \le 4, $$ \smallskip

\noindent are the $X$-coordinates of points of order 4 on the curve $E_3=E_3(\beta)$.

\bigskip

\section{Solutions of the Fermat equation.}

\noindent We can now put together the insights we have gained in Sections 2-5.  When $-d$ is square-free, we know that the curve $E_1(\alpha)$ has complex multiplication by the maximal order $R_K$ of $K$ and is defined over $\Sigma$.  Torsion points on this curve of order 2 or 4 will therefore generate low degree abelian extensions of $K$.  The form of the $X$-coordinates of points of order 2 or 4 will allow us to locate simple functions of $\alpha$ and $\beta$ in the field $\Sigma$.  \medskip

The Weierstrass normal form of the curve $E_1(\alpha)$ in (2.10) or (4.1) is

$$Y_1^2=4X_1^3-g_2X_1-g_3=4X_1^3-\bigg(\frac{16b^2-16b+1}{12}\bigg)X_1+\frac{64b^3+120b^2-24b+1}{216},$$

$$\textrm{with} \hspace{.1 in}Y_1=2Y+X+b, \hspace{.1 in} X_1 = X + \frac{4b+1}{12}, \quad b =\frac{1}{\alpha^4}.$$ \smallskip

\noindent From the theory of complex multiplication [33, p. 135] the ray class field $\Sigma_2$ over $K$ with conductor 2 is generated over $K$ by the values $j(E_1)$ and

$$h(P)=\frac{g_2g_3}{\Delta}\bigg (X(P) + \frac{4b+1}{12}\bigg),\eqno{(6.1)}$$ \smallskip

\noindent where $P$ runs over the points of order 2 on $E_1$.  However, the Euler-$\varphi$ function $\varphi(2)=\varphi(\wp_2) \varphi(\wp_2') = 1$ in $K$, so that $\Sigma_2 = \Sigma$.  By Proposition 3.1 we know that $g_2, g_3$ and $\Delta=b^4(1-16b)$ lie in $\Sigma$.  Letting $P$ be the last point in (4.8), we deduce that $X(P)=-(\beta^2-4)/(8\beta^2)$ lies in $\Sigma$, and it follows that $\beta^2 \in \Sigma$. \medskip

We now do the same for the curve

$$E_3(\beta): \hspace{.1 in}  Y^2 + XY + 4bY = X^3 + 16bX^2 +6bX + b-4b^2, \hspace{.1 in} b = \frac{1}{\beta^4},$$ \smallskip
whose $j$-invariant we know to be equal to 
$$\frac{(\beta^8-256\beta^4+4096)^3}{\beta^{16}(16-\beta^4)}=\frac{(\alpha^8-256\alpha^4+4096)^3}{\alpha^{16}(16-\alpha^4)}=j\bigg(\frac{w}{4}\bigg), \eqno{(6.2)} $$ \smallskip

\noindent by virtue of (2.9) and the defining condition for $\beta$.  Thus, $E_3(\beta)$ has complex multiplication by $R_K$.  The Weierstrass normal form of $E_3(\beta)$ has coefficients in $\Sigma$, so once again the $X$-coordinates of points of order 2 on $E_3(\beta)$ must lie in $\Sigma$.   From (5.2) these coordinates are
$$X=-1/4, \hspace{.05 in} -8b \pm 2\sqrt{16b^2-b} = -\frac{8}{\beta^4} \pm 2 \sqrt{\frac{16-\beta^4}{\beta^8}}= -\frac{8}{\beta^4} \pm \frac{8i}{\alpha^2 \beta^2}.$$ \smallskip
It follows that $i \alpha^2 \beta^2 \in \Sigma$, whence we have that $i \alpha^2 \in \Sigma$.  From Proposition 3.1 we now obtain \bigskip

\noindent {\bf Proposition 6.1.} With notation as in Proposition 3.1, the number

$$i \alpha^2 = y(w)= \frac{\eta(w/4)^4}{\eta(w)^4}$$ \smallskip

\noindent lies in the Hilbert class field $\Sigma$, along with the number $\beta^2$ defined by $16\alpha^4 + 16\beta^4=\alpha^4 \beta^4$. \bigskip

\noindent {\bf Remark.} This proposition improves upon a result of Schertz [31], according to which $(\eta(w/4)/\eta(w))^8 \in \Sigma$.  See also [32, p. 159]. \bigskip

Consider next the point $P=(-\beta_1\beta_2\beta_3,\beta_1^2\beta_2^2\beta_3)$, which has order 4 on $E_1(\alpha)$.  Its $X$-coordinate,

$$X(P)=-\beta_1 \beta_2 \beta_3 = -\frac{\beta+2i}{2\beta} \frac{\beta-2}{2\beta} \frac{\beta-2i}{2\beta}=-\frac{\beta^2+4}{4\beta^2} \frac{\beta-2}{2\beta},$$ \smallskip

\noindent lies in the field $\Sigma_4$.  It is not hard to see that $\Sigma_4 = \Sigma(i)$.  First, we have $\varphi(4)=\varphi(\wp_2^2 \wp_2'^2)=4$, and the units $\pm 1$ are incongruent (mod 4) in $R_K$, so $[\Sigma_4: \Sigma]=2$.  (See [20].)  Since $X(P) \in \Sigma_4$, Proposition 6.1 implies that $ \frac{\beta-2}{2\beta}$, and therefore $\beta$, lies in $\Sigma_4$.  But the formulas in Proposition 4.1 imply that $\beta_1=(\beta+2i)/(2\beta)$ also lies in $\Sigma_4$, so $i \in \Sigma_4$.  The field $\Sigma$ does not contain $i=\sqrt{-1}$ because the prime $2$ is not ramified in $\Sigma$.  Hence,  $[\Sigma(i): \Sigma]=2$ shows that $\Sigma_4 = \Sigma(i)$.  \medskip

We set $\beta = r+is$, with $r,s \in \Sigma$ and use $\beta^2 = r^2 - s^2 + 2rsi \in \Sigma$ to conclude that $rs=0$.  Thus, either $\beta$ or $i\beta$ lies in $\Sigma$.  Since $\beta$ has only been determined up to multiplication by a power of $i$, we may assume that $\beta \in \Sigma$. \medskip

Finally, we examine the $X$-coordinate of the following point $Q$ of order 4 on $E_3(\beta)$:

$$X(Q)=-\frac{1}{2} \left(\frac{\alpha}{\alpha+\beta\zeta_8^7}\right)=-\frac{1}{2} \left(\frac{\zeta_8 \alpha}{\zeta_8\alpha+\beta}\right). $$ \smallskip

\noindent Since $\beta \in \Sigma$ and $E_3(\beta)$ has complex multiplication by $R_K$ we deduce from this expression that $\zeta_8 \alpha \in \Sigma_4$.  Then $i \alpha^2 \in \Sigma$ and $\zeta_8 \alpha = u + iv$ with $u,v \in \Sigma$ implies that $uv=0$; so either $\zeta_8 \alpha \in \Sigma$ or $\zeta_8^3 \alpha \in \Sigma$.  Set $\alpha_1 = \zeta_8^j \alpha \in \Sigma$, $ j = 1$ or $3$.  Then $\alpha^4=-\alpha_1^4$ so the equation
$$16\alpha_1^4 - 16 \beta^4 = \alpha_1^4 \beta^4$$ \smallskip
\noindent shows that $16X^4 - 16Y^4 = X^4 Y^4$ has a solution in $\Sigma$.  Further, $(\alpha,\beta)$ is a solution of $Fer_4$ in $\Sigma_4(\zeta_8)=\Sigma(\zeta_8)$. \bigskip

\noindent {\bf Theorem 6.2.} The numbers $\beta$ and $\zeta_8^j \alpha$, for $j = 1$ or $j=3$, lie in the Hilbert class field $\Sigma$.  The point $(X,Y) = (\zeta_8^j \alpha, \beta)$ is a nontrivial solution of the equation
$$16X^4 = X^4Y^4 + 16Y^4$$ \smallskip
\noindent in the field $\Sigma$.  Thus, $x^4 + y^4 = z^4$ has a solution in the Hilbert class field of $Q(\sqrt{-d})$ whenever the discriminant $-d \equiv 1$ (mod 8).  \bigskip

We will convert $(X,Y) = (\zeta_8^j \alpha, \beta)$ into a solution of the form (1.2) using the following lemma.   \bigskip

\noindent {\bf Lemma 6.3.}  If the prime ideal factorization of $2$ in $R_K$ is $(2)=\wp_2 \wp_2'$, where $\wp_2=(2,w)$ with $w=(v+\sqrt{-d})/2$ (as in Section 3), then $(\alpha)=\wp_2'^2$ in $\Sigma(\zeta_8)$ and $(\beta) = (2)\wp_2'=\wp_2 \wp_2'^2$ in $\Sigma$.   \medskip

\noindent {\it Proof.} The assertion about $(\alpha)$ follows from (2.7) and classical results, since by (2.7) the ideal generated by $(\alpha^{12})$ is
$$(\alpha^{12})=\left(\frac{\eta(w/4)}{\eta(w/2)} \frac{\eta(w/2)}{\eta(w)} \right)^{24}=\left(2^{12}\frac{\Delta(\wp_2^2)}{\Delta(\wp_2)}\right) \left( 2^{12}\frac{\Delta(\wp_2)}{\Delta(1)}\right) = \wp_2'^{12} \wp_2'^{12}=\wp_2'^{24}.$$
See [19, Satz 10] or [32, p. 114].  For the second assertion, the equation $\alpha^4(\beta^4-16)=16\beta^4$ implies that $\wp_2' | \beta$ and $16(\alpha^4 + \beta^4) = \alpha^4 \beta^4$ implies that $\wp_2 | \beta$.  Thus, $2 | \beta$, and
$$ \left(\frac{ \beta}{2} \right)^4 -1 = \left( \frac{\beta}{\alpha} \right)^4 \eqno{(6.3)}$$ \smallskip
\noindent implies that $\alpha | \beta$.  The last equation implies further that $\beta/2$ and $\beta/\alpha$ have no prime divisor in common.  In particular, $\beta$ is not divisible by any prime divisor of an odd prime in ${\bf Z}$.  Also, $\beta/2$ cannot be divisible by a prime divisor of $\wp_2$ because $\beta/\alpha$ is, and $\beta/\alpha$ cannot be divisible by a prime divisor of $\wp_2'$ because $\beta/2$ is.  Hence, $(\beta) = 2\wp_2'$ in $\Sigma$ and $(\beta/2)^2$ is a unit times $\alpha$.  $\square$ \bigskip

We may re-write (6.3) in the form

$$ \left(\frac{ \beta}{2} \right)^4 + \left( \frac{\beta}{\zeta_8^j \alpha} \right)^4 = 1,\eqno{(6.4)}$$ \smallskip

\noindent where $(\beta/2)=\wp_2'$ and $(\beta/(\zeta_8^j \alpha)) = \wp_2$.  Thus our solution of the Fermat quartic comes from ideal generators $\pi=\beta/(\zeta_8^j \alpha)$ and $\xi= \beta/2$ of $\wp_2$ and $\wp_2'$ in the Hilbert class field.   This proves the first assertion of Theorem 1.1 of the Introduction and shows that the solution (1.1) is one of an infinite family.  \medskip

All the same arguments apply in case $-d=-f^2d_1$, where $f$ is odd and $d_1\equiv 7$ (mod 8) is square-free.  In that case the field $\Sigma$ is replaced by the ring class field $\Omega_f$ of conductor $f$ over $K$, and we use Proposition 3.3 along with the following result from the theory of complex multiplication.  \bigskip

\noindent {\bf Proposition 6.4.}  Assume the elliptic curve $E$ is defined over the field $\Omega_f$ and has complex multiplication by the order $\textsf{R}_{-d}$ of $K$.  If $\mathfrak{m}$ is an integral ideal of $\textsf{R}_{-d}$  with $(f, \mathfrak{m})=1$, then the $X$-coordinates of points in $E[\mathfrak{m}]$ generate the field $\Sigma_{\mathfrak{m}} \Omega_f$ over $\Omega_f$, where $\Sigma_{\mathfrak{m}}$ is the ray class field of conductor $\mathfrak{m}$ over $K$.  \bigskip

This is an easy consequence of Theorem 2 of [13].  For our argument we take $\mathfrak{m}=(4)=4\textsf{R}_{-d}$.  Then the condition $(4,f)=1$ implies by this proposition that the $X$-coordinates of points of either $E_1(\alpha)[4]$ or $E_3(\beta)[4]$ generate the field $\Sigma_4 \Omega_f=\Omega_f(i)$ over $\Omega_f$.  This field is the same as $\Omega_{4f}$, since a degree argument shows that $\Sigma_4=\Omega_4$ and $\Omega_4 \Omega_f = \Omega_{4f}$.  The above arguments show that $\beta$ and $\zeta_8^j \alpha $ lie in $\Omega_f$ for $j=1$ or $3$, and lead to the following theorem.  \bigskip

\noindent {\bf Theorem 6.5.} Let $-d=-f^2d_1$, with odd $f$ and square-free $d_1 \equiv 7$ (mod 8).  Further, let $(2)=\wp_2 \wp_2'$ be the factorization of $2$ into conjugate prime ideals $\wp_2=(2, (v+\sqrt{-d})/2)$ and $\wp_2'=(2, (v-\sqrt{-d})/2)$ in $R_K$, corresponding to the factorization of the ideal (2) in the order $\textsf{R}_{-d}$ of discriminant $-d$ and conductor $f$ in $K$ .  In the ring class field $\Omega_f$  of $K$ the generators $\pi_f=\beta/(\zeta_8^j \alpha)$ and $\xi_f=\beta/2$ of the principal ideals $\wp_2 R_{\Omega_f}$ and $\wp_2'R_{\Omega_f}$ satisfy
$$ \pi_f^4 + \xi_f^4 = 1.\eqno{(6.5)}$$
Each of $\pi_f$ and $\xi_f$ generates the ring class field $\Omega_f$ over $K$.  $\square$ \bigskip

The last assertion of this theorem follows from (6.2), since $j(w/4)$ generates $\Omega_f$ over $K$.  We will show in Section 8 that each of the numbers $\pi_f$ and $\xi_f$ generates $\Omega_f$ over $\mathbb{Q}$.  This theorem shows that we have a lattice of solutions to the quartic Fermat equation over $K$, since the numbers $\pi_f$ satisfy $\pi_f | \pi_{f'}$ in $\Omega_{f'}$ if $f | f'$.  In fact, $\pi_{f'}/\pi_f$ is a unit in $\Omega_{f'}$. \medskip

For later use we require the following more general form of Lemma 6.3, whose proof carries over word for word. (See [10, pp. 32-33].) \bigskip

\noindent {\bf Lemma 6.6.}  Assume $d=f^2d_1$, where $f$ is odd and $d_1 \equiv 7$ (mod 8) is square-free. If $(2)=\wp_2 \wp_2'$ in $R_K$, where $\wp_2=(2,w)$ with $w=(v+\sqrt{-d})/2$ (as in Section 3), then $(\alpha)=\wp_2'^2$ in $\Omega_f(\zeta_8)$ and $(\beta) = (2)\wp_2'=\wp_2 \wp_2'^2$ in $\Omega_f$.   \bigskip  

\section{Application to the Legendre normal form.}

\noindent We now apply the results of Section 6 to the $\lambda$-parameters of the Legendre normal form
$$E_\lambda: Y^2 = X(X-1)(X-\lambda),$$
for which the curve $E_\lambda$ has complex multiplication by the order $\textsf{R}_{-d}$ of $K$, where $-d=-f^2d_1 \equiv 1$ (mod 8).  Values of $\lambda$ satisfying this requirement may be found as roots of the equation

$$2^8 (\lambda^2-\lambda+1)^3-j(w/2)(\lambda^2-\lambda)^2=0,\eqno{(7.1)}$$ \smallskip

\noindent since $j(w/2)$ is a root of the class equation $H_{-d}(x)=0$ (see (2.8) and the first paragraph of Section 3).  Moreover, any $\lambda$ for which $E_\lambda$ has complex multiplication by $\textsf{R}_{-d}$ is a conjugate of a root of (7.1) over $K$.  Using (2.8) it is easy to see that $\lambda_1=\alpha^4/16$ is one of the roots of this equation.  It follows that the other roots are related to $\alpha^4/16$ by the anharmonic group. \medskip

\noindent (i) We note that $\displaystyle \lambda_1=-\left( \frac{\zeta_8^j \alpha}{2} \right)^4$  for $j=1$ or $3$.  Thus $-\lambda_1$ and $-1/\lambda_1$ are fourth powers in the field $\Omega_f$. \medskip

\noindent (ii) The value $\displaystyle \lambda_2 = 1 - \lambda_1 = 1 - \frac{\alpha^4}{16} = \frac{16-\alpha^4}{16}=-\frac{\alpha^4}{\beta^4}=\left( \frac{\zeta_8^j \alpha}{\beta} \right)^4$, so $\lambda_2$ and $1/\lambda_2$ are fourth powers in $\Omega_f$. \medskip

\noindent (iii) The value $\displaystyle \lambda_3 = 1-\frac{1}{\lambda_1} = \frac{\alpha^4 - 16}{\alpha^4} = \left( \frac{2}{\beta} \right)^4$, so $\lambda_3$ and $1/\lambda_3$ are fourth powers in $\Omega_f$. \medskip

Now we use Lemma 6.6.  Since $(\lambda_1)=(\wp_2'/\wp_2)^4$, it is clear that neither $\lambda_1$ nor $1/\lambda_1$ are algebraic integers, while $1/\lambda_2$ and $1/\lambda_3$ are algebraic integers. Thus we have: \bigskip

\noindent {\bf Theorem 7.1.} Let $\lambda$ be a complex number for which the Legendre normal form $E_\lambda$ has complex multiplication by the order $\textsf{R}_{-d}$ of $K$, where $-d=-f^2d_1 \equiv 1$ (mod 8).  Then $\lambda$ lies in the ring class field $\Omega_f$ of $K$ with conductor $f$.  If either $\lambda$ or $1/\lambda$ is an algebraic integer, then $\lambda$ is a fourth power in the field $\Omega_f$.  If neither $\lambda$ nor $1/\lambda$ is an algebraic integer, then $-\lambda$ is a fourth power in $\Omega_f$.  \bigskip

The above computations and Lemma 6.6 also show that the only prime ideal divisors of $\Omega_f$ dividing the numerator or denominator of $\lambda$ or $\lambda-1$ in this theorem are prime divisors of 2.

\section{Further properties of the solution.}

\noindent We investigate further properties of the numbers $\alpha$ and $\beta$ by making the substitution $\displaystyle \alpha^4=(16\beta^4)/(\beta^4-16)$ in the $j$-invariant

$$j(w) = \frac{(\alpha^8-16\alpha^4+16)^3}{\alpha^8-16\alpha^4}$$ \smallskip

\noindent from (2.8).  This substitution gives that

$$j(w)= \frac{(\beta^8+224\beta^4+256)^3}{\beta^4(\beta^4-16)^4}.\eqno{(8.0)}$$ \smallskip

\noindent Thus, $\beta$ is a root of the polynomial

$$F_d(x)= \left( x^4(x^4-16)^4 \right)^{h(-d)} H_{-d}\left(  \frac{(x^8+224x^4+256)^3}{x^4(x^4-16)^4} \right).\eqno{(8.1)}$$ \smallskip

\noindent Let $B(x)=B_d(x)$ denote the minimal polynomial of $\beta$ over $\mathbb{Q}$.  We now prove: \bigskip

\noindent {\bf Proposition 8.1.}  The degree of $\beta$ over $\mathbb{Q}$ is $2h(-d)$, so that $\Omega_f = \mathbb{Q}(\beta)$.  The numbers $\beta$ and $-\beta$ are not conjugate over $\mathbb{Q}$, and $\beta^4$ has degree $2h(-d)$ over $\mathbb{Q}$. \medskip

\noindent {\it Proof.} Note that $\mathbb{Q}(\beta)$ contains $\mathbb{Q}(j(w))$ and is contained in $\Omega_f$, so $\beta$ has degree $n$ over $\mathbb{Q}$ satisfying $h(-d) \le n \le 2h(-d)$.  If $n=h(-d)$, then $\beta$ has to be real, by the choice of $w = (v+\sqrt{-d})/2$.  This implies in turn that $j(w/4)$ is real, by (6.2). However, $j(w)=j(\textsf{R}_{-d})$ is the $j$-invariant associated to the ideal class 1 in the ring class group (mod $f$) of $K$, while $j(w/4)=j(\wp_{2,-d}^2)$, where $\wp_{2,-d}=\wp_2 \cap \textsf{R}_{-d}$.   Since $\overline{ j(\mathfrak{a})}=j(\mathfrak{a}^{-1})$ for any proper ideal $\mathfrak{a}$ in $\textsf{R}_{-d}$, $j(w/4)$ can be real if and only if $\wp_2^2$ has order 1 or 2 in the ring class group (mod $f$), so that the order of $\wp_2$ must divide 4.  Now $\wp_2^4 \sim 1$ implies that $4 \cdot 2^4 = x^2 +d y^2$ with $x \equiv y$ (mod 2) and $x \equiv y \equiv 2$ (mod 4) if $x$ and $y$ are even.  This equation has a solution if and only if $d=7, 15, 39, 55$ or $63$.  We have the following class equations for these discriminants (see [14, vol. III] and [36]):

\begin{eqnarray*}
H_{-7}(x) = x+ 3375,\hspace{2 in} \\
H_{-15}(x) = x^2 +191025x - 121287375, \hspace{1 in}\\
H_{-39}(x)=x^4+331531596x^3-429878960946x^2+109873509788637459x \\
+20919104368024767633,\\
H_{-55}(x) = x^4+13136684625x^3-20948398473375x^2+172576736359017890625x\\
-18577989025032784359375,\\
H_{-63}(x)=x^4+67515199875x^3-193068841781250x^2+4558451243295023437500x\\
-6256903954262253662109375.
\end{eqnarray*}

\noindent Factoring the polynomial

$$x^{16h(-d)}(16-x^4)^{h(-d)}H_{-d} \left(\frac{(x^8-256x^4+4096)^3}{x^{16}(16-x^4)} \right)$$

\noindent (see (6.2)) in each of these cases shows that \bigskip

$\displaystyle B_7(x) = x^2 +2x+8$ or $x^2-2x+8$; \medskip

$\displaystyle B_{15}(x) = x^4+8x^3+20x^2+16x+64$ or $ x^4-8x^3+20x^2-16x+64$; \medskip

$\displaystyle B_{39}(x) = x^8+12x^7+168x^6+480x^5+848x^4+1728x^3+1536x^2+4096$ 

\hspace{.5 in} or $x^8-12x^7+168x^6-480x^5+848x^4-1728x^3+1536x^2+4096$; \medskip

$\displaystyle B_{55}(x) = x^8+12x^7+312x^6-672x^5+848x^4-2112x^3-768x^2+3072x+4096$

\hspace{.5 in} or $x^8-12x^7+312x^6+672x^5+848x^4+2112x^3-768x^2-3072x+4096$. \bigskip

$\displaystyle B_{63}(x) = x^8+40x^7+440x^6-800x^5+784x^4-2560x^3-2560x^2+5120x+4096$

\hspace{.5 in} or $x^8-40x^7+440x^6+800x^5+784x^4+2560x^3-2560x^2-5120x+4096$. \bigskip

\noindent Hence, the degree of $\beta$ over $\mathbb{Q}$ is $2h(-d)$ in all cases.  \medskip

\noindent If $\beta$ were conjugate to $-\beta$ over $\mathbb{Q}$, then $\beta^2$ would have degree $h(-d)$ over $\mathbb{Q}$, and since $j(w)$ and $j(w/4)$ are rational expressions in $\beta^2$, the same argument as above would show that $\beta^2$ is real and $d = 7, 15, 39, 55$ or $63$.  The above polynomials $B_d(x)$ show in these cases that $\beta$ and $-\beta$ have different minimal polynomials.  A similar argument shows that $\beta^4$ cannot have degree $h(-d)$ over $\mathbb{Q}$, so its degree must be $2h(-d)$.  This proves the proposition.  $\square$ \bigskip

As in [27, p. 256] we set $\displaystyle r(x) = \frac{(x^8+224x^4+256)^3}{x^4(x^4-16)^4}$ and use the fact that the algebraic extension $ \mathbb{Q}(i,x)/\mathbb{Q}(i,r(x))$ is normal with a Galois group $\hat{G}_{24}$ isomorphic to the octahedral group.  Since $\beta$ generates the field $\Omega_f$ over $\mathbb{Q}$ and $\Omega_f \bigcap \mathbb{Q}(i) = \mathbb{Q}$ (2 is unramified in $\Omega_f$), it is clear that $B(x)$ is irreducible over $\mathbb{Q}(i)$.  Now the group $\hat{G}_{24}$ is generated by the substitutions

$$x \rightarrow ix, \quad x \rightarrow \frac{2x+4}{x-2}.$$ \smallskip

\noindent Since the polynomial $H_{-d}(z)$ is irreducible over $\mathbb{Q}(i)$, it corresponds to a prime divisor $\mathfrak{p}$ of the field $\mathbb{Q}(i,z)$.  The irreducible factors of $F_d(x)$ over $k=\mathbb{Q}(i)$ correspond to the prime divisors of the field $k(x), z = r(x)$, which extend $\mathfrak{p}$.  As noted above, the extension $k(x)/k(z)$ is finite and normal with Galois group $\hat{G}_{24}$.  Hence, $\hat{G}_{24}$ is transitive on the irreducible factors of $F_d(x)$ over $k$.  Since these factors have degree $2h(-d)$ over $k$, by Proposition 8.1, the orbit of $B(x)$ has 12 elements, and the stabilizer of $B(x)$ has order 2. \medskip

We now let $S$ denote the subgroup of $\hat{G}_{24}$ consisting of the linear fractional mappings

$$x \rightarrow \pm x, \quad \pm \frac{4}{x}, \quad \pm \frac{2(x+2)}{x-2}, \quad \pm \frac{2(x-2)}{x+2}. \eqno{(8.2)}$$ \smallskip

\noindent Note that $S$ is a 2-Sylow subgroup of $\hat{G}_{24}$ and that it is isomorphic to the dihedral group of order 8.  It consists of all the linear fractional maps in $\hat{G}_{24}$ which don't involve $i=\sqrt{-1}$ as a coefficient.  \bigskip

\noindent {\bf Proposition 8.2.} The stabilizer of $B(x) = B_d(x)$ in $\hat{G}_{24}$ is one of the two sugroups

$$S_1=\{(x \rightarrow x), \left(x \rightarrow \frac{2(x+2)}{x-2}\right) \} \quad \textrm{or} \quad S_2=\{ (x \rightarrow x), \left(x \rightarrow \frac{-2(x-2)}{x+2} \right) \}.\eqno{(8.3)}$$ \smallskip

\noindent One of these subgroups is the stabilizer of $B_d(x)$ and the other is the stabilizer of $B_d(-x)$. \medskip

\noindent {\it Proof.} The polynomial $B(x)$ is a normal polynomial, since its root $\beta$ generates the normal extension $\Omega_f/\mathbb{Q}$.  The stabilizer of $B(x)$ must therefore be a subgroup of $S$, because the other elements of $\hat{G}_{24}$ map the roots of $B(x)$ to numbers which do not lie in $\Omega_f$.  For example, the mapping

$$x \rightarrow \frac{-2i(ix+2)}{ix-2} = \frac{-2i(x-2i)}{x+2i}$$ \smallskip

\noindent cannot map $\beta$ to an element of $\Omega_f$, because otherwise

$$-2i \beta-4 =  \beta \gamma +2i \gamma, \quad \textrm{with} \quad \gamma \in \Omega_f,$$ \smallskip

\noindent implies $\gamma = -\beta$ and $\beta = \pm 2 \in \mathbb{Q}$, which is not the case.  Now the maps $x \rightarrow \frac{-2(x+2)}{x-2}$ and $x \rightarrow \frac{2(x-2)}{x+2}$ cannot fix $B(x)$ because they have order 4 and the stabilizer has order 2.  Furthermore, the maps $x \rightarrow \pm 4/x$ cannot fix $B(x)$ because the principal ideal generated by $\beta$ in the ring of integers of $\Omega_f$ is $2 \wp_2' = \wp_2 \wp_2'^2$ and the principal ideal generated by $4/\beta$ is $\wp_2$ (by Lemma 6.6).  Hence $\beta$ and $\pm 4/\beta$ are divisible by different numbers of prime ideals in $\Omega_f$ and therefore cannot be conjugate over $\mathbb{Q}$.  Finally, Proposition 8.1 shows that $\beta$ and $-\beta$ are not conjugate over $\mathbb{Q}$, so only two nontrivial mappings remain in (8.2) as possibilities.  If one of these lies in the stabilizer of $B(x)$, then the other lies in the stabilizer of $B(-x)$.  This proves the proposition.  $\square$ \bigskip

\noindent  {\bf Corollary.} Assume that the stabilizer of $B_d(x)$ is the group $S_1$ in (8.3).  The roots of the minimal polynomial $b_d(x)=B_d(2x)/2^{2h(-d)}$ of $\beta/2$ over $\mathbb{Q}$ are mapped into themselves by the mapping $x \rightarrow \frac{x+1}{x-1}$. \bigskip

It follows from this corollary and $(\beta/2)=\wp_2'$ that $\left( (\beta+2)/(\beta-2) \right)=\wp_2$, because $\beta/2+1$ and $\beta/2-1$ are relatively prime to $\wp_2'$.  This gives the following result.  \bigskip

\noindent {\bf Proposition 8.3.} Choose the sign of $\beta$ so that the stabilizer of $B_d(x)$ is the group $S_1$ in (8.3).  Then the number $\displaystyle \gamma=\frac{\beta(\beta+2)}{4(\beta-2)}$ is a unit in $\Omega_f$ with degree $h(-d)$ over $\mathbb{Q}$, and $\mathbb{Q}(\gamma)=\mathbb{Q}(j(w))$.  \medskip

\noindent {\it Proof.} The fact that $\gamma$ is a unit follows from the ideal factorization

$$\left( \frac{\beta}{2} \right) \left( \frac{\beta+2}{\beta-2} \right) = \wp_2' \wp_2 = (2).$$ \smallskip

\noindent Since the map $\sigma: \beta/2 \rightarrow (\beta+2)/(\beta-2)$ has order 2, it is clear that $\gamma$ is invariant under $\sigma$.  Its degree over $\mathbb{Q}$ is therefore at most $h(-d)$.  If $\gamma$ had degree $k < h(-d)$ over $\mathbb{Q}$, then $\beta$ would satisfy an equation of degree $2k$ over $\mathbb{Q}$, which contradicts Proposition 8.1.  Moreover, the automorphism $\beta \rightarrow 2(\beta+2)/(\beta-2)$ fixes $j(w)$ by (8.0),  which implies $\mathbb{Q}(j(w)) = \mathbb{Q}(\gamma)$, since $j(w)$ and $\gamma$ have the same degree over $\mathbb{Q}$. $\square$ \bigskip

The units in Proposition 8.3 are similar to the units discussed in [36], since they have degree $h(-d)$ over $\mathbb{Q}$ and generate $\Omega_f$ over $K$, and their minimal polynomials have small coefficients (though generally not as small as the coefficients of the polynomials $W_{-d}(x)$ in [36]).  See Theorem 8.6 below and the examples in Section 12.  \medskip

Note that the rational function $r(x)$ can be expressed as:

$$r(x) = \frac{(x^8+224x^4+256)^3}{x^4(x^4-16)^4}=\frac{(4z^2+1)^3(4z^2-8z+1)^3}{z^4(2z-1)^4}, \quad z =\frac{x(x+2)}{4(x-2)},$$ \smallskip
so we have the formula
$$j(w) = \frac{(4\gamma^2+1)^3(4\gamma^2-8\gamma+1)^3}{\gamma^4(2\gamma-1)^4}.$$ \bigskip

\noindent {\bf Proposition 8.4.}  The numbers $\beta/2$ and $\pm \beta/(\zeta_8^j \alpha)$ are conjugate over $\mathbb{Q}$. \medskip

\noindent {\it Proof.}  From the arguments which precede Theorem 7.1 we know that the numbers

$$\lambda_1=-\left( \frac{\zeta_8^j \alpha}{2} \right)^4, \quad \frac{1}{\lambda_2}=\frac{1}{1-\lambda_1}=\left( \frac{\beta}{\zeta_8^j \alpha} \right)^4, \quad \frac{1}{\lambda_3}=\left( \frac{\beta}{2} \right)^4\eqno{(8.4)}$$ \smallskip

\noindent are roots of the polynomial

$$L_d(x)=(x^2-x)^{2h(-d)}H_{-d}\left( \frac{2^8(x^2-x+1)^3}{(x^2-x)^2}\right).$$ \smallskip

\noindent Furthermore, the numbers in (8.4) have degree $2h(-d)$ over $\mathbb{Q}$, by Proposition 8.1 and the fact that $\alpha^4$ and $\beta^4$ generate the same extension of $\mathbb{Q}$.  From Lemma 6.6 it is clear that the minimal polynomials of $\lambda_1, 1/\lambda_3$, and $\lambda_3$ over $\mathbb{Q}$ are all distinct, because these numbers do not have conjugate prime ideal factorizations:

$$(\lambda_1)=\frac{\wp_2'^4}{\wp_2^4}, \quad \left(\frac{1}{\lambda_3} \right) = \wp_2'^4, \quad (\lambda_3)=\frac{1}{\wp_2'^4}.$$ \smallskip

\noindent  Denoting the respective minimal polynomials by $g_i(x)$ we have the factorization $L_d(x)=g_1(x)g_2(x)g_3(x)$ since the degrees match.  The roots of $g_2(x)$ are the roots of $L_d(x)$ which are algebraic integers; the roots of $g_3(x)$ are the roots of $L_d(x)$ whose reciprocals are algebraic integers; and the roots of $g_1(x)$ are the remaining roots of $L_d(x)$.  Since $1/\lambda_2$ and $1/\lambda_3$ are algebraic integers, they are both roots of the polynomial $g_2(x)$.  Hence, $\left( \frac{\beta}{\zeta_8^j \alpha} \right)^4$ and $\left( \frac{\beta}{2} \right)^4$ are conjugate over $\mathbb{Q}$.  It follows that $\beta/2$ and $\pm  \beta/(\zeta_8^j \alpha)$ are conjugate over $\mathbb{Q}$. $\square$  \bigskip

\noindent {\bf Remarks.}  1. Having chosen the sign of $\beta$ so that the corollary to Proposition 8.2 holds, we may replace the root of unity $\zeta_8^j$ by $\zeta_8^{j+4}$, if necessary, to guarantee that $\xi=\beta/2$ and $\pi=\beta/(\zeta_8^j \alpha)$ are conjugates, where $j \in \{1, 3, 5, 7 \}$.  The residue class of $j$ (mod 8) is uniquely determined by this requirement.  \medskip

\noindent 2. It is not hard to see that $\pm \pi \pm \xi \neq 1$.   If, for example, $\pi + \xi = 1$, then $\pi = 1-\xi$ and $\xi$ are conjugates, so the minimal polynomial $b_d(x)$ over $\mathbb{Q}$ of $\xi$ satisfies $b_d(1-x)=b_d(x)$.  Then the composition of the maps $\sigma=(x \rightarrow 1-x)$ and $\tau=\left( x \rightarrow (x+1)/(x-1) \right)$ would leave the roots of $b_d(x)$ invariant.  But this composition is $\sigma \tau = \left( x \rightarrow  -2/(x-1) \right)$, which has infinite order in the group of M\"obius transformations.  This is not possible if $d > 7$, since the only fixed points of powers of $\sigma \tau$ have degree 2 over $\mathbb{Q}$.   A similar argument applies to show that $\pm \pi \pm \xi \neq 1$.  This remark proves the second assertion of Theorem 1.1.  It can be shown in the same way that the point $(\pi, \xi)$ does not lie on {\it any} rational line.  \bigskip

\noindent {\bf Proposition 8.5.} Choose the sign of $\beta$ and the root of unity $\zeta_8^j$ so that $\displaystyle \frac{\beta+2}{\beta-2}$ and $\displaystyle \pi = \frac{\beta}{\zeta_8^j \alpha}$ are conjugate to $\xi=\beta/2$ over $\mathbb{Q}$.  If $\displaystyle \tau=\left( \frac{\Omega_f/K}{\wp_2} \right)$ is the automorphism associated to the ideal $\wp_2$ by the Artin map for $\Omega_f/K$, then

$$\pi=\frac{\beta^{\tau^{-2}}+2}{\beta^{\tau^{-2}}-2}=\frac{\xi^{\tau^{-2}}+1}{\xi^{\tau^{-2}}-1}.$$ \smallskip
\noindent Equivalently,
$$\xi=\frac{\pi^{\tau^{2}}+1}{\pi^{\tau^{2}}-1}.$$ \smallskip

\noindent {\bf Remark.} This proposition verifies the last assertion of Theorem 1.1. \medskip

\noindent {\it Proof.} Let $\psi$ be the automorphism in $\Gamma=Gal(\Omega_f/\mathbb{Q})$ for which

$$\xi^\psi = \left( \frac{\beta}{2} \right)^\psi=\pi = \frac{\beta}{\zeta_8^j \alpha}$$

\noindent and $\sigma \in \Gamma$ the automorphism for which

$$\xi^\sigma =  \left( \frac{\beta}{2} \right)^\sigma = \frac{\beta+2}{\beta-2}.$$ \smallskip

\noindent These are both well-defined automorphisms of $\Omega_f=\mathbb{Q}(\beta)$ over $\mathbb{Q}$ and they both have order 2 since they interchange $\wp_2$ and $\wp_2'$.  Furthermore, we have

$$(\beta^4)^\psi=-16 \frac{\beta^4}{\alpha^4}=16-\beta^4$$
and
$$\beta^\sigma=\frac{2(\beta+2)}{\beta-2}.$$ \smallskip

\noindent By (8.2), (8.0), and (6.2) we have

$$j(w)^{\sigma \psi}=\left( \frac{(\beta^8+224\beta^4+256)^3}{\beta^4(\beta^4-16)^4} \right)^{\sigma \psi}=\left( \frac{(\beta^8+224\beta^4+256)^3}{\beta^4(\beta^4-16)^4} \right)^{\psi}$$
$$=\frac{(\beta^8-256\beta^4+4096)^3}{\beta^{16}(16-\beta^4)} =j\left(\frac{w}{4} \right).$$ \smallskip

\noindent Hence, $\sigma \psi$ is an automorphism in $Gal(\Omega_f/K)$ taking $j(w)$ to $j(w/4)$.  On the other hand,

$$j\left( \frac{w}{4} \right)^{\tau^2}=j(\wp_{2,-d}^2)^{\tau^2} \equiv j(\wp_{2,-d}^2)^4 \equiv j(\wp_{2,-d}^{-2} \wp_{2,-d}^2) = j(w) \ (\textrm{mod} \ \wp_2),$$ \smallskip

\noindent by Hasse's congruence [19, Satz 11], [10, p.34]; and therefore $\displaystyle j\left( \frac{w}{4} \right)^{\tau^2}=j(w)$, since $p=2$ does not divide the discriminant of $H_{-d}(x)$ (see [8]).  It follows that $\sigma \psi=\tau^{-2}$ and therefore $\psi=\sigma \tau^{-2}$.  This implies the assertion of the proposition. $\square$  \bigskip

\noindent {\bf Theorem 8.6.} If the sign of $\beta$ is chosen as in Proposition 8.3, so that $\beta$ and $\displaystyle \frac{2(\beta+2)}{\beta-2}$ are conjugates over $\mathbb{Q}$, then: \medskip

a) the quantities $\beta^3(\beta^2-4)$ and $-\beta(\beta^2+4)$ are fourth powers in $\Omega_f$; \medskip

b) the unit $\displaystyle \gamma=\frac{\beta(\beta+2)}{4(\beta-2)}$ is a square in $\mathbb{Q}(\gamma)$; and  \medskip

c) the quantity $\beta(\beta-2)/8$ is a unit, with $N_{\mathbb{Q}(\gamma)}(\beta(\beta-2)/8) = \gamma$.  \medskip

d) $\displaystyle 2\gamma-1=\frac{\beta^2+4}{2(\beta-2)}$ is a unit. \bigskip

\noindent {\it Proof.} From Proposition 8.5 we have that

$$\pi=\frac{\beta}{\alpha_1}=\frac{\beta^{\tau^{-2}}+2}{\beta^{\tau^{-2}}-2},$$
so Theorem 1.1 implies the equation
$$(\beta^{\tau^{-2}}+2)^4 + \left( \frac{\beta}{2} \right)^4 (\beta^{\tau^{-2}}-2)^4 = (\beta^{\tau^{-2}}-2)^4.$$ \smallskip
From this it follows that
$$ \left( \frac{\beta}{2} \right)^4 (\beta^{\tau^{-2}}-2)^4 = (\beta^{\tau^{-2}}-2)^4-(\beta^{\tau^{-2}}+2)^4 = -16 \beta^{\tau^{-2}}(\beta^{2\tau^{-2}}+4).$$ \smallskip
Applying the automorphism $\tau^2$ to both sides gives
$$\left( \frac{\beta^{\tau^2}}{2} \left(\frac{\beta}{2}-1\right)\right)^4 = -\beta(\beta^2+4).\eqno{(8.5)}$$ \smallskip
Now we apply the automorphism $\sigma$ taking $\beta$ to $\displaystyle \frac{2(\beta+2)}{\beta-2}$:
$$\left( \frac{\beta^{\tau^2 \sigma}}{2} \left( \frac{\beta^\sigma}{2}-1\right)\right)^4 = -\beta^\sigma((\beta^\sigma)^2+4)=\frac{-16(\beta+2)(\beta^2+4)}{(\beta-2)^3}.\eqno{(8.6)}$$ \smallskip
The right side of this equation is equal to
$$\frac{-16(\beta+2)(\beta^2+4)}{(\beta-2)^3}=\frac{16\beta^3(\beta^2-4) \times [-\beta(\beta^2+4)]}{\beta^4(\beta-2)^4}.\eqno{(8.7)}$$ \smallskip

\noindent Equations (8.5) and (8.6) now clearly imply that $\beta^3(\beta^2-4)$ is a fourth power in $\Omega_f$, as claimed.  This proves part a). \medskip

It is immediate that $\gamma$ is a square in $\Omega_f$.  To prove part b) we have to show that $\sqrt{\gamma}$ is invariant under the automorphism $\sigma$.  From (8.5)-(8.7) we have

$$\sqrt{\gamma}=\frac{\beta}{8} \frac{\beta^{2\tau^2 \sigma} (\beta^\sigma-2)^2}{\beta^{2\tau^2}(\beta-2)},$$
so
$$\sqrt{\gamma}^\sigma=\frac{\beta^\sigma}{8} \frac{\beta^{2\tau^2} (\beta-2)^2}{\beta^{2\tau^2 \sigma}(\beta^\sigma-2)}=\frac{1}{\sqrt{\gamma}} \frac{\beta \beta^\sigma}{64} (\beta-2)(\beta^\sigma-2)=\frac{1}{\sqrt{\gamma}} \frac{\beta(\beta+2)}{4(\beta-2)}=\sqrt{\gamma}.$$
This proves the assertion.  \medskip

The last equation makes it clear that $\gamma$ is the norm of the quantity $\displaystyle \frac{\beta (\beta-2)}{8}$ in the extension $\Omega_f/\mathbb{Q}(\gamma)$.  To see that this number is also a unit, recall from Lemma 6.6 that $(\beta)= \wp_2 \wp_2'^2$, so that $(\beta-2)=\wp_2 \wp_2' \mathfrak{a}$, where $(\mathfrak{a},\wp_2')=1$.  Now note that
$$(\beta^\sigma-2)(\beta-2) = \left( \frac{2\beta+4}{\beta-2}-2 \right) (\beta-2) = 8.$$ \smallskip
\noindent Thus, the norm to $L=\mathbb{Q}(\gamma)$ of $(\beta/2-1)=\mathfrak{a}$ is $(2)=\wp_2 \wp_2'$.  If a prime ideal $\mathfrak{p}$ of $\textsf{R}_{\Omega_f}$ divides $\mathfrak{a}$, then its conjugate $\mathfrak{p}^\sigma$ divides $\wp_2'$.  But $\mathfrak{p}^2$ cannot divide $\mathfrak{a}$ because $N_{\Sigma/L}(\mathfrak{a})=(2)$ is square-free in $L$.  Now $\mathfrak{a} \mathfrak{a}^\sigma = \wp_2 \wp_2'$ and $\mathfrak{a} | \wp_2$ implies that $\mathfrak{a}=\wp_2$.  Hence, $(\beta-2)=\wp_2^2 \wp_2'$, and this implies that $\displaystyle \frac{\beta (\beta-2)}{8}$ is a unit.  This proves part c).  \medskip

To prove d), first use $(\beta-2)=\wp_2^2 \wp_2'$, $(\beta+2)=(4(\beta-2)/\beta)$, and Lemma 6.6 to deduce that $(\beta+2)=\wp_2^3 \wp_2'$.  Then use $\beta^4-16=16\beta^4/\alpha^4$ to see that $(\beta^4-16) = \wp_2^8 \wp_2'^4$.  Putting these facts together gives that $(\beta^2+4)=\wp_2^3 \wp_2'^2 = (2(\beta-2))$.  This proves the assertion.  $\square$ \medskip

\noindent {\bf Corollary 1.} \noindent The point $(\beta, 2(\beta-2)\sqrt{\gamma})$ is an integral point on the elliptic curve $E: Y^2 = X(X^2-4)$ with coordinates in the ring class field $\Omega_f$ of $K$ with conductor $f$. \bigskip

\noindent {\bf Corollary 2.} \noindent The point $\displaystyle \left(-\beta,\frac{1}{16} \beta^{2\tau^2} (\beta-2)^2 \right)$ is an integral point on the elliptic curve $\bar E: Y^2 = X(X^2+4)$ with coordinates in the ring class field $\Omega_f$. \bigskip

\noindent Note that the curves in Corollary 1 and 2 are curves (64B) and (32B), respectively, in Cremona's tables [6].  \bigskip

\noindent {\bf Corollary 3.} The algebraic integer $-\beta$ is the sum of two squares in the ring $R_{\Omega_f}$.  Moreover, $\beta$ is a sum of two squares in $\Omega_f$ if and only if the order of $\wp_{2,-d}=\wp_2 \cap \textsf{R}_{-d}$ in the class group of $\textsf{R}_{-d}$ is even, or equivalently, iff the order of $\wp_2$ in the ring class group (mod $f$) is even. \medskip

\noindent {\it Proof.} By part a) of the theorem, we have
$$-8\beta^3=\beta^3(\beta^2-4)-\beta^3(\beta^2+4)=A^4+\beta^2B^4,$$
for some integers $A, B$ in $\Omega_f$.  Thus,
$$-\beta=\left(\frac{A^2+\beta B^2}{4 \beta} \right)^2+\left(\frac{-A^2+\beta B^2}{4\beta} \right)^2,$$
where the ideals $(A)^4=\wp_2^8 \wp_2'^8$ and $(B)^4=\wp_2^4 \wp_2'^4$ by the computations in the above proof; hence the expressions inside the squares are algebraic integers.  It follows that $\beta$ is a sum of two squares in $\Omega_f$ if and only if $-1$ is a sum of two squares, which holds if and only if the level ({\em Stufe}) of $\Omega_f$ is 2.  The characterization of the level $s(\Omega_f)$ in terms of the order of $\wp_{2,-d}$ in $Pic(\textsf{R}_{-d})$ follows from [18, Satz 14] (see also [12]), since the absolute degree of the prime divisors of $2$ in $\Omega_f$ is equal to the order of $\wp_2$ in the ring class group (mod $f$) of $K$.
$\square$  \bigskip

\noindent {\bf Corollary 4.} In the ring $R_{\Omega_f}$ the algebraic integer $j(w)$ is a cube times the fourth power of a unit. \medskip

\noindent {\it Proof.} This follows from Part d) of Theorem 8.6 and the formula
$$j(w) = \frac{(4\gamma^2+1)^3(4\gamma^2-8\gamma+1)^3}{\gamma^4(2\gamma-1)^4}.$$
\noindent $\square$ \bigskip

\noindent {\bf Remark.} We also note the formula
$$j(w) = \frac{(2\gamma+i)^3(2\gamma-i)^3(2\gamma-2+\sqrt{3})^3(2\gamma-2-\sqrt{3})^3}{\gamma^4(2\gamma-1)^4},$$ \bigskip
which expresses $j(w)$ as $\gamma^{-4}$ times a product of powers of quantities of the form $(2\gamma-\varepsilon)$, where $\varepsilon$ is a unit.  \bigskip

To prove that $\alpha_1=\zeta_8^j \alpha$ is a square in $\Omega_f$ we start with the following.  \bigskip

\noindent {\bf Theorem 8.7.} If $\alpha$ is any algebraic number for which the curve $E_1(\alpha)$ (see (2.10)) has complex multiplication by the order $\textsf{R}_{-d}$ in $K=\mathbb{Q}(\sqrt{-d})$, where $-d =-f^2d_1\equiv 1$ (mod $8$), then $\alpha$ lies in the class field $\Sigma_{16} \Omega_f$ over $K$, where $\Sigma_{16}$ is the ray class field of conductor $16$ over $K$. \bigskip

\noindent {\it Proof.} Assume that $E_1=E_1(\alpha)$ has complex multiplication by $\textsf{R}_{-d}$.  From (8.0) and (8.1) and the discussion preceding (8.2), it is clear that any algebraic number $\beta$ for which $(\alpha, \beta)$ lies on $Fer_4$ must lie in the field $\Omega_f(i)$.  This is because $\beta$ is then a root of the polynomial $F_d(x)$ in (8.1), and its minimal polynomial is in the orbit of the polynomial $B(x)=B_d(x)$ under the action of the group $\hat G_{24}$.  Hence, $b=1/\alpha^4=1/16-1/\beta^4$ lies in $\Omega_f(i)$ also, and therefore $E_1$ is defined over this field.  It follows from Proposition 6.4 (see [13]) that the $X$-coordinates of points of order $16$ on $E_1$ lie in $\Sigma_{16} \Omega_f$, since $i \in \Sigma_4 \subset \Sigma_{16}$.  \medskip

Now we use the fact that there is a cyclic isogeny $\rho: E_1 \rightarrow E_3(\alpha)$, where $\rho= \phi \circ \psi$ is given by the equations immediately preceding equation (2.12).  Let $Q$ be the point on $E_3=E_3(\alpha)$ whose $X$-coordinate is
$$\xi=-\frac{1}{2} \gamma_1 = - \frac{1}{2} \frac{\beta}{\beta+\zeta_8 \alpha},$$ \smallskip
from Proposition 5.1 (switching the roles of $\alpha$ and $\beta$).  Let $P$ be a point on $E_1$ for which $\rho(P)=Q$.  Then $[4]Q=O$ on $E_3$ implies that $\rho([4]P) = O$ on $E_3$, and therefore, since $\rho$ has degree $4$, that $P \in E_1[16]$.  Since the coefficients of the rational function $\phi_1 \circ \psi_1$ depend only on $b$, they lie in the field $\Omega_f(i)$.  Hence, $\xi = \phi_1 \circ \psi_1(X(P))$ shows that $\xi \in \Sigma_{16} \Omega_f$ and therefore $\zeta_8 \alpha \in \Sigma_{16} \Omega_f$. \medskip

On the other hand, it is easy to see that $\zeta_8 \in \Sigma_8 \subset \Sigma_{16}$, by the transference theorem ({\em Verschiebungssatz}) of class field theory (see [20]).  By that theorem, $K(\zeta_8)$ is class field over $K$ for the ideal group $H$ consisting of the ideals $\mathfrak{a}$ of $K$ (prime to $2$) whose norms to $\mathbb{Q}$ are $1$ (mod $8$).  However, any principal ideal in $R_K$ with a generator congruent to $1$ (mod $8$) satisfies this condition.  This shows that the ideal group corresponding to $\Sigma_8$ lies in $H$, and therefore $K(\zeta_8) \subset \Sigma_8$.  (A similar argument shows of course that $\zeta_{16} \in \Sigma_{16}$ also.)  \medskip

It follows that $\alpha \in \Sigma_{16} \Omega_f$.  $\square$ \bigskip

Computations on Maple suggest the following stronger result. \bigskip

\noindent {\it Conjecture.} If $\alpha$ is an algebraic number for which the curve $E_1(\alpha)$ has complex multiplication by the order $\textsf{R}_{-d}$ in $K=\mathbb{Q}(\sqrt{-d})$, where $-d =-f^2d_1\equiv 1$ (mod $8$), then $\alpha$ lies in the ring class field $\Omega_{16} \Omega_f =\Omega_{16f}$ over $K$.  \bigskip

Next, we need the following result.  \bigskip

\noindent {\bf Theorem 8.8.} (See [26].)  Let $\displaystyle j=j(\alpha)=\frac{(\alpha^8-16\alpha^4+16)^3}{\alpha^4(\alpha^4-16)}$
be the $j$-invariant of the curve $E_1(\alpha)$, where $\alpha$ is an {\it indeterminate}.  Also, let $k$ be any field whose characteristic is not 2 that contains a primitive $8$-th root of unity $\zeta_8=(1+i)/\sqrt{2}$.  Then the normal closure of the algebraic extension $k(\alpha)/k(j)$ is the function field
$$N = k\left(\zeta_{16} \beta^{1/2}, \zeta_{16} \beta^{1/4}(\beta^2-4)^{1/4}, \beta^{1/4}(\beta^2+4)^{1/4} \right),$$ 
where $(\alpha, \beta)$ satisfies the equation $Fer_4: 16\alpha^4+16\beta^4=\alpha^4 \beta^4$.  Moreover, $k(\alpha, \beta) \subset N$.  \medskip

\noindent {\it Proof.} Let $F(x,j)$ denote the polynomial
$$F(x,j)=(x^2-16x+16)^3-jx(x-16).$$
We are looking for the splitting field $N$ of $F(x^4,j)=0$ over $k(j)$.  Note first that
$$F(8 \pm \sqrt{y+48},j)=y^3-jy+16j=G(y,j).$$
From (8.0) we have
$$j=\frac{(\alpha^8-16\alpha^4+16)^3}{\alpha^4(\alpha^4-16)}=\frac{(\beta^8+224\beta^4+256)^3}{\beta^4(\beta^4-16)^4}.$$
With this substitution for $j$ in $G(y,j)$, the latter polynomial factors over $k(\beta)$ as follows:
$$G(y,j)=\frac{1}{\beta^4(\beta^4-16)^4} \left((\beta^6-8\beta^4+16\beta^2)y+\beta^8+224\beta^4+256 \right)$$
$$\ \ \ \ \times \left((\beta^8-32\beta^4+256)y-16\beta^8-3584\beta^4-4096 \right)$$
$$\times \left((\beta^6+8\beta^4+16\beta^2)y-\beta^8-224\beta^4-256 \right).$$ \smallskip

\noindent Denote the respective roots of $G(y,j)=0$ by
$$\xi_1=-\frac{\beta^8+224\beta^4+256}{\beta^2(\beta^2-4)^2},$$
$$\xi_2=16\frac{\beta^8+224\beta^4+256}{(\beta^4-16)^2},$$
$$\xi_3=\frac{\beta^8+224\beta^4+256}{\beta^2(\beta^2+4)^2}.$$ \smallskip

\noindent Now a straightforward calculation on Maple shows that the roots of $F(x,j)=0$ are:
$$8 \pm \sqrt{\xi_1+48} =\frac{i(\beta-2i)^4}{\beta(\beta^2-4)}, \ \ \frac{-i(\beta+2i)^4}{\beta(\beta^2-4)};$$
$$8 \pm \sqrt{\xi_2+48} =\frac{16\beta^4}{\beta^4-16}, \ \  \frac{-256}{\beta^4-16};$$
$$8 \pm \sqrt{\xi_3+48} =\frac{(\beta+2)^4}{\beta(\beta^2+4)}, \ \ \frac{-(\beta-2)^4}{\beta(\beta^2+4)}.$$ \smallskip

\noindent The third root listed above is $8+\sqrt{\xi_2+48}=\alpha^4$.  Noting that $\zeta_8 \in k$, this shows that the splitting field of $F(x^4,j)=0$ is generated over $k(j)$ by the elements
$$\beta, \ \ \zeta_{16}\beta^{1/4}(\beta^2-4)^{1/4}, \ \ (\beta^4-16)^{1/4}, \ \ \beta^{1/4}(\beta^2+4)^{1/4}.$$ \smallskip
Since
$$\frac{\left(\zeta_{16}\beta^{1/4}(\beta^2-4)^{1/4} \right) \left(\beta^{1/4}(\beta^2+4)^{1/4} \right)}{(\beta^4-16)^{1/4}} = \zeta_{16}\beta^{1/2},$$
this proves the assertion.   $\square$
\bigskip

Now take $k=\Sigma_{16} \Omega_f$ in Theorem 8.8.  As noted in the proof of Theorem 8.7, this field satisfies the hypothesis of Theorem 8.8.   \medskip

If $j_0$ is any root of the class equation $H_{-d}(x)=0$, and $P_{j_0}$ is the prime divisor of the function field $k(j)$ corresponding to the numerator divisor of $j-j_0$, then any $\alpha_0$ for which $j(\alpha_0)=j_0$ lies in $k=\Sigma_{16} \Omega_{f}$, by Theorem 8.7; because the $j$-invariant  of $E_1(\alpha_0)$ is equal to $j_0$, this curve has complex multiplication by $\textsf{R}_{-d}$.  Hence the polynomial
$$F(x^4,j_0)=(x^8-16x^4+16)^3-j_0 x^4(x^4-16)=0,$$ \smallskip
which is the specialization for $j=j_0$ of the minimal polynomial of $\alpha$ over $k(j)$, splits completely over $k$.  Therefore, the prime divisor $P_{j_0}$ splits into prime divisors of degree 1 in $k(\alpha)$, and it follows that $P_{j_0}$ splits into prime divisors of degree 1 in the normal closure $N$ of $k(\alpha)/k(j)$.  \medskip

This has the following consequence.  Let $\bar P$ be any prime divisor of $N$ lying over $P_{j_0}$.  Then the residue class field of $N$ mod $\bar P$ is just the field $k$.  Hence, the residue mod $\bar P$ (i.e., the specialization) of any element in $N$ which is integral for $\bar P$ lies in $k$.  In particular, if $(\alpha_0, \beta_0)$ is a solution of $Fer_4$ for which $\zeta_8^j \alpha_0$ and $\beta_0$ lie in $\Omega_f$, then because $k(\alpha,\beta) \subset N$ there is a prime divisor $\bar P$ of $N$ for which
$$\alpha \equiv \alpha_0, \quad \beta \equiv \beta_0 \quad (\textrm{mod} \ \bar P),$$ \smallskip
and consequently the residue of $\zeta_{16} \beta^{1/2} \in N$ also lies in $k=\Sigma_{16} \Omega_f$.  Switching back to our usual notation, we have: \bigskip

\noindent {\bf Proposition 8.9.} If $(\alpha, \beta)$ is a solution of $Fer_4$ for which $\zeta_8^j \alpha$ and $\beta$ lie in $\Omega_f$, then the quantity $\beta^{1/2} \in \Sigma_{16} \Omega_f$ generates an abelian extension of $K$. \bigskip

We can now prove \bigskip

\noindent {\bf Theorem 8.10.} If $\beta$ and $\zeta_8^j$ are chosen as in Proposition 8.5, then $\alpha_1=\zeta_8^j \alpha$ is a square in $\Omega_f$. \bigskip

\noindent {\it Proof.} From Proposition 8.5, we know that

$$\frac{\beta}{ \alpha_1} = \frac{\beta^{\tau^{-2}}+2}{\beta^{\tau^{-2}}-2},$$
so that
$$\alpha_1=\frac{\beta(\beta^{\tau^{-2}}-2)}{\beta^{\tau^{-2}}+2}.$$
This equation implies further that
$$\alpha_1=\frac{\beta}{\beta^{\tau^{-2}}} \frac{\beta^{\tau^{-2}}(\beta^{\tau^{-2}}-2)(\beta^{\tau^{-2}}+2)}{(\beta^{\tau^{-2}}+2)^2}.$$ \smallskip
By Theorem 8.6, the quantity
$$\beta^{\tau^{-2}}(\beta^{\tau^{-2}}-2)(\beta^{\tau^{-2}}+2)=f(\beta)^{\tau^{-2}},$$
with $f(\beta)=\beta(\beta^2-4)$, is a square in $\Omega_f$. Hence, $\alpha_1=\zeta_8^j \alpha$ is a square in $\Omega_f$ if and only if the unit $\displaystyle \beta^{1-\tau^{-2}}=\frac{\beta}{\beta^{\tau^{-2}}}$ is a square in $\Omega_f$. \medskip

\noindent Now we appeal to Proposition 8.9.  By that result $L=\Omega_f(\beta^{1/2})$ is an abelian extension of $K$.  There is an isomorphism of $L$ over $K$ which extends the automorphism $\tau^{-2} \in Gal(\Omega_f/K)$ on $\Omega_f/K$.  Since $L$ is normal over $K$, this is an automorphism, which we again denote by $\tau^{-2}$.  Hence,

$$L=L^{\tau^{-2}}=\Omega_f \left(\left(\beta^{1/2}\right)^{\tau^{-2}}\right)=\Omega_f \left(\left(\beta^{\tau^{-2}}\right)^{1/2}\right).$$ \smallskip

\noindent The field $L$ is a quadratic extension of $\Omega_f$, so by Kummer theory we have that
$$\beta=\beta^{\tau^{-2}} \eta^2, \ \ \eta \in \Omega_f.$$ \smallskip
Hence, $\beta^{1-\tau^{-2}}$ is a square in $\Omega_f$. $\square$  \bigskip

\noindent {\bf Corollary.} If $-d =f^2 d_1 \equiv 1$ (mod $8$) and $w$ is chosen as in Proposition 3.3, then there is a unique value of $j \in \{1, 3, 5, 7 \}$ for which
$$\zeta_8^{(j-1)/2} \frac{\eta(w/4)}{\eta(w)} \in \Omega_f.$$  \medskip

\noindent {\bf Remark.} Theorem 10.6 gives a formula for the value of $j$ in this corollary, in the case that $w=(v+\sqrt{-d})/2$ and $v=1$ or $3$.  \bigskip

We can derive one more interesting consequence from Theorem 8.6 and its Corollary 1.  We start with the elementary observation that on the curve
$$E: \ Y^2 = X(X^2-4)$$
the sum of the points $P_d=\left(\beta,2(\beta-2) \sqrt{\gamma} \right)$ and $(2,0)$ is
$$P_d+(2,0)=\left(\beta,2(\beta-2) \sqrt{\gamma} \right)+(2,0)=\left( \frac{2(\beta+2)}{\beta-2}, \frac{-16 \sqrt{\gamma}}{\beta-2} \right)=-P_d^{\sigma},$$
since $2(\beta^\sigma-2)=16/(\beta-2)$ and $\sqrt{\gamma}^\sigma=\sqrt{\gamma}$, by the proof of Theorem 8.6b).  This gives that
$$P_d^\sigma=(2,0)-P_d.\eqno{(8.8)}$$
Now let $f=1$, so that $\Omega_1=\Sigma$ is the Hilbert class field of $K$, and let $Q_K \in E(\Sigma)$ be the trace of the point $P_d$ to the field $K$:
$$Q_K=\sum_{\sigma_1 \in G}{P_d^{\sigma_1}}, \ \ \ G = Gal(\Sigma/K).\eqno{(8.9)}$$
It is clear that $Q_K$ has coordinates in $K$, since it is invariant under the group $G$.  We shall prove that $Q_K$ is a nontrivial point of $E(K)$, i.e. that $Q_K$ does not lie in $E(\mathbb{Q})= \{O, (0,0), (2,0), (-2,0) \}$, whenever the class number $h(K)$ of the field $K$ is odd.  \bigskip

\noindent {\bf Theorem 8.11.} If $d=p \equiv 7$ (mod $8$) is a prime, then the point $Q_K$ is a nontrivial point on $E: \ Y^2=X(X^2-4)$ which is defined over the field $K=\mathbb{Q}(\sqrt{-p})$. \medskip

\noindent {\it Proof.}  First we apply the automorphism $\sigma$ to the point $Q_K$:
$$Q_K^\sigma=\sum_{\sigma_1 \in G}{P_d^{\sigma_1 \sigma}}=\sum_{\sigma_1 \in G}{P_d^{\sigma \sigma_1^{-1}}}$$
$$=\sum_{\sigma_1 \in G}{(2,0)}-\sum_{\sigma_1 \in G}{P_d^{\sigma_1^{-1}}}$$
by (8.8) and (8.9), so that
$$Q_K^\sigma=[h(K)](2,0)-Q_K.\eqno{(8.10)}$$
Now assume that $Q_K$ has coordinates in $\mathbb{Q}$.  Then the last equation yields
$$2Q_K=[h(K)](2,0) = (2,0),\eqno{(8.11)}$$
since the class number of $K$ is odd.  If $Q_K=(x,y)$, then the doubling formula for $E$ gives that
$$\frac{1}{4} \frac{(x^2+4)^2}{x^3-4x}=2,$$
and therefore $x=2 \pm 2\sqrt{2} \not \in \mathbb{Q}$.  Hence, equation (8.11) is impossible and $Q_K \not \in E(\mathbb{Q})$.  \ \ $\square$  \bigskip

\noindent {\bf Remark.}  When $d=23$, the point $Q_K$ is computed to be
$$Q_K=\left( \frac{-7+\sqrt{-23}}{9}, \ \frac{-50+2 \sqrt{-23}}{27} \right).$$
Note that $[2]Q_K=(49/36, -77 \sqrt{-23}/216)$, so that there is a solution over $\mathbb{Q}$ of the twisted curve $-23Y^2=X(X^2-4)$.  A similar circumstance holds in the general case, since (8.10) implies that $([2]Q_K)^\sigma=-[2]Q_K$ and therefore $[2]Q_K=(x,y \sqrt{-p})$, with $x, y \in \mathbb{Q}$.  An elementary argument similar to the above shows that $[2]Q_K \neq (0,0), (-2,0)$, so that $xy \neq 0$.  This gives: \bigskip

\noindent {\bf Corollary.} If $p \equiv 7$ (mod $8$) is a prime, the equation $E_{p}: \ pY^2=X(X^2-4)$ has a nonzero rational solution.  \bigskip

Note that the above proof of Theorem 8.11 and its corollary is algebraic, except for the use of the properties of the Schl\"afli functions in Section 2.  (See [9] for Deuring's algebraic derivation of the theorems of complex multiplication.)  \bigskip

For a general $d$ with $d\equiv 7$ (mod $8$) it is easy to see that the point $P_d=(\beta,2(\beta-2) \sqrt{\gamma})$ has infinite order on the curve $E$.  This is because $E$ has complex multiplication by $\mathbb{Z}[i]$, and is defined over $\mathbb{Q}$, so that the squares of the $X$-coordinates of torsion points on $E$ generate abelian extensions of $k=\mathbb{Q}(i)$.  If $P_d$ had finite order, then $k(X(P_d)^2)=k(\beta^2)=\mathbb{Q}(i,\beta^2)=\Sigma(i)$, by Proposition 8.1, so $\Sigma(i)=\Sigma_4$ would be abelian over $k$.  Since $\Sigma \cap k=\mathbb{Q}$, we have that $Gal(\Sigma_4/k) \cong Gal(\Sigma/\mathbb{Q})$, so $\Sigma$ would be abelian over $\mathbb{Q}$.  In that case, $\Sigma$ coincides with its genus field, which corresponds in the sense of class field theory to the group of squares in the class group of $K$.  Hence, the class group would have exponent 2.  Then $\wp_2^2 \sim 1$ implies that $16=x^2+dy^2$ for $x, y \in \mathbb{Z}$, $x \equiv y$ (mod $2$), and $xy \neq 0$.  It follows that $d=7$ or $15$.  In both of these cases it can be checked directly that $P_d$ has infinite order on $E$.  This proves part a) of the following theorem.  \bigskip

\noindent {\bf Theorem 8.12.} a) For any $d \equiv 7$ (mod $8$), the point $P_d=\left(\beta,2(\beta-2) \sqrt{\gamma} \right)$ on $E$ has infinite order. \medskip

\noindent b) If $d \equiv 7$ (mod $8$), then the only torsion points on the curve $E$ in $K=\mathbb{Q}(\sqrt{-d})$ are the points in $E(\mathbb{Q})=E[2]$.  \medskip

\noindent c) If the point $Q_K \not \in \{O, (0,0), (2,0), (-2,0) \}$, then $Q_K$ has infinite order on $E$. \medskip

\noindent {\it Proof of b).} Assume that the point $Q \in E(K)$ has finite order, and that $Q \not \in E(\mathbb{Q})$.  Then the coordinates $(x,y)$ of $Q$ are algebraic integers, since the prime $2$ splits in $K$ (see [34, Thm. 7.1, p. 240]).  Furthermore, we know $y^2 \mid -4^4$, so $y \mid 16 \cong \wp_2^4 \wp_2'^4$, by the arguments of [34, Cor. 7.2, p. 240].  Now $y$ cannot be a rational integer, since the polynomial $X^3-4X-y^2$ is irreducible in the variable $X$ over $\mathbb{Q}$ for any integer $y$ dividing $16$, so its roots have degree 3 and cannot lie in $K$.  Then the power of $\wp_2$ dividing $y$ is not the same as the power of $\wp_2'$ dividing $y$: without loss of generality $(y/2^r) = \wp_2^e \neq (1)$, with $1 \le e+r \le 4$.  It follows that $2^{e+2}=u^2+dv^2$, with $u \equiv v \equiv 1$ (mod $2$).  This implies that $d \le 63$, so $d \in \{7, 15, 23, 31, 39, 47, 55 \}$, from which $47$ can be discarded, because the order of $\wp_2$ in the class group of $\mathbb{Q}(\sqrt{-47})$ is $5$.  This yields only the following solutions and their conjugates:
$$x=\frac{-1+\sqrt{-7}}{2}, \ \ y= \pm \frac{5-\sqrt{-7}}{2}=\pm \left(\frac{-1+\sqrt{-7}}{2} \right)^3;$$
and
$$x=1-\sqrt{-7}, \ \ y=\pm 2^2 \left(\frac{1+\sqrt{-7}}{2} \right).$$
In both cases $X(2Q)=1/4$, so these points have infinite order on $E$.  This proves parts b) and c).   $\square$  \bigskip

\begin{table}
  \centering 
  \caption{Points $Q_K$ on $E(K)$ for $15 \le d \le 103$. \medskip} \label{ }

\begin{tabular}{|c|c|c|c|}
\hline
  &  &  &  \\
$d$	&   {$Q_K$}  & $d$  &  {$Q_K$}  \\
  &  &  &  \\
\hline
% after \\ : \hline or \cline{col1-col2} \cline{col3-col4} ...
   &   &  &  \\
 15  &  $(-3,\sqrt{-15})$  &  $63$  &  $O$  \\
  23 &   $\left( \frac{-7+\sqrt{-23}}{9}, \ \frac{-50+2 \sqrt{-23}}{27} \right)$ & 71 & $\left(\frac{12809-8183\sqrt{-71}}{29241}, \frac{13245170+2373070\sqrt{-71}}{5000211} \right)$ \\
  31 &  $\left(\frac{151+23\sqrt{-31}}{49}, \frac{1060+460\sqrt{-31}}{343} \right)$  & 79 & $\left(\frac{-761-49\sqrt{-79}}{961}, \frac{53660+980\sqrt{-79}}{29791} \right)$  \\
  39 &  $\left(\frac{1}{13}, -\frac{15}{169} \sqrt{-39} \right)$ & 87 & $(-27,15\sqrt{-87})$ \\
  47 & $\left(\frac{23+17\sqrt{-47}}{49}, \frac{-900+204\sqrt{-47}}{343} \right)$ &  95 & $\left( \frac{9}{5}, \frac{3}{25} \sqrt{-95} \right)$  \\
  55 & $\left(\frac{1}{5}, -\frac{3}{25} \sqrt{-55} \right)$ & 103 & $\left(\frac{-1103+161\sqrt{-103}}{1681}, \frac{-151810+5474\sqrt{-103}}{68921} \right)$ \\
   &  &  &  \\
\hline
\end{tabular}

\end{table}

Based on Theorem 8.12 and the examples in Table 1, we make the following conjecture. \bigskip

\noindent {\it Conjecture.} If $d$ is a positive, square-free integer with $d\equiv 7$ (mod $8$) and $K=\mathbb{Q}(\sqrt{-d})$, then the point $Q_K \in E(K)$ defined by (8.9) has infinite order on the curve $E: \ Y^2=X(X^2-4)$. \bigskip

This conjecture holds for all $d$ of the required type for $d < 100$.  If true, the conjecture implies there are $\gtrsim \frac{1}{\pi^2}X$ positive, square-free integers  $d \le X$ for which the quadratic twist $E_{-d}$ has rank $\ge 1$.  (See Landau [23, p. 635] and [4].)  We note that this curve is not one of the curves discussed by Byeon [4], since it has no nontrivial rational $3$-torsion and conductor $N=64$.  Since the twist $E_{-d} \cong E_d$ is isomorphic to the curve $Y^2=X(X^2-4d^2)$, a proof of the above conjecture would give an unconditional proof that square-free positive integers $2d \equiv 14$ (mod $16$) are congruent numbers.  (See Koblitz, {\it Introduction to Elliptic Curves and Modular Forms}, p. 92.)  \bigskip

\bigskip

\section{Connection with supersingular curves.}

\noindent We compare Theorem 7.1 with the result of Theorem 1.2(i) of [27], according to which the values of $\bar \lambda$ for which the Legendre normal form 
$$E_{\bar \lambda}: \quad Y^2=X(X-1)(X-\bar \lambda)$$ \smallskip
\noindent is supersingular in characteristic $p>2$ are fourth powers in the finite field $\mathbb{F}_{p^2}$.  \medskip

If $\mathfrak{p}$ is a prime ideal of $\Omega_f$ lying above the odd, rational prime $p$ for which the Legendre symbol $\left(-d/p \right) = 0$ or $-1$, then the degree of $\mathfrak{p}$ over $\mathbb{Q}$ is either 1 or 2, so that $R_{\Omega_f}/\frak{p} \cong \mathbb{F}_{p}$ or $\mathbb{F}_{p^2}$.  Since $-1$ is a fourth power in $\mathbb{F}_{p^2}$, then for any $\lambda$ satisfying the hypothesis of Theorem 7.1 the residue class $\bar \lambda = \lambda$ (mod $\frak{p}$) will certainly be a fourth power in $\mathbb{F}_{p^2}$.  By Deuring's reduction theory [7], any supersingular curve $E_{\bar \lambda}$ arises by reduction (mod $\mathfrak{p}$) from a curve $E_\lambda$ with complex multiplication by $\textsf{R}_{-d}$, if $\sqrt{-d}$ injects into its endomorphism ring $End(E_{\bar \lambda})$.  To show that Theorem 1.2(i) of [27] can be derived from Theorem 7.1, we prove the following theorem.  \bigskip

\noindent {\bf Theorem 9.1.} If $E_{\bar \lambda}$ is a supersingular curve in characteristic $l \ge 3$, then there exists a positive integer $d \equiv 7$ (mod 8) for which $\sqrt{-d}$ injects into $End( E_{\bar \lambda})$.  \medskip

\noindent Equivalently, let $\textsf{M}$ be any maximal order of the definite quaternion algebra $D_l$ which is ramified only at $l$ and the infinite prime $p_\infty$.  Then some square-root $\sqrt{-d}$ with positive $d \equiv 7$ (mod 8) injects into $\textsf{M}$.  \bigskip

\noindent {\bf Remark.} Note that if $\sqrt{-d} \rightarrow \mu$ injects into the maximal order $\textsf{M}$, then $K=\mathbb{Q}(\sqrt{-d})$ is isomorphic to a maximal commutative subfield $\mathbb{Q}(\mu)$ of $D_l$, and is therefore a splitting field for $D_l$.  By the well-known criterion for splitting fields of central simple algebras [re, Thm. 32.15, p. 278], it is clear that the prime $l$ does not split in $K$, i.e., $\displaystyle \left(\frac{-d}{l} \right)=0$ or $-1$.  Conversely, if this condition holds, then $K$ is a splitting field for $D_l$ and is isomorphic to a maximal subfield of $D_l$ (see [29, Cor. 28.10, p. 240]). \bigskip

\noindent {\it Proof of Theorem 9.1.}  The assertion of Theorem 9.1 is equivalent to the statement that for some square-free integer $d_1 \equiv 7$ (mod $8$) and some embedding $\mathbb{Q}(\sqrt{-d_1}) \rightarrow K' = \mathbb{Q}(\mu) \subset D_l$, the conductor $t$ of the order $\mathcal{O}=\textsf{M} \cap \mathbb{Q}(\mu)$ in the field $K'$ is either odd or exactly divisible by $2$.  In the first case $t\sqrt{-d_1}=\sqrt{-t^2d_1}$ injects into $\textsf{M}$ and $d=t^2d_1 \equiv 7$ (mod $8$); and in the second case the same situation holds with $t$ replaced by $t/2$.  \medskip

We shall use the following characterization of maximal orders in the quaternion algebra $D_l$ first given by Dorman [11], and extended and clarified by Lauter and Viray [25, pp.13-19].  Let $l \ge 3$ be a prime and $d_1=p \equiv 7$ (mod $8$) a prime for which $\displaystyle \left(\frac{-p}{l} \right)=-1$.  Consider the quaternion algebra $\displaystyle D_l=\left( \frac{-lq, -p}{\mathbb{Q}} \right)$, where $q$ is a prime for which $-lq \equiv 1$ (mod $8$) and $\displaystyle \left(\frac{-p}{q} \right) = +1$.  The Hilbert symbol conditions

$$\left( \frac{-lq, -p}{r} \right)=+1, \ \textrm{for all primes} \ r \neq l,$$
and
$$\left( \frac{-lq, -p}{l} \right) =  \left( \frac{-lq, -p}{p_\infty} \right)=-1,$$ \smallskip

\noindent imply that $D_l$ is ramified only at $l$ and the infinite prime $p_\infty$.  Furthermore, setting $L=\mathbb{Q}(\sqrt{-p})$, $D_l$ is isomorphic to the algebra
$$\mathbb{B}=\{[\alpha, \beta]:=\left(\begin{array}{cc}\alpha & \beta \\-l q \bar \beta & \bar \alpha \end{array}\right) | \ \alpha, \beta \in L \},$$
with the embedding $L \rightarrow \mathbb{B}$ defined by $\alpha \in L \rightarrow [\alpha, 0] \in \mathbb{B}$.  \medskip

Now let $\mu = [\sqrt{-p}, 0] \in \mathbb{B}$ and let $\textsf{M}$ be any maximal order in $\mathbb{B}$.  Assume that $\mathcal{O}=\textsf{M} \cap \mathbb{Q}(\mu)$ has conductor $t$, where $\mu^2=-p$.  If $t$ is odd or exactly divisible by $2$, then by the above remarks there is nothing to show.  So we assume $4 \mid t$, and we let $d=-t^2 p$ be the discriminant of $\mathcal{O}$.  Since $q$ splits in $L$, we can write $q\mathcal{O}=\mathfrak{q} \bar \mathfrak{q}$, for a prime ideal $\mathfrak{q}$ of $\mathcal{O}$.  We may also assume $(q,t)=1$, since we can replace $q$ by any other odd prime satisfying the same conditions to obtain an isomorphic algebra.  \medskip

By [25, Lemma 6.10] the maximal order $\textsf{M}$ is conjugate in $\mathbb{B}$ (by an element in $L^\times$) to the maximal order $R(\mathfrak{a},\lambda)$ determined as follows.  The ideal $\mathfrak{a}$ is some integral, invertible ideal of $\mathcal{O}$, relatively prime to the conductor $t$, and $\lambda$ is an element of $\mathcal{O}$ for which \medskip

(1) $\lambda \mathfrak{q}^{-1} \bar \mathfrak{a} \mathfrak{a}^{-1} \subseteq \mathcal{O}$, and \smallskip

(2) $\textsf{N}(\lambda) \equiv -lq$ (mod $d$).  \medskip

\noindent With these definitions, and putting $\mathfrak{D} = (1/\sqrt{d})$, the set
$$R(\mathfrak{a}, \lambda) = \{ [\alpha, \beta]  \ | \  \alpha \in \mathfrak{D}^{-1}, \beta \in \mathfrak{q}^{-1} \mathfrak{D}^{-1} \bar \mathfrak{a} \mathfrak{a}^{-1}, \alpha - \lambda \beta \in \mathcal{O} \}$$
\noindent is a maximal order in $\mathbb{B}$ with $R(\mathfrak{a},\lambda) \cap \mathbb{Q}(\mu)=\mathcal{O}$.  Now we take $\mu_1=[0, \beta] \in R(\mathfrak{a}, \lambda)$, where $\beta \in \bar \mathfrak{a}$ is chosen as follows.  Suppose the ideal $\bar \mathfrak{a}$ has the integral basis $\{a, \omega_1 \}$ in $\mathcal{O}$, where $a=\textsf{N}(\bar \mathfrak{a})$ is prime to $t$, $\omega_1=x+\frac{ty}{2} \sqrt{-p}$, with $x,y \in \mathbb{Z}$, and $a \mid \textsf{N}(\omega_1)$.  Then we let
$$\beta=ra+2s(x+\frac{ty}{2} \sqrt{-p}) \in \bar \mathfrak{a}, \ \ r \equiv 1 \ (\textrm{mod} \ 2), \ \ r, s \in \mathbb{Z}.$$
Then $\textsf{N}(\beta)= (ra+2sx)^2+(sty)^2 p \equiv 1$ (mod $8$), since $4|t$.  It follows that the minimal polynomial of $\mu_1=[0,\beta]$ is $X^2+lq \textsf{N}(\beta)$, where
$d_2=lq \textsf{N}(\beta) \equiv 7$ (mod $8$), by the choices of $q$ and $\beta$.  Then $\mu_1^2=-d_2$, so $\sqrt{-d_2}$ injects into $R(\mathfrak{a},\lambda)$, and therefore also into the isomorphic  maximal order $\textsf{M}$.  This proves the theorem.  $\square$  \bigskip

This theorem verifies the remarks in the Introduction immediately following Theorem 1.3.

\bigskip

\section{The solutions as values of modular functions.}

\noindent The solution $(\pi_f, \xi_f)$ of the Fermat equation that we have given in Theorems 1.1 and 6.5 can be represented in terms of modular functions.   First recall that the Schl\"afli functions $\mathfrak{f}(\tau), \mathfrak{f}_1(\tau), \mathfrak{f}_2(\tau)$ are defined by

$$\mathfrak{f}(\tau)=e^{-\frac{\pi i}{24}} \frac{\eta\left( \frac{\tau + 1}{2} \right) }{\eta(\tau)}, \quad \mathfrak{f}_1(\tau)= \frac{\eta\left( \frac{\tau }{2} \right) }{\eta(\tau)},\quad \mathfrak{f}_2(\tau)= \sqrt{2} \hspace{.05 in} \frac{\eta(2\tau) }{\eta(\tau)},$$ \smallskip

\noindent where $\eta(\tau)$ is the Dedekind $\eta$-function ([5, p. 256], [32, p. 148]).  These functions have the infinite product representations

$$\mathfrak{f}(\tau)=q^{-\frac{1}{48}} \prod_{n=1}^\infty{(1+q^{n-\frac{1}{2}})},\quad \mathfrak{f}_1(\tau)=q^{-\frac{1}{48}} \prod_{n=1}^\infty{(1-q^{n-\frac{1}{2}})},$$ 

$$\mathfrak{f}_2(\tau)=\sqrt{2} \hspace{.05 in} q^{\frac{1}{24}} \prod_{n=1}^\infty{(1+q^n)}, \quad q = e^{2 \pi i \tau}.$$ \smallskip

\noindent From the defining formulas (2.7) and (1.3) we have

$$\alpha^4=-y(w)^2=-\mathfrak{f}_1(w/2)^8 \mathfrak{f}_1(w)^8,$$

$$-\frac{\beta^4}{\alpha^4}=\frac{-16}{\alpha^4-16}=\frac{16}{\mathfrak{f}_1(w/2)^8 \mathfrak{f}_1(w)^8+16}.\eqno{(10.1)}$$ \smallskip

\noindent Now we use the relations $\mathfrak{f}_1(w)\mathfrak{f}_2(w/2)=\sqrt{2}$ and $\mathfrak{f}_1(w/2)^8+\mathfrak{f}_2(w/2)^8=\mathfrak{f}(w/2)^8$ from [35, pp. 114-115] to write the last expression in (10.1) in the form

$$\frac{16}{\mathfrak{f}_1(w/2)^8 \mathfrak{f}_1(w)^8+16}=\frac{\mathfrak{f}_2(w/2)^8}{\mathfrak{f}_1(w/2)^8+\mathfrak{f}_2(w/2)^8}=\frac{\mathfrak{f}_2(w/2)^8}{\mathfrak{f}(w/2)^8}.$$ \smallskip

\noindent Hence, we have 

$$\pi = \frac{\beta}{\zeta_8^j \alpha}=i^a\frac{\mathfrak{f}_2(w/2)^2}{\mathfrak{f}(w/2)^2}, \quad w = \frac{v+\sqrt{-d}}{2}, \quad v^2 \equiv -d \hspace{.05 in} (\textrm{mod} \hspace{.05 in} 16).\eqno{(10.2)}$$ \smallskip

\noindent It follows that

$$\frac{\beta^4}{16}=1+\frac{\beta^4}{\alpha^4}=1-\frac{\mathfrak{f}_2(w/2)^8}{\mathfrak{f}(w/2)^8}=\frac{\mathfrak{f}_1(w/2)^8}{\mathfrak{f}(w/2)^8}.$$ \smallskip

\noindent This gives

$$\xi = \frac{\beta}{2}=i^b\frac{\mathfrak{f}_1(w/2)^2}{\mathfrak{f}(w/2)^2}. \eqno{(10.3)}$$ \smallskip

\noindent From numerical calculations it appears that the factor $i^b$ in this formula is given by

$$i^b=i^{-v}=\cases{i, &if $v=3$ \hspace{.05 in} \textrm{and} \hspace{.05 in} $d \equiv 7$ (mod $16$),\cr
	-i, &if $v=1$ \hspace{.05 in} \textrm{and} \hspace{.05 in} $d \equiv 15$ (mod $16$).}\eqno{(10.4)}$$ \smallskip
We will prove this at the end of this section.\medskip

From [36, p. 1646] or [35, p. 179] we know that $X=\mathfrak{f}^{24}(w/2), -\mathfrak{f}_1^{24}(w/2), -\mathfrak{f}_2^{24}(w/2)$ are the roots of
$$0=(X-16)^3-j(w/2)X=(X-16)^3-\frac{(\beta^8-16\beta^4+256)^3}{\beta^8(\beta^4-16)^2} X,$$
where we have used (2.8) and the defining relation between $\alpha$ and $\beta$.  The right side of the last equation factors, its roots being
$$X=-\frac{2^{12}}{\beta^4(\beta^4-16)}, \ \ \frac{\beta^8}{\beta^4-16}, \ \ -\frac{(\beta^4-16)^2}{\beta^4}.$$
By Lemma 6.6 and the arguments at the end of the proof of Theorem 8.6, the ideal factorizations of these three expressions in $\beta$ are, respectively, $1, \wp_2'^{12}$, and $\wp_2^{12}$.  It follows from (10.3) that we must have
$$\mathfrak{f}^{24}(w/2)=-\frac{2^{12}}{\beta^4(\beta^4-16)}, \ \ \mathfrak{f}_1^{24}(w/2) =-\frac{\beta^8}{\beta^4-16},$$
$$\mathfrak{f}_2^{24}(w/2)=\frac{(\beta^4-16)^2}{\beta^4}.$$ \medskip

Now, from the proof of Theorem 8.6 we have

$$-\beta(\beta^2+4)=\left(\frac{\beta^{\tau^2}}{2} \left(\frac{\beta}{2}-1 \right)\right)^4,\eqno{(10.5a)}$$
$$\beta^3(\beta^2-4)=\left(\frac{4\beta \beta^{\tau^2 \sigma}}{\beta^{\tau^2}(\beta-2)} \right)^4.\eqno{(10.5b)}$$
These formulas give that
$$\mathfrak{f}_1^{24}(w/2)=\frac{\beta^{12}}{\beta^3(\beta^2-4) \times -\beta(\beta^2+4)}=\frac{\beta^{12}}{(\beta \beta^{\tau^2 \sigma})^4},$$
so that
$$\mathfrak{f}_1^6(w/2)=\varepsilon_1 \frac{\beta^3}{\beta \beta^{\tau^2 \sigma}}=\varepsilon_1 \frac{\alpha_1 \beta}{2},\eqno{(10.6a)}$$
for some $4$-th root of unity $\varepsilon_1$.  (This uses the fact that $\beta^\psi=2\beta/\alpha_1$ for the automorphism $\psi=\tau^2 \sigma$ in the proof of Theorem 8.5.)  Similarly, we find that
$$\mathfrak{f}^6(w/2)=\varepsilon_2 \frac{2^3}{\beta \beta^{\tau^2 \sigma}}=\varepsilon_2 \frac{4\alpha_1}{\beta^2}=\varepsilon_2 \frac{2}{\pi \xi},\eqno{(10.6b)}$$
where $\varepsilon_2^4=1$.  These equations and (10.3) imply that $\varepsilon_1=\varepsilon_2 i^{b}$.  From (10.2) and (10.3) we also have that
$$\frac{\pi \xi}{2}=\frac{\beta^2}{4\zeta_8^j \alpha}=i^{a+b}\left( \frac{ \mathfrak{f}_1(w/2) \mathfrak{f}_2(w/2) }{\sqrt{2} \mathfrak{f}(w/2)^2} \right)^2=\frac{i^{a+b}}{\mathfrak{f}(w/2)^6},$$    
where we have used the relation $\mathfrak{f}(\tau) \mathfrak{f}_1(\tau) \mathfrak{f}_2(\tau) = \sqrt{2}$ [35, p. 114].  Comparing with (10.6b) implies that
$$\varepsilon_2=i^{a+b}, \ \ \varepsilon_1=i^{a+2b}.$$ \medskip

It is clear from (10.6b) that $\mathfrak{f}(w/2)$ is a unit, and from (10.6a) that $(\mathfrak{f}_1^2(w/2))=\wp_2'$. \medskip

We now use a result from the paper [36] of Yui and Zagier.  We recall the definition from [36, p. 1647] of the Weber singular modulus for the primitive, positive definite quadratic form $Q(x,y)=ax^2+bxy+cy^2$ with discriminant $-d$: if $\tau_Q$ is the unique root of $Q(x,1)=0$ with positive imaginary part, then

$$f_w(Q)=\cases{ \ \ \ \zeta^{b(a-c-ac^2)} \ \mathfrak{f}(\tau_Q), \ \ \textrm{if} \ (a,c) \equiv (0,0) \ (\textrm{mod} \ 2);\cr 
\varepsilon_d \zeta^{b(a-c-ac^2)} \ \mathfrak{f}_1(\tau_Q), \ \textrm{if} \ (a,c) \equiv (0,1) \ (\textrm{mod} \ 2);\cr
\varepsilon_d \zeta^{b(a-c+a^2 c)} \ \mathfrak{f}_2(\tau_Q), \ \textrm{if} \ (a,c) \equiv (1,0) \ (\textrm{mod} \ 2);}$$ \smallskip

\noindent where $\varepsilon_d=(-1)^{(-d-1)/8}$ and $\zeta=\zeta_{48}=e^{2\pi i/48}$.  For $Q$ we first take the form $Q_2(x,y)=2x^2-vxy+\left(\frac{v^2+d}{8}\right)y^2$, for which $Q_2(w/2,1)=0$, with $\displaystyle w=\frac{v+\sqrt{-d}}{2}$.  The Weber singular modulus for $Q_2$ is

$$f_w(Q_2)=\zeta^{-v(2-(v^2+d)/8-(v^2+d)^2/32)} \mathfrak{f}(w/2). \eqno{(10.7)}$$ \smallskip

\noindent If $3$ does not divide $d$, then by the proposition of [36, p. 1647] we have $f_w(Q_2) \in \mathbb{Q}(j(w/2)) \subset \Omega_f$.  We require a modest extension of this result.
\bigskip

\noindent {\bf Proposition 10.1.} a) If $(d,3)=1$, then for any primitive quadratic form $Q(x,y)=ax^2+bxy+cy^2$ of discriminant $-d \equiv 1$ (mod $8$), the value $f_w(Q)$ depends only on the $SL_2(\mathbb{Z})$ equivalence class $\mathcal{A}=[Q]$ of $Q$ and lies in $\mathbb{Q}(j(\tau_Q))$. \medskip

b) If $3 \mid d$, then the value $f_w(Q)^3$ depends only on the $SL_2(\mathbb{Z})$ equivalence class $\mathcal{A}=[Q]$ of $Q$ and lies in $\mathbb{Q}(j(\tau_Q))$.  \medskip

\noindent {\it Proof.}  Part a) is proved in [36].  For part b), it is only necessary to modify the proof in [36] slightly.  The proof that $f_w(Q)^3$ is an invariant of the class $\mathcal{A}=[Q]$ is the same as the proof in [36], except that only congruences (mod $16$) need to be considered instead of congruences (mod $48$).  For the second part of the proof, choose $Q=[a,b,c]$ in $\mathcal{A}=[Q]$ so that $(a,6)=1$ is odd and $b \equiv -a$ (mod $48$).  Then with the same notation as in [36], $a^2 \equiv 1$ (mod $24$) implies that
$$f_w(\mathcal{A})^3=\pm (\zeta^{b(a-c+a^2c)} \mathfrak{f}_2(\tau_Q))^3=\pm (\zeta^{-1} \mathfrak{f}_2(\tau_Q))^3=\pm \left(\frac{\sqrt{2}}{\mathfrak{f}(2\tau_Q-1)}\right)^3.$$
Now $\varpi=2\tau_Q-1$ is a root of $\tilde Q(x,1)=0$ for the quadratic form
$$\tilde Q =Ax^2+2Bxy+Cy^2, \ \ A = a, \ \ B=a+b, \ \ C=a+2b+4c.$$
(This corrects the value given for $C$ in the last paragraph of the proof in [36, p. 1648].)  By Thm. 6.4.1 of [32, p. 148], which is applicable since $A$ is odd and $B \equiv 0$ (mod $16$), it follows that $(\mathfrak{f}(\varpi)/\sqrt{2})^3 \in \mathbb{Q}(j(\varpi))$.  But the discriminant of $\tilde Q$ is $-4d$, so $j(\varpi)$ generates a subfield of the ring class field $\Omega_{2f}$ of $K$, which coincides with $\Omega_f$.  Hence, $f_w(\mathcal{A})^3=f_w(Q)^3 \in \Omega_f$ has degree at most $h(-d)$.  Since $\mathbb{Q}(j(\tau_Q)) \subseteq \mathbb{Q}(f_w(Q)^3)$ and $j(\tau_Q)$ has degree $h(-d)$, it follows that $\mathbb{Q}(f_w(Q)^3)=\mathbb{Q}(j(\tau_Q))$, which is what we needed to show.  $\square$

\bigskip

We apply this proposition by raising (10.7) to the $6$-th power:
$$f_w(Q_2)^6=i^{-v(1-(v^2+d)/16-(v^2+d)^2/64)} \mathfrak{f}^6(w/2)=i^{-v(1-(v^2+d)/16)} \cdot i^{a+b} \frac{4\alpha_1}{\beta^2},$$
where the second equality follows from (10.6b) and $\displaystyle \frac{(v^2+d)^2}{64} \equiv 0$ (mod $4$).  Since $i \notin \Omega_f$ and both $4\alpha_1/\beta^2$ and $f_w(Q_2)^6$ are squares in $\Omega_f$, this implies
$$\frac{4\alpha_1}{\beta^2} =f_w(Q_2)^6,\eqno{(10.8)}$$
and
$$a+b \equiv v \left(1-\frac{v^2+d}{16} \right) \hspace{.05 in} (\textrm{mod} \hspace{.05 in} 4).\eqno{(10.9)}$$
This proves the following theorem. \bigskip

\noindent {\bf Theorem 10.2.} If $(d,3)=1$ and $f_w(Q_2)$ is defined by (10.7), then
$\displaystyle \frac{4\alpha_1}{\beta^2}$ is the $6$-th power of the unit $f_w(Q_2)$ in $\Omega_f$.  If $3 \mid d$, then $\displaystyle \frac{4\alpha_1}{\beta^2}$ is the square of the unit $f_w(Q_2)^3$ in $\Omega_f$.   \bigskip

Note that $4\alpha_1/\beta^2=2^3/(\beta \beta^{\tau^2 \sigma})$ lies in the fixed field of the automorphism $\psi=\tau^2 \sigma$, and so has degree at most $h(-d)$ over $\mathbb{Q}$.  Since $j(w/2)$ is a rational function of $\mathfrak{f}^{24}(w/2)$, it follows from (10.6b) that $\mathbb{Q}(j(w/2)) \subseteq \mathbb{Q}(4\alpha_1/\beta^2)$ and therefore $\mathbb{Q}(j(w/2)) = \mathbb{Q}(4\alpha_1/\beta^2)$, since the degree of $j(w/2)$ equals $h(-d)$.  Thus, 
$$\mathbb{Q}(j(w/2))=\mathbb{Q}(4\alpha_1/\beta^2) = \mathbb{Q}(f_w(Q_2)^e) = Fix(\psi),$$
where $e=1$ or $3$ and $Fix(\psi)$ denotes the fixed field of $\psi=\tau^2 \sigma=\tau \sigma \tau^{-1}$ inside $\Omega_f$. This agrees with the fact that $\mathbb{Q}(j(w))=Fix(\sigma)$ (see the proof of Proposition 8.3) and that $\mathbb{Q}(j(w/2))=\mathbb{Q}(j(w))^{\tau^{-1}}$.\medskip

\noindent {\bf Remark.} Equation (10.8) clarifies the relationship between the generators of $\Omega_f$ considered by Yui and Zagier in [36] and the algebraic numbers $\alpha_1=\zeta_8^j \alpha$ and $\beta$ that are the main focus of this paper.  It also makes it clear that $4\alpha_1/\beta^2$ is a natural unit to consider in the case that $3$ does divide $d$.  See Case {\bf B} in [36, p. 1660] and the example $d=159$ in Section 12.  \medskip

Next, we consider the quadratic form $Q_1(x,y)=x^2-vxy+\left(\frac{v^2+d}{4}\right)y^2$, for which $Q_1(w,1)=0$.  The formulas of Yui and Zagier give that
$$f_w(Q_1)=(-1)^{(-d-1)/8} \zeta^{-v} \mathfrak{f}_2(w).\eqno{(10.10)}$$ \smallskip
The value $f_w(Q_1)$ (or $f_w(Q_1)^3$) is real since $j(w)$ is real, as the $j$-invariant of the principal class in $\textsf{R}_{-d}$.  To find the connection with our generators, we use (8.0) to solve the equation
$$0=(X-16)^3-j(w)X=(X-16)^3-\frac{(\beta^8+224\beta^4+256)^3}{\beta^4(\beta^4-16)^4} X.$$
Once again, this cubic factors, giving the roots
$$X=-\frac{(\beta^2-4)^4}{\beta^2(\beta^2+4)^2}, \ \ \frac{(\beta^2+4)^4}{\beta^2(\beta^2-4)^2}, \ \ -\frac{2^{12}\beta^4}{(\beta^4-16)^2}.$$
Using Lemma 6.6 and the computations at the end of the proof of Theorem 8.6, we find that these three roots generate the respective ideals in $R_{\Omega_f}$:

$$\left(-\frac{(\beta^2-4)^4}{\beta^2(\beta^2+4)^2}\right)=\wp_2^{12}, \ \ \left(\frac{(\beta^2+4)^4}{\beta^2(\beta^2-4)^2}\right) =(1), \ \ \left(-\frac{2^{12}\beta^4}{(\beta^4-16)^2}\right)=(\wp_2')^{12}.$$ \smallskip

Yui and Zagier prove in [36] that $f_w(Q_1)^2$ and $f_w(Q_2)^2$ are conjugates over $\mathbb{Q}$ (the same argument applies to $f_w(Q_1)^6$ and $f_w(Q_2)^6$ when $3 \mid d$), so $\mathfrak{f}_2(w)$ must be a unit, which forces
$$\mathfrak{f}_2^{24}(w)=-\frac{(\beta^2+4)^4}{\beta^2(\beta^2-4)^2}=-\frac{\beta^4(\beta^2+4)^4}{\beta^6(\beta^2-4)^2}.$$
Using (10.5) yields that
$$\mathfrak{f}_2^{3}(w)=\frac{\zeta_{16}^r}{2^6} \frac{\left(\beta^{\tau^2} (\beta-2) \right)^3}{\beta \beta^{\tau^2 \sigma}},\eqno{(10.11)}$$ \smallskip

\noindent for some {\it odd} integer $r$, since the $8$-th power of $\zeta_{16}^r$ must give $-1$.  On the other hand, we have from (10.6b) and (10.8) that $\displaystyle \beta \beta^{\tau^2 \sigma}=2^3/f_w(Q_2)^6$.  Putting this into (10.11) yields
$$\mathfrak{f}_2(w)=\frac{\zeta_{48}^{r'}}{2^3} \beta^{\tau^2} (\beta-2) f_w(Q_2)^2,$$\smallskip

\noindent where $r'$ is an odd integer, and therefore
$$f_w(Q_1)=(-1)^{(-d-1)/8} \zeta_{48}^{r'-v} \times \frac{1}{2^3} \beta^{\tau^2} (\beta-2) f_w(Q_2)^2.$$
If $d \not \equiv 0$ (mod $3$), this implies that $r' \equiv v$ (mod $24$), since the cube roots of unity are not contained in $\Omega_f$, and therefore
$$f_w(Q_1)=(-1)^{(-d-1)/8+(r'-v)/24} \frac{1}{2^3} \beta^{\tau^2} (\beta-2) f_w(Q_2)^2.\eqno{(10.12a)}$$
If $3 \mid d$, then on cubing we obtain instead that $r' \equiv v$ (mod $8$) and
$$f_w(Q_1)^3=(-1)^{(-d-1)/8+(r'-v)/8} \left(\frac{1}{2^3} \beta^{\tau^2} (\beta-2) \right)^3 f_w(Q_2)^6.\eqno{(10.12b)}$$
These equations express relations between the units $f_w(Q_1), f_w(Q_2)$, and $\xi^{\tau^2} (\xi-1)/2$.  Note that the last unit is one-half the $Y$-coordinate of the point $(-\beta, \xi^{\tau^2} (\xi-1))$ on the curve $Y^4=X(X^2+4)$.  (See (8.5).)    \medskip

There are similar formulas for the quadratic forms
$$Q_c(x,y)=cx^2-vxy+\left(\frac{v^2+d}{4c}\right)y^2,$$
$$Q_{2c}(x,y)=2cx^2-vxy+\left(\frac{v^2+d}{8c}\right)y^2,$$
where $(c,6d)=1$, $v$ is odd, and $16 c \mid (v^2+d)$.  These are the forms corresponding to the ideals $\mathfrak{c}=(c, cw_1)$ and $\mathfrak{c} \wp_{2,-d}=(2c,cw_1)$, where $w_1=w/c=(v+\sqrt{-d})/2c$.  For these forms the Yui-Zagier formulas give

$$f_w(Q_c)=(-1)^{(-d-1)/8} \zeta^{-v(c+(c^2-1)(v^2+d)/(4c))}\mathfrak{f}_2(w/c)$$
$$=(-1)^{(-d-1)/8}\zeta^{-vc}\mathfrak{f}_2(w/c), \ \ (c,6d)=1;\eqno{(10.13)}$$
$$f_w(Q_{2c})=\zeta^{-v(2c-(v^2+d)/(8c)-(v^2+d)^2/(32c))} \mathfrak{f}(w/(2c))$$
$$=\zeta^{-2vc(c^2-(v^2+d)/16-(v^2+d)^2/64)} \mathfrak{f}(w/(2c))$$
$$=\zeta^{-2vc(1-(v^2+d)/16-(v^2+d)^2/64)} \mathfrak{f}(w/(2c)),\eqno{(10.14)}$$

\noindent where the next to last equality follows from multiplying inside and outside the parenthesis in the exponent by $c$, which is valid since $c^2 \equiv 1$ (mod 24). \medskip

Further, let $\tau_c=\left(\frac{\Omega_f/K}{R_K \mathfrak{c}} \right)$ (note that $\tilde \mathfrak{c}=R_K \mathfrak{c}$ is the ideal of $R_K$ corresponding to the ideal $\mathfrak{c}$ of $\textsf{R}_{-d}$), so that
$$j(w)^{\tau_c^{-1}}=j(\mathfrak{c})=j(w_1)$$
and
$$j(w/2)^{\tau_c^{-1}}=j(\wp_{2,-d} \mathfrak{c})=j(w_1/2).$$
Then, defining
$$\alpha_c=\alpha_1^{\tau_c^{-1}}, \ \ \beta_c=\beta^{\tau_c^{-1}},$$
we have all the same formulas for $\alpha_c, \beta_c$ from Sections 6-8 that we have for $\alpha_1,\beta$.  In place of the automorphism $\sigma$ in those formulas we need to take its conjugate $\sigma_c=\tau_c \sigma \tau_c^{-1}$, since
$$\frac{2(\beta_c+2)}{\beta_c-2}=\beta^{\sigma \tau_c^{-1}}=\beta_c^{\tau_c \sigma \tau_c^{-1}}.$$
Furthermore, the formulas of this section all hold when $w$ is replaced by $w_1=w/c$. \medskip

Now we apply the Reciprocity Law of Shimura [32, p. 123 and p. 72] to the algebraic numbers $\mathfrak{f}(w/2)$ and $\mathfrak{f}_2(w)$ in (10.7) and (10.10).  For this we note that the coefficients in the $q$-expansions of $\mathfrak{f}(\tau)$ and $\mathfrak{f}_2(\tau)$ are in $\mathbb{Q}$ and $\mathbb{Q}(\sqrt{2})$, respectively, and these are modular functions for $\Gamma(48)$ (see [36] and [32, p. 148]).  Taking the matrix $C$ in the Reciprocity Law to be
$$C=\left(\begin{array}{cc}1 & 0 \\0 & c\end{array}\right), \ \ \textrm{with} \ \ cC^{-1}=\left(\begin{array}{cc}c & 0 \\0 & 1\end{array}\right)=-TCT,$$
and $T=\left(\begin{array}{cc}0 & 1 \\-1 & 0\end{array}\right)$, we have that $C(w, 1)^t$ is a basis of the ideal $\mathfrak{c}$, and
$$\mathfrak{f} \circ cC^{-1}=\mathfrak{f}, \ \ \mathfrak{f}_2 \circ cC^{-1}= \mathfrak{f}_2,$$
in Schertz's notation [32, p. 72].  (The second formula holds because $\mathfrak{f}_2(-1/\tau)=\mathfrak{f}_1(\tau)$ and the latter has a $q$-series with rational coefficients.)  Now the Frobenius automorphism for the conjugate ideal $\bar \mathfrak{c}$ satisfies $(\Omega_f/K,R_K \bar \mathfrak{c}) =\tau_c^{-1}$.  If $L$ is an abelian extension of $K$ containing 
the numbers $\zeta=\zeta_{48}$, $\mathfrak{f}(w/2)$, and $\mathfrak{f}_2(w)$ ($L=\Sigma_{48f}(\zeta_{48})$ would suffice by [32, Theorem 5.2.1]), the Frobenius automorphism $(L/K, R_K \bar \mathfrak{c})$ extends the automorphism $\tau_c^{-1}$ on $L/K$.  Denoting $(L/K, R_K \bar \mathfrak{c})$ by $\tau_c^{-1}$, the Reciprocity Law implies
$$\mathfrak{f}(w/2)^{\tau_c^{-1}}=\mathfrak{f}(w/(2c))=\mathfrak{f}(w_1/2), \ \ \mathfrak{f}_2(w)^{\tau_c^{-1}}= \mathfrak{f}_2(w/c)= \mathfrak{f}_2(w_1).$$
Furthermore, $\zeta^{\tau_c^{-1}}=\zeta^{Norm(\bar \mathfrak{c})}=\zeta^c$ and equations (10.7), (10.10), (10.13), and (10.14) give
$$f_w(Q_2)^{\tau_c^{-1}}=f_w(Q_{2c}), \ \ f_w(Q_1)^{\tau_c^{-1}}=f_w(Q_{c}).\eqno{(10.15)}$$

Now it is easy to see that every ideal class in $\textsf{R}_{-d}$ contains an ideal $\mathfrak{c}$ with $(Norm(\mathfrak{c}),6d)=1$, given by $\mathfrak{c}=(c,(v+\sqrt{-d})/2)$, with $16c \mid (v^2+d)$.  Furthermore, $(1,(v+\sqrt{-d})/2) = (1)$, so that $f_w(Q_1)$ is the singular modulus for the principal class, as above.  From (10.15) it follows that $f_w(Q_1)$ is conjugate to the Weber singular modulus of every other ideal class.  \medskip

This verifies the conjecture of [36, p. 1648] that the numbers $f_w(\mathcal{A})=f_w(Q)$ are algebraic conjugates of each other, as $\mathcal{A}$ varies over ideal classes (or equivalence classes of quadratic forms), in the case that $3$ does not divide $d$.  By Proposition 10.1b) and the above computations we get a similar statement when $3 \mid d$.  (Note that $f_w(Q_2)^e$, for $e=1$ or $3$, has $h(-d)$ conjugates since $\mathbb{Q}(f_w(Q_2)^e)=\mathbb{Q}(j(w/2))$.)  \bigskip

\noindent {\bf Theorem 10.3.} a) If $(d,3)=1$ and $\mathcal{O}=\textsf{R}_{-d}$, then the polynomial
$$W_{-d}(X)=\prod_{\mathcal{A} \in Pic(\mathcal{O})} {(X-f_w(\mathcal{A}))}$$
is the minimal polynomial of $f_w(Q_1)$ over $\mathbb{Q}$. \medskip

b) If $3 \mid d$, then
$$\tilde W_{-d}(X)=\prod_{\mathcal{A} \in Pic(\mathcal{O})} {(X-f_w(\mathcal{A})^3)}$$
is the minimal polynomial of $f_w(Q_1)^3$ over $\mathbb{Q}$.   \bigskip

If we now apply the automorphism $\tau_c^{-1}$ to (10.12a) and (10.12b), we obtain: \bigskip

\noindent {\bf Theorem 10.4.} The Weber singular moduli for the quadratic forms $Q_c=[c,-v,\frac{v^2+d}{4c}]$ and $Q_{2c}=[2c,-v,\frac{v^2+d}{8c}]$ corresponding to the ideals $\mathfrak{c}$ and $\wp_{2,-d} \mathfrak{c}$, with $c=Norm(\mathfrak{c})$ and $(c,6d)=1$, are related by the formulas
$$f_w(Q_c)= \frac{\epsilon}{2^3} \beta_c^{\tau^2} (\beta_c-2) f_w(Q_{2c})^2, \ \ \textrm{if} \ (d,3)=1,$$
$$f_w(Q_c)^3= \frac{\epsilon}{2^9} (\beta_c^{\tau^2} (\beta_c-2))^3 f_w(Q_{2c})^6, \ \ \textrm{if} \ 3 \mid d,$$
where $\epsilon = (-1)^{(-d-1)/8+e(r'-v)/24}$ ($e=1$ or $3$) is independent of $c$ and $\beta_c=\beta^{\tau_c^{-1}}$, with $\tau_c=(\Omega_f/K, R_K \mathfrak{c})$.  \bigskip

As a corollary of this theorem, we have: \bigskip

\noindent {\bf Theorem 10.5.} Assume that $(d,3)=1$.  If $\mathcal{A}$ represents the ideal class containing $\mathfrak{c}$ in the order $\textsf{R}_{-d}$, and $\mathcal{T}$ is the class containing the ideal $\wp_{2,-d}=\wp_2 \cap \textsf{R}_{-d}$, then
$$2\frac{f_w(\mathcal{A})}{f_w(\mathcal{TA})^2}=\frac{\epsilon}{4}\beta_c^{\tau^2}(\beta_c-2)$$
is the $Y$-coordinate of the point $\Phi=(-\beta_c,\frac{\epsilon}{4}\beta_c^{\tau^2}(\beta_c-2))$ on the curve $Y^4=X(X^2+4)$. \bigskip

The formulas in (10.15) show that the Weber singular moduli $f_w(Q)$ transform under the automorphisms of $Gal(\Omega_f/K)$ in the same way that the corresponding $j$-invariants transform.  Therefore, $j(w)^{\tau^{-1}}=j(w/2)$ and $j(w)^{\tau^{-2}}=j(w/4)$ imply that
$$f_w(Q_1)^{\tau^{-1}}=f_w(Q_2), \ \ f_w(Q_1)^{\tau^{-1} \tau_c^{-1}} = f_w(Q_{2c}),$$
$$f_w(Q_1)^{\tau^{-2}}=f_w(Q_4), \ \ Q_4(x,y)=4x^2-vxy+\frac{v^2+d}{16}.$$
\medskip 

We are now ready to verify the formula (10.4).  \bigskip

\noindent {\bf Theorem 10.6.} Equations (10.3)-(10.4) hold.  Furthermore, $\displaystyle \zeta_8^{(j-1)/2} \frac{\eta(w/4)}{\eta(w)} \in \Omega_f$, where
$$j \equiv  \cases{ \frac{3d+19}{8}=3\left(\frac{d+1}{8} \right)+2  \hspace{.05 in} (\textrm{mod} \hspace{.05 in} 8), &if $v=3$ \hspace{.01 in} \textrm{and} \hspace{.01 in} $d \equiv 7$ (mod $16$),\cr
	\frac{d+25}{8} = \frac{d+1}{8} +3 \hspace{.05 in} (\textrm{mod} \hspace{.05 in} 8), &if $v=1$ \hspace{.01 in} \textrm{and} \hspace{.01 in} $d \equiv 15$ (mod $16$).}
	 \eqno{(10.16)}$$ \medskip

\noindent {\it Proof.} We first prove that (10.4) holds. Note that $\beta$ is determined among the fourth roots of $\beta^4$ by the condition that $\displaystyle \beta^\sigma=\frac{2(\beta+2)}{\beta-2}=\bar \beta$ is the complex conjugate of $\beta$.  This is because $j(w)$ is fixed by $\sigma$: since $j(w)=j(\mathcal{O})$ is the $j$-invariant corresponding to the principal class $\mathcal{O}=\textsf{R}_{-d}$, it is real, and because $Fix(\sigma)=\mathbb{Q}(j(w))$ it follows that $\sigma$ is complex conjugation on $\Omega_f$.  This implies that the equation
$$\beta \bar \beta=2(\beta+\bar \beta)+4$$
holds.  Now this condition holds for $\beta$, and we claim that it cannot also hold for $-\beta, i\beta$, or $-i \beta$.  If it also held for $-\beta$ in place of $\beta$ then we would have that $\beta+\bar \beta=0$, which is impossible, since $\bar \beta$ is a $\mathbb{Q}$-conjugate of $\beta$, but $-\beta$ is not (see Proposition 8.1).  If it held for $i \beta$ in place of $\beta$, then we would have $2(\beta+\bar \beta)=2i(\beta-\bar \beta)$, so $i=(\beta + \bar \beta)/(\beta-\bar \beta) \in \Omega_f$, which is impossible.  The same argument applies to $-i \beta$.  Therefore, the factor $i^b$ in (10.3) is determined by the condition that
$$\frac{i^b \mathfrak{f}_1(w/2)^2+\mathfrak{f}(w/2)^2}{i^b \mathfrak{f}_1(w/2)^2-\mathfrak{f}(w/2)^2}=\frac{i^{-b} \mathfrak{f}_1(- \bar w/2)^2}{\mathfrak{f}(-\bar w/2)^2},\eqno{(10.17)}$$
where we have used the fact that $\mathfrak{f}$ and $\mathfrak{f}_1$ have $q$-expansions with rational coefficients. \medskip

Now we let $q_1=q^{1/2}=e^{\pi i w/2}=e^{\pi i(v+i\sqrt{d})/4}$, so that $\bar q_1=e^{-\pi i v/2} q_1 = i^{-v} q_1$, and write
$$q^{\frac{1}{48}}\mathfrak{f}_1(w/2)=A(q_1)=\prod_{n=1}^\infty{(1-q_1^{2n-1})}$$
and
$$q^{\frac{1}{48}}\mathfrak{f}(w/2)=A(-q_1)=\prod_{n=1}^\infty{(1+q_1^{2n-1})}.$$
With this notation, (10.17) is equivalent to the relation
$$(i^b A(q_1)^2+A(-q_1)^2)A(-i^{-v} q_1)^2= (A(q_1)^2-i^{-b} A(-q_1)^2)A(i^{-v} q_1)^2.\eqno{(10.18)}$$
Now, for $v=3$, say, and an infinite sequence of integers $d \equiv 7$ (mod $16$), we know that (10.18) holds for a fixed choice of the factor $i^b$.  For this sequence of $d$'s we have $q_1 \rightarrow 0$ as $d \rightarrow \infty$.  Since $A(q_1)$ is holomorphic in $q_1$ for $|q_1| < 1$, we thus have an identity in $q_1$:
$$(i^b A(q_1)^2+A(-q_1)^2)A(-i q_1)^2= (A(q_1)^2-i^{-b} A(-q_1)^2)A(i q_1)^2.$$
Now consider the first two terms on both sides of this equation, with $\lambda=i^b$:
$$1+\lambda+ (-2\lambda+2+2i(\lambda+1))q_1+ O(q_1^2) =1-\frac{1}{\lambda}+(-2-\frac{2}{\lambda}-2i(1-\frac{1}{\lambda}))q_1+O(q_1^2).$$
Equating coefficients of $q_1$ gives that $\lambda=i$, and therefore $i^b=i^{-v}$ when $v=3$.  This proves the identity
$$(i A(q_1)^2+A(-q_1)^2)A(-i q_1)^2= (A(q_1)^2+i A(-q_1)^2)A(i q_1)^2,$$
which is equivalent to the assertion that the left side of this equation is an even function of $q_1$.  Taking the complex conjugate of this equation now shows that when $v=1$, (10.18) holds with $i^b=-i=i^{-v}$.  This proves (10.4).  \medskip

Taking $v=1$ or $3$ and $b=-v$ in (10.9) implies 
$$a \equiv v \left(2-\frac{v^2+d}{16} \right) \hspace{.05 in} (\textrm{mod} \hspace{.05 in} 4),$$
or
$$a \equiv  \cases{ \frac{-3d+5}{\ 16} \hspace{.05 in} (\textrm{mod} \hspace{.05 in} 4), &if $v=3$ \hspace{.01 in} \textrm{and} \hspace{.01 in} $d \equiv 7$ (mod $16$),\cr
	\frac{-d+31}{16} \hspace{.05 in} (\textrm{mod} \hspace{.05 in} 4), &if $v=1$ \hspace{.01 in} \textrm{and} \hspace{.01 in} $d \equiv 15$ (mod $16$).}
	 \eqno{(10.19)}$$ \smallskip
	 
\noindent Now

$$\zeta_8^j \alpha=2\frac{\beta/2}{ \beta/\zeta_8^j \alpha}=2i^b\frac{\mathfrak{f}_1(w/2)^2}{\mathfrak{f}(w/2)^2} \times i^{-a}\frac{\mathfrak{f}(w/2)^2}{\mathfrak{f}_2(w/2)^2}=2i^{b-a}\frac{\mathfrak{f}_1(w/2)^2}{\mathfrak{f}_2(w/2)^2}.\eqno{(10.20)}$$ \smallskip

\noindent Expressing the last quotient in this equation in terms of the $\eta$-function gives

$$\zeta_8^j \alpha=2i^{b-a} \frac{\eta(w/4)^2}{\eta(w/2)^2} \times \frac{\eta(w/2)^2}{2\eta(w)^2}=i^{b-a}\frac{\eta(w/4)^2}{\eta(w)^2}.$$

\noindent From the definition of $\alpha$ in (2.7) we conclude that $\zeta_8^{j-1} = i^{b-a}$, whence we have that $j \equiv 2(b-a) + 1$ (mod 8).  Then (10.4) and (10.19) yield (10.16).
$\square$ \bigskip

The above proof shows that identity (10.18) holds with $i^b=i^{-v}$ in all cases, with $w=(v+\sqrt{-d})/2$, whether or not $v$ is restricted to be $1$ or $3$.  Thus we have: \bigskip

\noindent {\bf Theorem 10.7.} If $w=(v+\sqrt{-d})/2$, with $16c \mid (v^2+d)$ and $(c,6d)=1$, then
$$\frac{\beta_c}{2}=\frac{\beta^{\tau_c^{-1}}}{2}=i^{-vc}  \frac{\mathfrak{f}_1(w/(2c))^2}{\mathfrak{f}(w/(2c))^2}.$$
The numbers $\beta_c$ are class invariants.   Namely, if $\mathfrak{c}_1=(c_1,(v_1+\sqrt{-d})/2)$ and $\mathfrak{c}_2=(c_2,(v_2+\sqrt{-d})/2)$ are ideals of $\textsf{R}_{-d}$ with norms $c_1$ and $c_2$ satisfying
$$16c_i \mid (v_i^2+d), \ \ (c_1c_2,6d)=1, \ \ v_1 \equiv v_2 \ (\textrm{mod} \ 8),$$ then $\beta_{c_1}=\beta_{c_2}$ if and only if $\mathfrak{c}_1 \sim \mathfrak{c}_2$ in the ring $\textsf{R}_{-d}$. \medskip

\noindent {\it Proof.}  The formula for $\beta_c$ follows from the above remarks.  For the proof that the numbers $\beta_c$ are class invariants, let $\mathfrak{c}_1$ and $\mathfrak{c}_2$ be as in the statement of the theorem.  Note first that the condition $v_1 \equiv v_2$ (mod $8$) guarantees that the ideal $\wp_2=(2,(v_i+\sqrt{-d})/2)$ (for $i=1, 2$) is consistently determined.  Actually, it would be enough to assume that $v_1 \equiv v_2$ (mod $4$), but then
$$v_1^2 \equiv (v_2+4k)^2 \equiv v_2^2+8kv_2 \equiv v_1^2+8kv_2 \ (\textrm{mod} \ 16)$$
implies the stronger condition that $8 \mid (v_1-v_2)$.  With this condition we have
that
$$\frac{\beta}{2}=i^{-v_1} \frac{\mathfrak{f}_1((v_1+\sqrt{-d})/2)^2}{\mathfrak{f}((v_1+\sqrt{-d})/2)^2}=i^{-v_1} \frac{\mathfrak{f}_1((v_2+\sqrt{-d})/2+4k)^2}{\mathfrak{f}((v_2+\sqrt{-d})/2+4k)^2}.$$
By the transformation formulas for the Schl\"afli functions (see [36, p. 1647]), this becomes
$$\frac{\beta}{2}=i^{-v_2} \frac{\zeta^{-8k}\mathfrak{f}_1((v_2+\sqrt{-d})/2)^2}{\zeta^{-8k}\mathfrak{f}((v_2+\sqrt{-d})/2)^2}=i^{-v_2} \frac{\mathfrak{f}_1((v_2+\sqrt{-d})/2)^2}{\mathfrak{f}((v_2+\sqrt{-d})/2)^2},$$
so that the value of $\beta$ does not depend on the choice of $v_i$. If $\beta_{c_1}=\beta_{c_2}$, then by (8.0) and $j(w)^{\tau_{c_i}^{-1}}=j(\mathfrak{c}_i)$ we have $j(\mathfrak{c}_1)=j(\mathfrak{c}_2)$, which implies $\mathfrak{c}_1 \sim \mathfrak{c}_2$ in $\textsf{R}_{-d}$.  Conversely, assume that $\mathfrak{c}_1 \sim \mathfrak{c}_2$.  Then $\tau_{c_1}=\tau_{c_2}$ and $\beta_{c_1}=\beta^{\tau_{c_1}^{-1}}=\beta^{\tau_{c_2}^{-1}}=\beta_{c_2}$.  This proves the theorem.  $\square$  \bigskip

The proof of Theorem 10.6 also allows us to give an explicit formula for the unit $\gamma$ and its square-root. \bigskip

\noindent {\bf Proposition 10.8.}  If $w=(v+\sqrt{-d})/2$, with $16 \mid (v^2 + d)$, then
$$\sqrt{\gamma}=\pm \frac{1}{\sqrt{2}} \frac{\mathfrak{f}_1(w/2) \mathfrak{f}_1(-\bar w/2)}{\mathfrak{f}(w/2) \mathfrak{f}(-\bar w/2)}.$$
Thus, this expression lies in the field $\mathbb{Q}(j(w))$. \medskip

\noindent {\it Proof.} This follows from
$$\gamma=\frac{1}{2} \frac{\beta}{2} \frac{\beta+2}{\beta-2}=\frac{1}{2} i^{-v} \frac{\mathfrak{f}_1(w/2)^2}{\mathfrak{f}(w/2)^2} i^v \frac{\mathfrak{f}_1(- \bar w/2)^2}{\mathfrak{f}(- \bar w/2)^2},$$
by (10.3) and (10.17).

\section{Generators of $\Omega_f/K$ when $3 \mid d$.}

When $(d,3)=1$, the Weber singular moduli $f_w(Q)$ from [36] are generators of small height of the ring class fields $\Omega_f$.  However, when $3 \mid d$, the generators given by Yui and Zagier are not quite as good, in terms of the size of their discriminants, as is illustrated by the examples $d= 87, 159, 231$ given in [36, p. 1660].  These generators are essentially the numbers $f_w(Q)^3$ (when $3 \mid h(-d)$) and $f_w(Q)/\varepsilon^{\pm 1/3}$ (when $h(-d) \equiv \pm 1$ (mod $3$)), where $\varepsilon$ is a fundamental unit of the field $\mathbb{Q}(\sqrt{d/3})$.  However, the latter generators remained conjectural in [36].  \medskip

In the case that $3 \mid d$, Theorem 10.4 gives us a systematic way of finding generators for $\Omega_f$ of small height, independent of whether $3 \mid h(-d)$ or not.  We first multiply both sides of (10.12b) by $f_w(Q_2)^3$.  This gives
$$f_w(Q_1)^3f_w(Q_2)^3= \frac{\epsilon}{2^9} (\beta^{\tau^2} (\beta-2))^3 f_w(Q_2)^9.$$
The right side of this equation is a cube in $\Omega_f$, so taking cube roots yields
$$f_w(Q_1)f_w(Q_2)= \frac{\epsilon \omega^l}{2^3} (\beta^{\tau^2} (\beta-2)) f_w(Q_2)^3,\eqno{(11.1)}$$
for some $l$, where $\omega=e^{2\pi i/3}$.  Thus, $f_w(Q_1)f_w(Q_2) \in \Omega_f$, and by (10.15) we have
$$f_w(Q_c)f_w(Q_{2c})= \frac{\epsilon \omega^{lc}}{2^3} (\beta_c^{\tau^2} (\beta_c-2)) f_w(Q_{2c})^3.$$
It is easy to see that $f_w(Q_1)f_w(Q_2)=f_w(Q_1)^{1+\tau^{-1}}$ has degree at most $h(-d)$ over $\mathbb{Q}$, since it is fixed by the involution $\sigma \tau^{-1}$:
$$(f_w(Q_1)f_w(Q_2))^{\sigma \tau^{-1}}=f_w(Q_1)^{(1+\tau^{-1})\sigma \tau^{-1}}= f_w(Q_1)^{\sigma (1+\tau^{-1})}=f_w(Q_1)^{1+\tau^{-1}}.$$ \smallskip
In addition, $\omega=e^{2 \pi i/3}$ is also fixed by $\sigma \tau^{-1}$, since $\omega^{\sigma \tau^{-1}}=\bar \omega^{\tau^{-1}}=\bar \omega^2 = \omega$.  \medskip

More generally, a similar computation shows that each of the numbers $f_w(Q_1) f_w(Q_c)$ is fixed by the involution $\sigma \tau_c^{-1}$, so all the pairwise products of roots of $\tilde W_{-d}(x)$ have degree at most $h(-d)$.  Thus, the irreducible factors of the resultant
$$R_d(x)=\textrm{Res}_y \left(y^{h(-d)} \tilde W_{-d}(x/y), \tilde W_{-d}(y) \right),\eqno{(11.2)}$$
one of which is the minimal polynomial $\tilde t(x)$ of $f_w(Q_1)^3f_w(Q_2)^3$ over $\mathbb{Q}$, have degree at most $h(-d)$, and $\tilde t(x)$ has the property that $\tilde t(x^3)$ also factors into the product of three irreducibles of degree $h(-d)$, as we show in the following theorem.  \bigskip

\noindent {\bf Theorem 11.1.} If $3 \mid d$ and $d>15$, the numbers $f_w(Q_c)f_w(Q_{2c})$ lie in $\Omega_f$ and form a complete system of conjugates over $\mathbb{Q}$ (and over $K=\mathbb{Q}(\sqrt{-d})$) of
$$f_w(Q_1)f_w(Q_2)=(-1)^{(-d-1)/8} \zeta^{-v(3-(v^2+d)/8-(v^2+d)^2/32)} \mathfrak{f}(w/2) \mathfrak{f}_2(w),$$
as $\mathfrak{c}=(c,(v+\sqrt{-d})/2)$ varies over ideals in the $h(-d)$ ideal classes in $\textsf{R}_{-d}$, with $c=Norm(\mathfrak{c})$, $(c,6d)=1$, and $16c \mid (v^2+d)$.  \medskip

\noindent {\it Proof.} The formula for $f_w(Q_1)f_w(Q_2)$ is immediate from (10.7) and (10.10).  It is also clear that the numbers $f_w(Q_c)f_w(Q_{2c})$ are all conjugate to $f_w(Q_1)f_w(Q_2)$, by (10.15).  We will show that these numbers are distinct, as the ideal $\mathfrak{c}$ varies over distinct ideal classes in $\textsf{R}_{-d}$.  Assume that $\mathfrak{c}$ is an ideal with $Norm(\mathfrak{c})=c$, satisfying the conditions of the theorem, for which
$$f_w(Q_1)f_w(Q_2)=f_w(Q_c)f_w(Q_{2c})=\left(f_w(Q_1)f_w(Q_2) \right)^{\tau_c^{-1}}.\eqno{(11.3)}$$
Raising (11.1) to the 4-th power gives
$$(f_w(Q_1)f_w(Q_2))^4= \frac{\omega^l}{2^4} \left(\frac{1}{4}\beta^{\tau^2} (\beta-2) \right)^4 f_w(Q_2)^{12}$$
$$=\frac{\omega^l}{2^4} (-\beta(\beta^2+4)) \left(\frac{4\alpha_1}{\beta^2} \right)^2=\frac{-\omega^l (\beta^2+4)}{\beta^3} \alpha_1^2.$$
Squaring this equation gives
$$(f_w(Q_1)f_w(Q_2))^8=\frac{\omega^{2l} (\beta^2+4)^2}{\beta^6} \frac{-16\beta^4}{\beta^4-16}=\frac{-16 \omega^{2l} (\beta^2+4)}{\beta^2(\beta^2-4)},$$
and (11.3) yields
$$\frac{\omega^{2l} (\beta^2+4)}{\beta^2(\beta^2-4)}=\frac{\omega^{2lc} (\beta_c^2+4)}{\beta_c^2(\beta_c^2-4)}.\eqno{(11.4)}$$
Assume first that $\omega^l =1$ or $\tau_c=(\Omega_f/K, \mathfrak{c})$ fixes $\omega$.  Then
$$\frac{(\beta^2+4)}{\beta^2(\beta^2-4)}=\frac{(\beta_c^2+4)}{\beta_c^2(\beta_c^2-4)}$$
implies that 
$$(\beta^2-\beta_c^2)\left( (\beta^2+4)\beta_c^2+4(\beta^2-4) \right)=0.$$
Now $\beta^2=\beta_c^2$ implies $\beta=\beta_c$ by Proposition 8.1.  Thus, if $\beta \neq \beta_c$, then $\beta_c^2=-4(\beta^2-4)/(\beta^2+4)$.  However, comparing ideal factorizations, we have
$$(\beta_c^2)=(\beta^2)^{\tau_c^{-1}}=\wp_2^2 \wp_2'^4, \ \ \left(\frac{-4(\beta^2-4)}{\beta^2+4}\right)=\frac{\wp_2^2 \wp_2'^2 \wp_2^5 \wp_2'^2}{\wp_2^3 \wp_2'^2}=\wp_2^4 \wp_2'^2.$$
Since these two ideals are distinct, this shows that $\beta=\beta_c$ and therefore $\tau_c = 1$.  \medskip

Now if $\tau_c$ satisfies (11.3), then so does $\tau_c^2$, which certainly fixes $\omega$.  The above argument shows that $\tau_c^2=1$, so any nontrivial automorphism satisfying (11.3) must have order $2$.  If there are two nontrivial automorphisms $\tau_{c_1}$ and $\tau_{c_2}$ satisfying (11.3), then so does the product, which fixes $\omega$, so again we have $\tau_{c_1} \tau_{c_2}=1$ and $\tau_{c_1}=\tau_{c_2}$.  Hence, at most one nontrivial automorphism $\tau_c$ can satisfy (11.3), it can only have order $2$, and it does not fix $\omega$.  \medskip

Assume that there is a nontrivial automorphism $\tau_c$ which satisfies (11.4), where $\omega^{lc} \neq \omega^l$.  We apply the automorphism $\sigma$ (complex conjugation) to both sides of (11.4), using the fact that $\beta_c^\sigma=\beta^{\tau_c^{-1} \sigma}=\beta^{\sigma \tau_c^{-1}}$:
$$\frac{\omega^l (\beta-2)^2(\beta^2+4)}{\beta(\beta+2)^2}=\frac{\omega^{lc} (\beta_c-2)^2(\beta_c^2+4)}{\beta_c(\beta_c+2)^2}.\eqno{(11.5)}$$
Dividing (11.4) by (11.5) yields the equation
$$\frac{\omega^l (\beta+2)}{\beta(\beta-2)^3}=\frac{\omega^{lc} (\beta_c+2)}{\beta_c(\beta_c-2)^3}.\eqno{(11.6)}$$
In order to find the common solutions of (11.4) and (11.6) we set
$$g_1(x,y)=(x^2+4)y^2(y^2-4)-\omega^{2l(c-1)}(y^2+4)x^2(x^2-4),$$
$$g_2(x,y)=(x+2)y(y-2)^3-\omega^{l(c-1)}(y+2)x(x-2)^3,$$
where the equations $g_1(x,y)=0, g_2(x,y)=0$ have the common solution $(\beta, \beta_c)$.  Taking the resultant of these equations gives
$$\textrm{Res}_y(g_1(x,y), g_2(x,y))=8(\omega^{l(c-1)}-1)(x+2) x^4 (x-2)^4 g_3(x),$$
where $g_3(x) \in \mathbb{Q}(\omega)[x]$ with
$$\textsf{Norm}_{\mathbb{Q}(\omega)/\mathbb{Q}}(g_3(x))=(x^8+224x^4+256)(x+2)^2 (x^2-2x+8)^4 (x^4-8x^3+20x^2-16x+64)^4.$$
Since $\beta$ must be a root of this polynomial, it is clear that $d=7$ or $d=15$.  (Compare with the polynomials $B_7(x)$ and $B_{15}(x)$ in the proof of Proposition 8.1, and note that $\beta$ is not a root of the $8$-th degree factor, since otherwise $\beta^4=-112 \pm 64 \sqrt{3}$ and $\sqrt{-1} \in \Omega_f$.)  Thus, no such automorphism $\tau_c$ exists when $d>15$, and the proof is complete.  $\square$  \bigskip

\noindent {\bf Corollary.} If $3 \mid d$ and $d>15$, the field $L$ generated by $f_w(Q_1)f_w(Q_2)$ over $\mathbb{Q}$ contains the field $\mathbb{Q}(\omega)$. \medskip

\noindent {\it Proof.} By the theorem, we have $[L:\mathbb{Q}]=h(-d)$, so $L=Fix(\sigma \tau^{-1})$.  $\square$  \bigskip

If $t_d(x)$ is the minimal polynomial of $f_w(Q_1)f_w(Q_2)$ over $\mathbb{Q}$, then $t_d(x)$ factors into two irreducible polynomials of degree $h(-d)/2$ over $\mathbb{Q}(\omega)$, and has no real roots.  Table 2 below lists the polynomial $t_d(x)$ and its discriminant for $39 \le d \le 495$ and $3 \mid d$.  Note also that the polynomial $\tilde t(x^3)$ mentioned above is the product of $t_d(x)$ and the minimal polynomials over $\mathbb{Q}$ of the numbers $\omega f_w(Q_1)f_w(Q_2)$ and $\omega^2 f_w(Q_1)f_w(Q_2)$, both of which have degree $h(-d)$.  \bigskip

\begin{table}
  \centering 
  \caption{ The minimal polynomial $t_{d}(x)$ of $f_w(Q_1)f_w(Q_2)$.}\label{ }

\noindent \begin{tabular}{|c|l|c|}
\hline
  &  &  \\
$d$	&   $t_{d}(x)$  &   $\textrm{disc}(t_d(x))$ \\
\hline
% after \\ : \hline or \cline{col1-col2} \cline{col3-col4} ...
  &   &   \\
 39  &  $x^4+2x^3+2x^2+x+1$  &  $3^2 13$  \\
 63  &  $x^4-x^3+3x^2-x+1$      &   $3^3 7$    \\
 87 &   $x^6-x^5+4x^4-4x^3+5x^2-3x+1$  &   $-3^3 29^2$ \\
 111 &  $x^8+2x^7+4x^6+8x^5+9x^4+7x^3+7x^2+4x+1$  &  $3^4 11^2 37^3$  \\
 135  &  $x^6+3x^5+3x^4+2x^3-3x+1$  &  $-3^7 5^2$ \\
 159  &   $x^{10}-x^9+6x^7-x^6-6x^5+6x^4+8x^3-2x^2-3x+1$  &  $-3^9 29^2 53^4$  \\
 183  &  $x^8+x^7+x^6+4x^5+12x^4+11x^3+x^2-x+1$  &  $3^4 5^2 47^2 61^3$  \\
 207  &  $x^6+6x^5+16x^4+22x^3+16x^2+5x+1$  &  $-3^3 23^2$  \\
 231  &  $x^{12}+4x^{11}+11x^{10}+17x^9+20x^8+16x^7+2x^6-7x^5$ &  $3^{6} 7^4 11^6 17^4 31^2 71^2$  \\
         & $+2x^4+13x^3+11x^2+2x+1$  &   \\
 255  & $x^{12}+6x^{11}+11x^{10}-18x^8-6x^7+19x^6+12x^5+6x^4$  & $3^6 5^{14} 13^2 17^6$  \\
   &  $+18x^3+14x^2+3x+1$  &  \\
 279  &  $x^{12}-4x^{11}+3x^{10}+x^9+7x^8-9x^7+x^6-6x^5+6x^4$  &  $3^9 17^4 23^2 31^5 79^2$  \\
   &  $+8x^3+16x^2+6x+1$  &  \\
 303  &  $x^{10}-5x^9+7x^8+2x^7-17x^6+23x^5+8x^4-42x^3$  &  $-3^{21} 101^4$  \\
    &  $+32x^2-9x+1$  &  \\
 327  &  $x^{12}-4x^{11}+15x^{10}-30x^9+48x^8-50x^7+44x^6-28x^5$  &  $3^{18}5^2 109^5 191^2$ \\  
         &  $+18x^4-6x^3+9x^2+7x+1$  &  \\
351  &  $x^{12}-3x^{11}+6x^{10}-7x^9+15x^8+9x^7+23x^6+30x^5$  &  $3^{14} 13^5 101^2 151^2 251^2$ \\  
         &  $+36x^4+19x^3+24x^2+9x+1$  &  \\
 375  &  $x^{10}-5x^8+5x^7+30x^6+24x^5-25x^4-30x^3+15x^2$  &  $-3^5 5^{12} 7^4 41^2 239^2$  \\
        &  $+5x+1$  &  \\
 399  &  $x^{16}-4x^{15}+8x^{14}-24x^{13}+82x^{12}-182x^{11}+252x^{10}$  &  $3^{14} 7^{12} 11^2 19^8 89^2$  \\
        &  $-223x^9+114x^8+5 x^7-72x^6+61x^5+7x^4-45x^3$  &  $\cdot 199^2 227^2$  \\
        &  $+35x^2-7x+1$  &  \\
  423  &  $x^{10}-7x^9+21x^8-30x^7+22x^6-8x^5+19x^4-28x^3$  &  $-3^5 5^2 13^2 31^2 41^4 47^4$  \\
        &  $+30x^2-2x+1$  &   \\
  447  &  $x^{14}-3x^{13}+9x^{12}-28x^{11}+44x^{10}-21x^9-8x^8-29x^7$  &  $-3^{23} 5^{14} 13^2 17^4 149^6$  \\
         &  $+97x^6-63x^5+11x^4-16x^3+33x^2-3x+1$  &  \\
   471  &  $x^{16}+8x^{14}+4x^{13}-24x^{11}-30x^{10}-17x^9+81x^8+189x^7$  &  $3^{20} 11^{12} 17^4 157^7$  \\
         &  $+199 x^6+162 x^5+86 x^4+49 x^3+19 x^2-7x+1$  &  $\cdot 271^2 311^2$  \\
   495  &  $x^{16}-2x^{15}+13x^{14}-27x^{13}+62x^{12}-59x^{11}+39x^{10}$  &  $3^{12} 5^{14} 11^8 23^4 29^2$  \\
         &  $+29x^9-35x^8-20x^7+129x^6-169x^5+146x^4-57x^3$  &  $\cdot 79^2 103^2 367^2$  \\
         &  $+16x^2+8x+1$  &  \\
  &  &  \\
\hline
\end{tabular}
\end{table}

\bigskip

\section{Examples.}

\noindent If we have the class equation $H_{-d}(x)$ we may compute the polynomial $b_d(x)$, the minimal polynomial of $\beta/2$, as the unique polynomial of degree $2h(-d)$ which divides

$$G_{d}(x)=x^{16h(-d)}(1-x^4)^{h(-d)} H_{-d}\left( \frac{16(x^8-16x^4+16)^3}{x^{16}(1-x^4)} \right) \eqno{(12.1)}$$ \smallskip

\noindent and which is stabilized by the map $x \rightarrow (x+1)/(x-1)$.  (The argument of $H_{-d}(x)$ in (12.1) is obtained from (6.2) by replacing $\beta$ by $2x$.)   The minimal polynomial $A_d(x)$ of $\zeta_8^j \alpha$ can be obtained in a similar way as a factor of degree $2h(-d)$ of the polynomial

$$\tilde G_{d}(x)=(x^8+16x^4)^{h(-d)} H_{-d}\left( \frac{(x^8+16x^4+16)^3}{x^8+16x^4} \right). \eqno{(12.2)}$$ \smallskip

\bigskip

If the class equation $H_{-d}(x)$ is not available, then as in [36] or [32] the minimal polynomials $b_d(x)$ of $\beta/2$ and $A_d(x)$ of $\alpha_1=\zeta_8^j \alpha$ can be computed using the infinite product representations of the Schl\"afli functions.  For example, one may use representatives $\mathfrak{c}_i=(c_i, (v_i+\sqrt{-d})/2)$ (with $16c_i|(v_i^2+d$)) of the various ideal classes and compute the minimal polynomial of $\beta/2$ over $K$ from Theorem 10.7:

$$m(x)=\prod_{i=1}^{h(-d)}{\left( x - i^{-v_ic_i}\frac{\mathfrak{f}_1(w_i/2)^2}{\mathfrak{f}(w_i/2)^2} \right)}, \quad w_i = \frac{v_i+\sqrt{-d}}{2c_i}. \eqno{(12.3)}$$ \smallskip

\noindent It is necessary to take $v_i \equiv v_k$ (mod $8$) for $i \neq k$ in this product, in order for the ideal $(\beta_{c_i})$ to be $\wp_2 \wp_2'^2$.  Then $b_d(x)=m(x) \bar m(x)$.  The polynomial $A_d(x)$ can be found in a similar way using (10.20), or can be determined as the factor of the polynomial
$$\tilde A(x)=\textrm{Res}_y(b_d(y),( 1-y^4)x^4-16y^4),$$
which has $\alpha_1=\zeta_8^j \alpha$ from (10.20) as a root.  \medskip

In Tables 3 and 4 we list the polynomials $b_d(x)$ and $A_d(x)$ for $15 \le d \le 103$.  In each case, the polynomial $A_d(x^2)=A_1(x)A_1(-x)$ factors, verifying Theorem 8.10 that $\zeta_8^j \alpha$ is a square in the relevant class field.  In Table 5 we list the minimal polynomial $q_d(x)$ of the unit $4\alpha_1/\beta^2$.  This polynomial may be computed as the unique irreducible factor of degree $h(-d)$ of the resultant
$$\textrm{Res}_y(A_d(xy),s_d(y)), \ \ \ s_d(y)=\textrm{Res}_x(b_d(x),y-x^2).\eqno{(12.4)}$$
By Theorems 10.2 and 10.3, the minimal polynomial $W_{-d}(x)$ of the unit $f_w(Q_1)$ in (10.10) divides $q_d(x^6)$ if $(d,3)=1$, and the minimal polynomial $\tilde W_{-d}(x)$ of $f_w(Q_1)^3$ divides $q_d(x^2)$ when $3 \mid d$.  \bigskip

\noindent {\bf Remarks.} 1. The polynomial $A_d(x)$ has the property that $x^{2h(-d)}A_d(4/x)=2^{2h(-d)}A_d(x)$.  The polynomials $b_d(x)$ and $A_d(x)$ are of course normal polynomials, so they provide interesting examples of the elementary theory discussed in [15].  \medskip

\noindent 2. Since it is a relatively simple matter to compute the polynomial $b_d(x)$ using (12.3), the class equation $H_{-d}(x)$ can be computed using the resultant
$$\textrm{Res}_y(b_d(y),y^{16}(1-y^4)x-16(y^8-16y^4+16)^3).$$
Computing this resultant mod $p$ allows one to compute $H_{-d}(x)$ (mod $p$) directly from the polynomial $b_d(x)$.  We used the polynomials $b_d(x)$ in Table 3 to compute the points $Q_K$ in Section 8. \medskip

\begin{table}
  \centering 
  \caption{The minimal polynomial $b_d(x)$ of $\beta/2$.}\label{ }

\begin{tabular}{|c|l|}
\hline
  &  \\
$d$	&   {$b_d(x)$} \\
\hline
% after \\ : \hline or \cline{col1-col2} \cline{col3-col4} ...
   &   \\
 15  &  $x^4-4x^3+5x^2-2x+4$ \\
  23 &   $x^6+x^5+9x^4-13x^3+18x^2-16x+8$ \\
  31 &  $x^6+7x^5+11x^4-15x^3+16x^2-20x+8$  \\
  39 &  $x^8-6x^7+42x^6-60x^5+53x^4-54x^3+24x^2+16$ \\
  47 & $x^{10}-15x^9+74x^8-90x^7+93x^6-187x^5+160x^4-156x^3+168x^2$ \\   
       & $-48x+32$ \\
  55 & $x^8+6x^7+78x^6-84x^5+53x^4-66x^3-12x^2+24x+16$ \\
  63 & $x^8+20x^7+110x^6-100x^5+49x^4-80x^3-40x^2+40x+16$ \\
  71 & $x^{14}-11x^{13}+195x^{12}-127x^{11}+473x^{10}-593x^9+489x^8-1285x^7$ \\
   & $+1858x^6-2880x^5+3320x^4-2656x^3+1792x^2-576x+128$ \\
  79 & $x^{10}-31x^9+290x^8-186x^7+5x^6-251x^5-56x^4-60x^3+256x^2$ \\
       & $+32x+32$ \\
  87 & $x^{12}+16x^{11}+395x^{10}+398x^9-357x^8-316x^7-155x^6-1058x^5$ \\
   & $+1332x^4-704x^3+800x^2-352x+64$ \\
  95 & $x^{16}+42x^{15}+508x^{14}-260x^{13}-78x^{12}+3160x^{11}-2072x^{10}+2204x^9$ \\
   & $+53x^8-9378x^7+12004x^6-16216x^5+16208x^4-10144x^3+6016x^2$ \\
    & $-2048x+256$ \\
  103  &  $x^{10}-21x^9+732x^8-290x^7-191x^6-369x^5-502x^4+40x^3+456x^2$  \\
     &  $+144x+32$  \\
     &   \\
\hline
\end{tabular}
\end{table}

\begin{table}
  \centering 
  \caption{ The minimal polynomial $A_d(x)$ of $\alpha_1=\zeta_8^j\alpha$.}\label{ }

\begin{tabular}{|c|l|}
\hline
  &   \\
$d$	&   {$A_d(x)$} \\
\hline
% after \\ : \hline or \cline{col1-col2} \cline{col3-col4} ...
   &   \\
 15  &  $x^4+x^3-3x^2+4x+16$ \\
  23 &   $x^6-2x^5+5x^4+3x^3+20x^2-32x+64$ \\
  31 &  $x^6-x^5+10x^4-17x^3+40x^2-16x+64$  \\
  39 &  $x^8-8x^7+26x^6-33x^5+29x^4-132x^3+416x^2-512x+256$ \\
  47 & $x^{10}-7x^9+11x^8+x^7+66x^6-241x^5+264x^4+16x^3+704x^2$ \\
       & $-1792x+1024$ \\
  55 & $x^8+5x^7+17x^6-15x^5-67x^4-60x^3+272x^2+320x+256$ \\
  63 & $x^8-9x^7+58x^6-240x^5+609x^4-960x^3+928x^2-576x+256$ \\
  71 & $x^{14}+3x^{13}-2x^{12}-18x^{11}+291x^{10}-602x^9+1069x^8-573x^7+4276x^6$ \\
   &  $-9632x^5+18624x^4-4608x^3-2048x^2+12288x+16384$ \\
  79 & $x^{10}+16x^9+96x^8+271x^7+406x^6+495x^5+1624x^4+4336x^3$ \\
  & $+6144x^2+4096x+1024$ \\
  87 & $x^{12}-5x^{11}+34x^{10}+20x^9+337x^8-35x^7+989x^6-140x^5$ \\
   & $+5392x^4+1280x^3+8704x^2-5120x+4096$ \\
  95 & $x^{16}+8x^{15}+82x^{14}+369x^{13}+1357x^{12}+3847x^{11}+9905x^{10}+20781x^9$ \\
   & $+42909x^8+83124x^7+158480x^6+246208x^5+347392x^4+377856x^3$ \\
    & $+335872x^2+131072x+65536$ \\
  103  &  $x^{10}-3x^9-21x^8+120x^7+1121x^6+3267x^5+4484x^4+1920x^3$  \\
    &  $-1344x^2-768x+1024$  \\
     &   \\
\hline
\end{tabular}
\end{table}

\begin{table}
  \centering 
  \caption{ The minimal polynomial $q_d(x)$ of $4\alpha_1/\beta^2$.}\label{ }

\begin{tabular}{|c|l|}
\hline
  &   \\
$d$	&   {$q_d(x)$} \\
\hline
% after \\ : \hline or \cline{col1-col2} \cline{col3-col4} ...
   &   \\
 15  &  $x^2-3x+1$ \\
  23 &   $x^3+2x^2+5x-1$ \\
  31 &  $x^3-x^2+10x-1$  \\
  39 &  $x^4+4x^3+2x^2-17x+1$ \\
  47 & $x^5+9x^4+27x^3+33x^2 +26x-1$ \\
  55 & $x^4-7x^3+29x^2-43x+1$ \\
  63 & $x^4-x^3+18x^2-64x+1$ \\
  71 & $x^7-9x^6+34x^5-42x^4+91x^3+50x^2+93x-1$ \\
  79 & $x^5-8x^4+63x^2+134x-1$ \\
  87 & $x^6+7x^5+46x^4+44x^3+217x^2-191x+1$ \\
  95 &  $x^8-4x^7+58x^6-95x^5+401x^4-445x^3-243x^2-263x+1$ \\
  103  &  $x^5+17x^4+119x^3+340x^2+361x-1$  \\
     &   \\
\hline
\end{tabular}
\end{table}

\noindent 3. In the case $d=55$ roots $\pi$ and $\xi$ of $b_{55}(x)=0$ satisfying $\pi^4+\xi^4=1$ can be computed using the fact that $b_{55}(x)$ is invariant under $x \rightarrow (x+1)/(x-1)$:

$$\pi =\frac{-3+\sqrt{5}+3\sqrt{-11}-\sqrt{-55}}{4}+ \left( \frac{11+4\sqrt{-11}+\sqrt{-55}}{44} \right) \sqrt{-138+62\sqrt{5}},$$
and
$$\xi=\frac{-3+\sqrt{5}-3\sqrt{-11}+\sqrt{-55}}{4}+ \left( \frac{-11+4\sqrt{-11}+\sqrt{-55}}{44} \right) \sqrt{-138+62\sqrt{5}}.$$ \smallskip
\medskip

\noindent 4. For $d=63$ a solution $(\pi_3, \xi_3)$ of $\pi_3^4+\xi_3^4=1$ is

$$\pi_3=\frac{-5-\sqrt{-3}+\sqrt{-7}+\sqrt{21}}{2}+\frac{1+\sqrt{-7}}{2} \sqrt{-9+2\sqrt{21}},$$
$$\xi_3=\frac{-5+\sqrt{-3}-\sqrt{-7}+\sqrt{21}}{2}-\frac{1-\sqrt{-7}}{2} \sqrt{-9+2\sqrt{21}}.$$ \smallskip
These numbers lie in the ring class field $\Omega_3$ over $K=\mathbb{Q}(\sqrt{-7})$ with conductor $f=3$.  We note that $\pi_3$ is divisible by $\pi_1=(1-\sqrt{-7})/2$ (not by $\xi_1=(1+\sqrt{-7})/2$!) and that $\pi_3/\pi_1$ is a unit in $\Omega_3$:
$$\frac{\pi_3}{\pi_1}=\frac{-6+3\sqrt{-3}-2\sqrt{-7}+\sqrt{21}}{4}+\frac{-3+\sqrt{-7}}{4} \sqrt{-9+2\sqrt{21}}.$$ \bigskip
\medskip

We illustrate the above discussion by working out the details of the case $d=159$, an example that was also considered in [36].  We use the following ideal basis quotients:

$$w_1=\frac{1+\sqrt{-159}}{2}, \ \ w_2=\frac{1+\sqrt{-159}}{2 \cdot 5}, \ \ w_3=\frac{9+\sqrt{-159}}{2 \cdot 5}, \ \ w_4=\frac{17+\sqrt{-159}}{2 \cdot 7},$$
$$w_5=\frac{25+\sqrt{-159}}{2 \cdot 7}, \ \ w_6=\frac{25+\sqrt{-159}}{2 \cdot 49}, \ \ w_7=\frac{33+\sqrt{-159}}{2 \cdot 13}, \ \ w_8=\frac{41+\sqrt{-159}}{2 \cdot 23},$$
$$w_9=\frac{89+\sqrt{-159}}{2 \cdot 101}, \ \ w_{10}=\frac{-87+\sqrt{-159}}{2 \cdot 23};$$ \smallskip

\noindent the formula (12.3) then yields the polynomial
$$b_{159}(x)=x^{20}+128x^{19}+4661x^{18}-12230x^{17}+40058x^{16}-14824x^{15}-30958x^{14}$$
$$+105436x^{13}-148835x^{12}+147288x^{11}-30615x^{10}-86502x^9+152488x^8-309376x^7$$
$$+268880x^6-255776x^5+217664x^4-85120x^3+50432x^2-12800x+1024.$$ \smallskip

\noindent 
From this we use the resultant $\tilde A(x)$ to compute
$$A_{159}(x)=x^{20}+5x^{19}+153x^{18}+89x^{17}+1401x^{16}-947x^{15}+5229x^{14}+3580x^{13}$$
$$+82934x^{12}-7773x^{11}+88893x^{10}-31092x^9+1326944x^8+229120x^7+1338624x^6$$
$$-969728x^5+5738496x^4+1458176x^3+10027008x^2+1310720x+1048576.$$ \smallskip

\noindent Now we compute the polynomial $q_{159}(x)$ using (12.4):
$$q_{159}(x)=x^{10}+x^9+141x^8+757x^7+3397x^6+10593x^5+18265x^4$$
$$+14032x^3+2454x^2-2501x+1,$$
and then $q_{159}(x^2)=\tilde W_{-159}(x) \tilde W_{-159}(-x)$, where
$$\tilde W_{-159}(x)=x^{10}-x^9+x^8+7x^7+63x^6+121x^5+219x^4+196x^3+146x^2+47x-1,$$
as in [36, p. 1661].  Next we compute the resultant $R_{159}(x)$ in (11.2):
$$R_{159}(x)=p_1(x) p_2(x)^2 p_3(x)^2 p_4(x)^2 p_5(x)^2 p_6(x)^2,$$
where all but one of the degrees of the $p_i(x)$ are $10$, and $\textrm{deg}(p_2(x))=5$, with
$$p_2(x)=x^5-8x^4+43x^3-101x^2+91x+1.$$
The polynomials
$$p_3(x)=x^{10}+17x^9+105x^8+312x^7+584x^6+642x^5+444x^4+245x^3+181x^2-21x+1$$
and
$$p_4(x)=x^{10}+4x^9+43x^8+86x^7-2x^6-2218x^5+5285x^4-5199x^3+2195x^2+27x+1$$
have the property that $p_3(x^3)$ and $p_4(x^3)$ split into $10$-th degree factors.  Of these six factors, the one with the smallest discriminant is the minimal polynomial of $\upsilon_1=\omega f_w(Q_1)f_w(Q_2)$, which is a factor of $p_3(x^3)$:
$$t_1(x)=x^{10}-4x^9+12x^8-21x^7+29x^6-27x^5+18x^4-7x^3+x^2+1,$$
and $\textrm{disc}(t_1(x))=-3^5 11^2 53^4$.  This discriminant is smaller than the discriminant $3^819^253^5$ of the most efficient polynomial
$$x^{10}-x^9+2x^8+7x^7+x^6-15x^5-5x^4+8x^3+5x^2-5x+1 \eqno{(12.5)}$$
given in [36, p. 1661], even though $t_1(x)$ has slightly larger coefficients.  (The polynomial $t_{159}(x)$ in Table 2 has smaller discriminant {\it and} smaller height than (12.5).)  We note that the polynomial in (12.5) is a factor of $\textrm{Res}_y(A_1(xy),A_1(y))$, where $A_1(x)$ is the minimal polynomial of $\sqrt{\alpha_1}$, and so its roots are quotients of conjugates of $\sqrt{\alpha_1}$.  This resultant is divisible by eight other 10th degree polynomials of relatively small height and discriminant.\medskip

Of course $t_1(x)$ is irreducible over $K$, so that we have $\Sigma=K(\upsilon_1)$. Note that $t_1(x)$ has no real roots, and factors over $\mathbb{Q}(\sqrt{-3})$ as the product of the polynomial
$$x^5+(2\omega-1)x^4+(-3\omega+1)x^3+(4\omega+1)x^2+(-2\omega-1)x+1$$
and its conjugate (with discriminant $53^2 \omega$).  This implies that an integral basis for $L=\mathbb{Q}(\upsilon_1)$ over $\mathbb{Q}$ is the set $\{\upsilon_1^i,\omega \upsilon_1^i, 0 \le i \le 4 \}$.  In addition, $p_2(x^3)$ has the convenient quintic factor
$$t_2(x)=x^5-2x^4+x^3-2x^2+4x+1,$$
where $\textrm{disc}(t_2(x))=3^4 53^2$.  Since $t_2(x)$ is irreducible over $K(\sqrt{-3})$ we have that $\Sigma=K(\upsilon_2, \sqrt{-3})$, where $\upsilon_2$ is a root of $t_2(x)$.  If $\upsilon_2$ is the real root of $t_2(x)$, then the real subfield of $\Sigma$ is $\Sigma^+=\mathbb{Q}(\upsilon_2,\sqrt{53})$.  \medskip

For the sake of comparison we note that the minimal polynomial of the unit $\gamma$ (see Theorem 8.6) is
$$u_{159}(x)=x^{10}+59x^9+918x^8-4696x^7+8545x^6-6567x^5+2029x^4-79x^3-63x^2+5x+1,$$
and $u_{159}(x^2)=r(x) r(-x)$, where
$$r(x)=x^{10}+x^9+30x^8-48x^7-39x^6+55x^5+29x^4-9x^3-7x^2-3x-1.$$
We have $\textrm{disc}(r(x))=3^{12} 31^2 47^2 53^5 79^2$, and so the root $\sqrt{\gamma}$ (taking the positive sign in Proposition 10.8) of $r(x)$ is not as efficient a generator of $\Omega_f/K$ as is the root $\upsilon_1$ of $t_1(x)$.\medskip

In conclusion, we see that the arithmetic of the Weber singular moduli yields natural generators of the ring class fields $\Omega_f$ over $K=\mathbb{Q}(\sqrt{-d})$, with small height and discriminant, even when $3$ divides the discriminant $-d$.  A short table of the minimal polynomials of such generators for discriminants divisible by $3$ is given in Table 2.

\bigskip

\section{References.}

\noindent [1] A. Aigner, \"Uber die M\"oglichkeit von $x^4+y^4=z^4$ in quadratischen K\"orpern, Jahresbericht der deutschen Mathematiker Vereinigung 43 (1934), 226-229.  \medskip

\noindent [2] B.J. Birch, Heegner Points: The Beginnings, in {\it Heegner Points and Rankin $L$-series}, H. Darmon and S.-W. Zhang, eds., MSRI Publications 49, Cambridge University Press, 2004, pp. 1-10.  \medskip 

\noindent [3] A. Bremner, D.J. Lewis, and P. Morton, Some varieties with points only in a field extension, Archiv Math. 43 (1984), 344-350. \medskip

\noindent [4] D. Byeon, Heegner points on elliptic curves with a rational torsion point, J. Number Theory 132 (2012), 3029-3036.  \medskip

\noindent [5] David A.Cox, {\it Primes of the Form $x^2+ny^2$; Fermat, Class Field Theory, and Complex Multiplication}, John Wiley and Sons, 1989. \medskip

\noindent [6] J.E. Cremona, {\it Algorithms for Modular Elliptic Curves}, Cambridge University Press, 1992.  \medskip

\noindent [7] M. Deuring, Die Typen der Multiplikatorenringe elliptischer Funktionenk\" orper, Abh. Math. Sem. Hamb. 14 (1941), 197-272. \medskip

\noindent [8] M. Deuring, Teilbarkeitseigenschaften der singul\"aren Moduln der elliptischen Funktionen und die Diskriminante der Klassengleichung, Commentarii math. Helvetici 19 (1946), 74-82.  \medskip

\noindent [9] M. Deuring, Die Struktur der elliptischen Funktionenk\"orper und die Klassenk\"orper der imagin\"aren quadratischen Zahlk\"orper, Math. Annalen 124 (1952), 393-426. \medskip

\noindent [10] M. Deuring, Die Klassenk\"orper der komplexen Multiplikation, Enzyklop\"adie der math. Wissenschaften I2, 23 (1958), 1-60. \medskip

\noindent [11] D. R. Dorman, Global orders in definite quaternion algebras as endomorphism rings for reduced CM elliptic curves, in: {\em Theorie des nombres}, J.-M. De Koninck and C. Levesque, eds., Proceedings of the International Number Theory Conference held at Universit\'e Laval (July 5-18, 1987), de Gruyter, Berlin, 1989, 108-116. \medskip

\noindent [12] B. Fein, B. Gordon, and J.H. Smith, On the representation of $-1$ as a sum of two squares in an algebraic number field, J. Number Theory 3 (1971), 310-315.  \medskip

\noindent [13] W. Franz, Die Teilwerte der Weberschen Tau-Funktion, J. reine angew. Math. 173 (1935), 60-64. \medskip

\noindent [14] R. Fricke, {\it Lehrbuch der Algebra}, II, III, Vieweg, Braunschweig, 1928. \medskip

\noindent [15] I. Gerst and J. Brillhart, On the prime divisors of polynomials, Amer. Math. Monthly 78 (1971), 250-266. \medskip

\noindent [16] B. Gross, Heegner points on $X_0(N)$, Chapter 5 in: {\it Modular Forms}, R.A. Rankin, ed., Ellis Horwood Limited and Halsted Press, Chichester, 1984, pp. 87-105. \medskip

\noindent [17] F. Hajir and F.R. Villegas, Explicit elliptic units, I, Duke Math. J. 90 (1997), 495-521. \medskip

\noindent [18] H. Hasse, Darstellbarkeit von Zahlen durch quadratische Formen in einem beliebigen algebraischen Zahlk\"orper, J. reine angew. Math. 153 (1924), 113-130. \medskip

\noindent [19] H. Hasse, Neue Begr\"undung der komplexen Multiplikation. I. Einordnung in die allgemeine Klassenk\"orpertheorie, J. reine angew. Math. 157 (1927), 115-139; paper 33 in {\it Mathematische Abhandlungen}, Bd. 2, Walter de Gruyter, Berlin, 1975, pp. 3-27. \medskip

\noindent [20] H. Hasse, {\it Vorlesungen \"uber Klassenk\"orpertheorie}, Physica-Verlag, W\"urzburg, 1933. \medskip

\noindent [21] F. Jarvis and P. Meekin, The Fermat equation over $\mathbb{Q}(\sqrt{2})$, J. Number Theory 109 (2004), 182-196. \medskip

\noindent [22] M. Klassen and P. Tzermias, Algebraic points of low degree on the Fermat quintic, Acta Arithmetica 82 (1997), 393-401. \medskip

\noindent [23] E. Landau, {\it Handbuch der Lehre von der Verteilung der Primzahlen}, vol. II, B.G. Teubner, Leipzig, 1909. \medskip

\noindent [24] P.S. Landweber, Supersingular curves and congruences for Legendre polynomials, in: P.S. Landweber (ed.), Elliptic Curves and Modular Forms in Topology, in: Lecture Notes in Math., vol. 1326, Springer, Berlin, 1988, 69-93. \medskip

\noindent [25] K. Lauter and B. Viray, On singular moduli for arbitrary discriminants, arXiv: 1206.6942v1, preprint, 2012.  \medskip

\noindent [26] Rodney Lynch, Arithmetic on Normal Forms of Elliptic Curves, Ph.D. Thesis, Indiana University -- Purdue University at Indianapolis (IUPUI), in preparation. \medskip

\noindent [27] P. Morton, Explicit identities for invariants of elliptic curves, J. Number Theory 120 (2006), 234-271. \medskip

\noindent [28] P. Morton, The cubic Fermat equation and complex multiplication on the Deuring normal form, Ramanujan J. of Math. 25 (2011), 247-275. \medskip

\noindent [29] I. Reiner, {\it Maximal Orders}, Academic Press, London, 1975.  \medskip

\noindent [30] R. Schertz, Die singul\"aren Werte der Weberschen Funktionen $\mathfrak{f}$, $\mathfrak{f}_1$, $\mathfrak{f}_2$, $\gamma_2$, $\gamma_3$, J. reine angew. Math. 286/287 (1976), 46-74. \medskip

\noindent [31] R. Schertz, Weber's class invariants revisited, J. de Th\'eorie de Nombres de Bordeaux 14 (2002), 325-343. \medskip

\noindent [32] R. Schertz, {\it Complex Multiplication}, New Mathematical Monographs, vol. 15, Cambridge University Press, 2010. \medskip

\noindent [33] J.H. Silverman, {\it Advanced Topics in the Arithmetic of Elliptic Curves}, in: Graduate Texts in Mathematics, vol. 151, Springer, New York, 1994. \medskip

\noindent [34] J.H. Silverman, {\it The Arithmetic of Elliptic Curves}, 2nd edition, in: Graduate Texts in Mathematics, vol. 106, Springer, New York, 2009. \medskip

\noindent [35] H. Weber, {\it Lehrbuch der Algebra}, vol. III, Chelsea Publishing Co., New York, reprint of 1908 edition. \medskip

\noindent [36] N. Yui and D. Zagier, On the singular values of the Weber modular functions, Math. Comp. 66 (1997), 1645-1662.  \bigskip

\noindent Indiana University - Purdue University at Columbus (IUPUC) \smallskip

\noindent 4601 Central Ave., Columbus, Indiana, 47203-1769 \smallskip

\noindent {\it e-mail}: rolynch@iupuc.edu \bigskip

\noindent Department of Mathematical Sciences \smallskip

\noindent Indiana University - Purdue University at Indianapolis (IUPUI) \smallskip

\noindent 402 N. Blackford St., LD 270, Indianapolis, Indiana, 46202 \smallskip

\noindent {\it e-mail}: pmorton@math.iupui.edu

\end{document}